\documentclass[11pt]{article}
\usepackage{amssymb,latexsym,amsmath,psfig}
\usepackage{graphicx}

\begin{document}

\newcommand{\bfi}{\bfseries\itshape}

\makeatletter

\@addtoreset{figure}{section}
\def\thefigure{\thesection.\@arabic\c@figure}
\def\fps@figure{h, t}
\@addtoreset{table}{bsection}
\def\thetable{\thesection.\@arabic\c@table}
\def\fps@table{h, t}
\@addtoreset{equation}{section}
\def\theequation{\thesubsection.\arabic{equation}}

\makeatother

\newtheorem{thm}{Theorem}[section]
\newtheorem{prop}[thm]{Proposition}
\newtheorem{lema}[thm]{Lemma}
\newtheorem{cor}[thm]{Corollary}
\newtheorem{defi}[thm]{Definition}
\newtheorem{rk}[thm]{Remark}
\newtheorem{exempl}{Example}[section]
\newenvironment{exemplu}{\begin{exempl}  \em}{\hfill $\surd$
\end{exempl}}

\newcommand{\comment}[1]{\par\noindent{\raggedright\texttt{#1}
\par\marginpar{\textsc{Comment}}}}
\newcommand{\todo}[1]{\vspace{5 mm}\par \noindent \marginpar{\textsc{ToDo}}\framebox{\begin{minipage}[c]{0.95 \textwidth}
\tt #1 \end{minipage}}\vspace{5 mm}\par}

\newcommand{\ea}{\mbox{{\bf a}}}
\newcommand{\eu}{\mbox{{\bf u}}}
\newcommand{\ueu}{\underline{\eu}}
\newcommand{\ueo}{\overline{u}}
\newcommand{\oeu}{\overline{\eu}}
\newcommand{\ew}{\mbox{{\bf w}}}
\newcommand{\ef}{\mbox{{\bf f}}}
\newcommand{\eF}{\mbox{{\bf F}}}
\newcommand{\eC}{\mbox{{\bf C}}}
\newcommand{\en}{\mbox{{\bf n}}}
\newcommand{\eT}{\mbox{{\bf T}}}
\newcommand{\eL}{\mbox{{\bf L}}}
\newcommand{\eR}{\mbox{{\bf R}}}
\newcommand{\eV}{\mbox{{\bf V}}}
\newcommand{\eU}{\mbox{{\bf U}}}
\newcommand{\ev}{\mbox{{\bf v}}}
\newcommand{\eve}{\mbox{{\bf e}}}
\newcommand{\uev}{\underline{\ev}}
\newcommand{\eY}{\mbox{{\bf Y}}}
\newcommand{\eK}{\mbox{{\bf K}}}
\newcommand{\eP}{\mbox{{\bf P}}}
\newcommand{\eS}{\mbox{{\bf S}}}
\newcommand{\eJ}{\mbox{{\bf J}}}
\newcommand{\eB}{\mbox{{\bf B}}}
\newcommand{\eH}{\mbox{{\bf H}}}
\newcommand{\leb}{\mathcal{ L}^{n}}
\newcommand{\eI}{\mathcal{ I}}
\newcommand{\eE}{\mathcal{ E}}
\newcommand{\hen}{\mathcal{H}^{n-1}}
\newcommand{\eBV}{\mbox{{\bf BV}}}
\newcommand{\eA}{\mbox{{\bf A}}}
\newcommand{\eSBV}{\mbox{{\bf SBV}}}
\newcommand{\eBD}{\mbox{{\bf BD}}}
\newcommand{\eSBD}{\mbox{{\bf SBD}}}
\newcommand{\ecs}{\mbox{{\bf X}}}
\newcommand{\eg}{\mbox{{\bf g}}}
\newcommand{\paromega}{\partial \Omega}
\newcommand{\gau}{\Gamma_{u}}
\newcommand{\gaf}{\Gamma_{f}}
\newcommand{\sig}{{\bf \sigma}}
\newcommand{\gac}{\Gamma_{\mbox{{\bf c}}}}
\newcommand{\deu}{\dot{\eu}}
\newcommand{\dueu}{\underline{\deu}}
\newcommand{\dev}{\dot{\ev}}
\newcommand{\duev}{\underline{\dev}}
\newcommand{\weak}{\stackrel{w}{\approx}}
\newcommand{\mild}{\stackrel{m}{\approx}}
\newcommand{\strong}{\stackrel{s}{\approx}}

\newcommand{\weakdown}{\rightharpoondown}
\newcommand{\opg}{\stackrel{\mathfrak{g}}{\cdot}}
\newcommand{\opn}{\stackrel{\mathfrak{n}}{\cdot}}
\newcommand{\tr}{\ \mbox{tr}}
\newcommand{\Ad}{\ \mbox{Ad}}
\newcommand{\ad}{\ \mbox{ad}}

\renewcommand{\contentsname}{ }

\title{Tangent bundles to sub-Riemannian groups}

\author{Marius Buliga \\
 \\
Institut Bernoulli\\
B\^{a}timent MA \\
\'Ecole Polytechnique F\'ed\'erale de Lausanne\\
CH 1015 Lausanne, Switzerland\\
{\footnotesize Marius.Buliga@epfl.ch} \\
 \\
\and
and \\
 \\
Institute of Mathematics, Romanian Academy \\
P.O. BOX 1-764, RO 70700\\
Bucure\c sti, Romania\\
{\footnotesize Marius.Buliga@imar.ro}}

\date{This version: 16.10.2003}

\maketitle


\newpage

\section{Introduction}

Classical calculus is a basic tool in analysis. We use it so often that we forget that its construction needed  considerable time and effort.

Especially in the last decade, the progresses made in the field of analysis in metric spaces make us reconsider this calculus. Along this line of thought,
all started with the definition of Pansu derivative \cite{pansu}  and its version of Rademacher theorem in Carnot groups. It is amazing that such a
basic notion can still lead to impressive results, like the rigidity of quasi-isometric embeddings.

The Gromov-Hausdorff convergence of metric spaces allows to define the notion of tangent space to a
metric space.  The tangent space to a metric space at a point is defined only up to isometry. For example, the tangent space to a n-dimensional Riemannian manifold, at a point, is $R^{n}$ with the Euclidean distance. 

For (almost) general  metric spaces Cheeger \cite{cheeger} constructed a tangent bundle. The tangent bundle constructed by Cheeger does not have as fiber the 
metric tangent space, in the case of regular sub-Riemannian manifolds.

The sub-Riemannian manifolds (and improper called sub-Riemannian geometry) form an important class of metric spaces which are not Euclidean at any scale. Basic references for sub-Riemannian geometry are Bella\"{\i}che \cite{bell} and Gromov
\cite{gromo}.

These spaces have more structure than just the metric one.  

The extra structure of sub-Riemannian manifolds
permitted to Margulis and Mostow \cite{marmos2} to construct a tangent bundle
to a sub-Riemannian manifold. Their previous paper \cite{marmos1} contains a  study of  the differentiability properties of quasi-conformal mappings between  regular sub-Riemannian manifolds. The central result is a Stepanov theorem 
in this setting, that is any quasi-conformal map from a Carnot-Carath\'eodory manifold to another is a.e. differentiable. Same result for quasi-conformal maps on Carnot groups has been first proved by Koranyi, Reimann \cite{kore}.

As a motivation of  the paper \cite{marmos2} the authors mention a comment of Deligne about the fact that in \cite{marmos1} they refer to "the tangent space" to a point, that is to a tangent bundle structure, which was not constructed in the first paper. They remedy this gap 
in the second paper \cite{marmos2}; their tangent bundle has fibers isometric with the tangent spaces.

Stepanov theorem is a hard generalisation of Rademacher theorem. 
A way to prove Stepanov theorem for sub-Riemannian spaces  is  to 
prove that  Pansu-Rademacher theorem  holds Lipschitz functions $f: A \subset M \rightarrow N$, where M, N are Carnot groups and $A$ is just measurable. This has been done for the first time  in Vodop'yanov, Ukhlov \cite{voduk}, in the case of Carnot groups. Their technique of proof  differs from the one used by  Pansu in \cite{pansu}. It  is closer to  the one used in Margulis, Mostow \cite{marmos1}, where it is proven that any quasi-conformal map from a Carnot-Carath\'eodory manifold to another is a.e. differentiable, without using Rademacher theorem.

This paper is a continuation of \cite{buliga1}, where it is argued that
sub-Riemannian geometry is in fact non-Euclidean analysis. This is more easy
to see when approaching the concept of a sub-Riemannian Lie group. As a starting
point we have a real Lie group $G$ which is connected and a vectorspace $D \subset
\mathfrak{g} = Lie(G)$ which bracket generates the Lie algebra $\mathfrak{g}$.
The vectorspace $D$ provides a left-invariant regular distribution on $G$ which is completely non-integrable. By an arbitrary choice of an Euclidean norm on $D$ we can endow
$G$ with a left invariant Carnot-Caratheodory distance. The machinery of metric spaces comes into action and tells that, as a metric space, $G$ has tangent space at any point and any such tangent space is isomorphic with the nilpotentization of $\mathfrak{g}$ with respect to $D$, denoted by $N(G)$.

$N(G)$ is a Carnot group, that is a simply connected nilpotent group  endowed with a one-parameter group of dilatations. The intrinsic calculus on $N(G)$ is based on the notion of derivative introduced by Pansu. We shall say therefore that a function is smooth if it is Pansu derivable. The same denomination will be used for functions between Carnot groups.

During the process of establishing a non-Euclidean analysis  problems begin to appear.
\begin{enumerate}
\item[a)] Right translations, the group operation  and the inverse map are not smooth.
\item[b)] If $G \not = N(G)$ then generically there is no atlas of $G$ over $N(G)$ with smooth transition of charts such that the group exponential is smooth. 
A noticeable exception is the contact case, that is the case where $N(G)$ is a Heisenberg group.
\item[c)] Mostow-Margulis (or any other known) construction of a tangent bundle
$TG$ of $G$ do not provide a group structure to $TG$.
\item[d)] there is no notion of higher order derivatives (except the horizontal
classical derivatives which have no direct intrinsic meaning).
\end{enumerate}

We would like to have a notion of sub-Riemannian Lie group such that a), b), c)
d) are solved. This is obviously a problem of good choice of "smoothness".

As a first attempt we present the formalism of uniform groups. We then concentrate on the typical example of a Lie group with left invariant distribution and give
several  notions of tangent bundle and  smoothness, all of them equivalent in the "commutative" case $D= \mathfrak{g}$. We exemplify our constructions in the case of Carnot groups and in the case of compact part of
real split groups.

{\bf Acknowledgements.}
I would first want to thank Martin Reimann and his team from the Mathematical Institute, University of Bern, not only for the invitation to participate to the Borel Seminar "Tangent spaces to metric spaces", but also for their constant interest in this type of work and for all the quality time spent with them.

My thanks go also to the members of the chairs of Geometric Analysis and Geometry,
especially to Tudor Ratiu, Marc Troyanov, Sergei Timoshin. I don't know  any other place which would be closer to perfect for doing mathematics than the Mathematics department of EPFL.

I have many reasons to thank Serguei Vodop'yanov. Among them is, it seems to me, the similar way to look to this mathematical subject. But also for his honesty,
for the opportunity to shear ideas  with a true mathematician and for his joy.

During the Borel seminar I had the occasion to meet Scott Pauls, which impressed
me with his fresh and highly rigorous ideas. I hope our discussions to continue.

And what would life be without the kind presence of Claudia?

\newpage

\tableofcontents

\newpage

\section{Overview of sub-Riemannian manifolds}

The name "sub-Riemannian manifold" is a compromise in the actual state of development of non-Euclidean analysis because it contains the word "manifold".
This denomination is therefore a bit misleading, because one does not have, up to the moment, a satisfactory notion of manifold, the Euclidean case excepted. Maybe a better name would by "model of sub-Riemannian manifold", which is rather long
and equally unclear at the moment.

\begin{defi}
A sub-Riemannian (SR) manifold is a triple $(M,H, g)$, where $M$ is a
connected manifold, $H$ is a subbundle of $TM$, named horizontal
bundle or distribution, and $g$ is a metric (Euclidean inner-product) on the
horizontal bundle.

A horizontal curve is a continuous, almost everywhere
differentiable curve, whose tangents lie in the horizontal bundle.

The length of a horizontal curve $c: [a,b] \rightarrow M$ is
$$l(c) \ = \ \int_{a}^{b} \sqrt{g_{c(t)}(\dot{c}(t), \dot{c}(t))} \mbox{ d} t$$

 The SR
manifold is called a Carnot-Carath\'eodory (CC) space if any two
points can be joined by a finite length horizontal curve.
\label{defsr}
\end{defi}

A CC space is a path metric space with the Carnot-Carath\'eodory
distance induced by the length $l$:
$$d(x,y) \ = \ \inf \left\{ l(c) \mbox{ : } c: [a,b] \rightarrow M \
, \ c(a) = x \ , \  c(b) = y \right\} $$

In the class of sub-Riemannian manifolds there is the  distinguished  subclass of regular ones. Regularity is a notion which concerns only the distribution $D$. Let us explain it.

Given the distribution $D$, a point $x \in M$ and a sufficiently small  open neighbourhood $x \in U \subset M$, one can define on $U$ a filtration of bundles
as follows. Define first the class of horizontal vectorfields on $U$:
$$\mathcal{X}^{1}(U,D) \ = \ \left\{ X \in \Gamma^{\infty}(TU) \mbox{ : }
\forall y \in U \ , \ X(u) \in D_{y} \right\}$$
Next, define inductively for all positive integers $k$:
$$ \mathcal{X}^{k+1} (U,D) \ = \ \mathcal{X}^{k}(U,D) \cup [ \mathcal{X}^{1}(U,D),
\mathcal{X}^{k}(U,D)]$$
Here $[ \cdot , \cdot ]$ denotes vectorfields bracket. We obtain therefore a filtration $\displaystyle \mathcal{X}^{k}(U,D) \subset \mathcal{X}^{k+1} (U,D)$.
Evaluate now this filtration at $x$:
$$V^{k}(x,U,D) \ = \ \left\{ X(x) \mbox{ : } X \in \mathcal{X}^{k}(U,D)\right\}$$
Because the distribution was supposed to be completely non-integrable, there are
$m(x)$, positive integer, and small enough $U$ such that $\displaystyle
V^{k}(x,U,D) = V^{k}(x,D)$ for all $k \leq m$ and
$$D_{x} V^{1}(x,D) \subset V^{2}(x,D) \subset ... \subset V^{m(x)}(x,D) = T_{x}M$$
We equally have
$$ \nu_{1}(x) = \dim V^{1}(x,D) < \nu_{2}(x) = \dim V^{2}(x,D) < ... < n = \dim M$$
Off course, $\nu_{1}(x)$ is constant over $M$, but generally $m(x)$, $\nu_{k}(x)$
may vary from a point to another.

\begin{defi}
The distribution $D$ is regular if $m(x)$, $\nu_{k}(x)$ are constant on the manifold $M$.
\label{dreg}
\end{defi}

For example a contact manifold and a Lie group with left invariant nonintegrable
distribution are examples of sub-Riemannian spaces with  regular distribution. In contrast, the Grushin plane distribution  is non regular.

The particular case that will be interesting for us is the following one.
Let $G$ be a real
 connected Lie group with Lie algebra $\mathfrak{g}$ and $D \subset \mathfrak{g}$
a vector space which generates the algebra. This means that the sequence
$$V_{1} \ = \ D \ \ , \ V^{i+1}  \ = \ V^{i} + [D,V^{i}]$$
provides a filtration of $\mathfrak{g}$:
\begin{equation}
V^{1} \ \subset \ V^{2} \ \subset \ ... \ V^{m} \ = \ \mathfrak{g}
\label{filtra}
\end{equation}
We set the distribution induced by $D$ to be
$$D_{x} \ = \ T L_{x} D$$
where $x \in G$ is arbitrary and $L_{x}: G \rightarrow G$, $L_{x}y \ = \ xy$
is the left translation by $x$. We shall use the same notation $D$ for the
induced distribution.

 This distribution is non-integrable (easy form of Chow theorem). For any (left or right invariant) metric defined on $D$
we have an associated  Carnot-Carath\'eodory distance.

\subsection{Carnot groups. Pansu derivative. Mitchell's theorem 2}

Carnot groups are particular examples of sub-Riemannian manifolds. They are especially important because they provide infinitesimal models for any sub-Riemannian manifold. Moreover, all the fundamental results of sub-Riemannian geometry are particularly easy to prove and understand in the case of Carnot groups.

\paragraph{Structure of Carnot groups. Mitchell's theorem 2}

\begin{defi}
A Carnot (or stratified nilpotent) group is a
connected simply connected group $N$  with  a distinguished vectorspace
$V_{1}$ such that the Lie algebra of the group has the
direct sum decomposition:
$$n \ = \ \sum_{i=1}^{m} V_{i} \ , \ \ V_{i+1} \ = \ [V_{1},V_{i}]$$
The number $m$ is the step of the group. The number
$$Q \ = \ \sum_{i=1}^{m} i \ dim V_{i}$$
is called the homogeneous dimension of the group.
\label{dccgroup}
\end{defi}

Because the group is nilpotent and simply connected, the
exponential mapping is a diffeomorphism. We shall identify the
group with the algebra, if is not locally otherwise stated.

The structure that we obtain is a set $N$ endowed with a Lie
bracket and a group multiplication operation given by the Baker-Campbell-Hausdorff formula.

This  formula has only a finite number of terms, because the algebra is nilpotent. This shows that the group operation
is polynomial. It is easy to see that the Lebesgue measure on the
algebra is (by identification of the algebra with the group) a bi-invariant measure.

We give further examples of such groups:

{\bf (1.)} $R^{n}$ with addition is the only commutative Carnot group. It generates the classical, or Euclidean, or commutative calculus.

{\bf (2.)} The Heisenberg group is the first non-trivial example.
This is the group $H(n) \ = \ R^{2n} \times R$ with the operation:
$$(x,\bar{x}) (y,\bar{y}) \ = \ (x+y, \bar{x} + \bar{y} +
\frac{1}{2} \omega(x,y) )$$
where $\omega$ is the standard symplectic form on $R^{2n}$.
The Lie bracket is
$$[(x,\bar{x}), (y,\bar{y})] \ = \ (0,\omega(x,y))$$
The direct sum decomposition of (the algebra of the) group is:
$$H(n) \ = \ V + Z \ , \ \ V \ = \ R^{2n} \times \left\{ 0 \right\}
\ , \ \ Z \ = \  \left\{ 0 \right\} \times R$$
$Z$ is the center of the algebra, the group has step 2 and
homogeneous dimension $2n+2$. This group generates the symplectic calculus, as shown in section 5.

{\bf (3.)} H-type groups. These are two step nilpotent Lie groups $N$
endowed with an inner product $(\cdot , \cdot)$, such that the
following {\it orthogonal} direct sum decomposition occurs:
$$N \ = \ V + Z$$
$Z$ is the center of the Lie algebra. Define now the function
$$J : Z \rightarrow End(V) \ , \ \ (J_{z} x, x') \ = \ (z, [x,x'])$$
The group $N$ is of H-type if for any $z \in Z$ we have
$$J_{z} \circ J_{z} \ = \  - \mid z \mid^{2} \  I$$
From the Baker-Campbell-Hausdorff formula we see that the group
operation is
$$(x,z) (x', z') \ = \ (x + x', z + z' + \frac{1}{2} [x,x'])$$
These groups appear naturally as the nilpotent part in the Iwasawa
decomposition of a semisimple real group of rank one. (see \cite{cdkr})

{\bf (4.)} The last example is the group of $n \times n$
upper triangular matrices, which is
nilpotent of step $n-1$.  This example is important because any Carnot group is isomorphic with a subgroup of a group of upper triangular matrices. However, 
it is not generally true that given a Carnot group $N$ with a compatible familly of dilatations, there is an injective morphism from  $N$ to a group of upper triangular matrices which commutes with dilatations (on $N$ we have the chosen dilatations and in the group of upper triangular matrices we have canonical dilatations).

Any Carnot group admits a one-parameter family of dilatations. For any
$\varepsilon > 0$, the associated dilatation is:
$$ x \ = \ \sum_{i=1}^{m} x_{i} \ \mapsto \ \delta_{\varepsilon} x \
= \ \sum_{i=1}^{m} \varepsilon^{i} x_{i}$$
Any such dilatation is a group morphism and a Lie algebra morphism.

In fact the class of Carnot groups is characterised by the existence of dilatations.

\begin{prop}
Suppose that the Lie algebra $\mathfrak{g}$ admits an one parameter group
$\varepsilon \in (0,+\infty) \mapsto \delta_{\varepsilon}$ of simultaneously
diagonalisable Lie algebra isomorphisms. Then $\mathfrak{g}$ is the algebra of
a Carnot group.
\end{prop}

We can always find Euclidean inner products on $N$ such that the
decomposition $N \ = \ \sum_{i=1}^{m} V_{i}$ is an orthogonal sum. Let us pick
such an inner product and denote by $\| \cdot \|$ the Euclidean norm associated to
it.

We shall endow the group $N$ with a structure of a sub-Riemannian
manifold now. For this take the distribution obtained from left
translates of the space $V_{1}$. The metric on that distribution is
obtained by left translation of the inner product restricted to
$V_{1}$.

If $V_{1}$ Lie generates (the algebra) $N$
then any element $x \in N$ can be written as a product of
elements from $V_{1}$.  A slight reformulation of Lemma 1.40, Folland, Stein \cite{fostein}:

\begin{lema}
Let $N$ be a Carnot group and $X_{1}, ..., X_{p}$ an orthonormal basis
for $V_{1}$. Then there is a
 a natural number $M$ and a function $g: \left\{ 1,...,M \right\}
\rightarrow \left\{ 1,...,p\right\}$ such that any element
$x \in N$ can be written as:
\begin{equation}
x \ = \ \prod_{i = 1}^{M} \exp(t_{i}X_{g(i)})
\label{fp2.4}
\end{equation}
Moreover, if $x$ is sufficiently close (in Euclidean norm) to
$0$ then each $t_{i}$ can be chosen such that $\mid t_{i}\mid \leq C
\| x \|^{1/m}$.
\label{p2.4}
\end{lema}

 This means that there is a horizontal curve
joining any two points. Hence the distance
$$d(x,y) \ = \ \inf \left\{ \int_{a}^{b} \| c^{-1} \dot{c} \| \mbox{
d}t \ \mbox{ : } \ c(a) = x , \ c(b) = y , \ c^{-1} \dot{c} \in
V_{1}
\right\}$$
is finite for any two $x,y \in N$.
The distance is obviously left invariant.

Associate to the CC distance $d$ the function:
$$\mid x \mid_{d} \ = \ d(0,x)$$
The map $\mid \cdot \mid_{d}$ looks like a norm. However it is
intrinsically defined and hard to work with. This is the reason for
the introduction of homogeneous norms.

\begin{defi}
A continuous function from $x \mapsto \mid x \mid$ from $G$ to
$[0,+\infty)$ is a homogeneous norm if
\begin{enumerate}
\item[(a)] the set $\left\{ x \in G \ : \ \mid x \mid = 1 \right\}$
does not contain $0$.
\item[(b)] $\mid x^{-1} \mid = \mid x \mid$ for any $x \in G$.
\item[(c)] $\mid \delta_{\varepsilon} x \mid = \varepsilon
\mid x \mid$ for
any $x \in G$ and $\varepsilon > 0$.
\end{enumerate}
\label{d2.3}
\end{defi}

Homogeneous norms exist. $\mid \cdot \mid_{d}$ is a homogeneous norm. For any $p \in [0,+\infty)$ a homogeneous norm  is:
$$\mid x \mid_{p} \ =  \ \left[ \sum_{i=1}^{m} \mid x_{i} \mid^{p/i} \right]^{1/p} $$
For $p= +\infty$ the corresponding norm is 
$$\mid x \mid_{\infty} \ =  \  \max \left\{ \mid x_{i} \mid^{1/i} \right\} $$

\begin{prop}
Any two homogeneous norms are equivalent. Let $\mid \cdot \mid$ be a
homogeneous norm. Then the set
$\left\{ x \mbox{ : } \mid x \mid = 1 \right\}$ is compact.

There is a constant $C> 0$ such that for any $x,y \in N$ we have:
$$ \mid xy \mid \ \leq \ C \ ( \mid x \mid + \mid y \mid ) $$
\label{p11}
\end{prop}

There are two important consequences of this simple proposition. The first
is a description of the behaviour of any homogeneous norm. The following proposition implies also that the Euclidean topology and uniformity are the same
as the topology and uniformity induced by the Carnot-Carath\'eodory distance.

\begin{prop}
Let $\mid \cdot \mid$ be a homogeneous norm. Then there are
constants $c, C > 0$ such that for any $x \in N$, $\mid x \mid < 1$,
we have
$$c \| x \| \ \leq \ \mid x \mid \ \leq \  C \| x \|^{1/m}$$
\label{pest}
\end{prop}

The second consequence of proposition \ref{p11} says that the balls in CC-distance look roughly like boxes (as it is the case
with the Euclidean balls). A box is a set
$$Box(r) \ = \ \left\{ x \ = \ \sum_{i=1}^{m} x_{i} \ \mbox{ : }
\| x_{i} \| \leq r^{i} \right\}$$

\begin{prop} ("Ball-Box theorem")
There are positive constants $c,C$ such that for any $r>0$ we have:
$$Box(c r) \ \subset \ B(0,r) \ \subset \  Box( C r)$$
\label{bbprop}
\end{prop}

As a consequence the Lebesgue measure is absolutely continuous with
respect to the Hausdorff measure $\mathcal{H}^{Q}$. Because of the
invariance with respect to the group operation it follows that the
Lebesgue measure is a multiple of the mentioned Hausdorff measure.

The following theorem is Theorem 2. from Mitchell \cite{mit}, in the particular case of Carnot groups.

\begin{thm}
The ball $B(0,1)$ has Hausdorff dimension $Q$.
\end{thm}

\paragraph{Proof.}
We know that the volume of a ball with radius $\varepsilon$ is
$c \varepsilon^{Q}$.

Consider a maximal filling of $B(0,1)$ with balls of radius
$\varepsilon$. There are $N_{\varepsilon}$ such balls in the
filling; an upper bound for this number is:
$$ N_{\varepsilon} \ \leq \  1/\varepsilon^{Q}$$
The set of concentric balls of radius $2\varepsilon$ cover
$B(0,1)$; each of these balls has diameter smaller than $4
\varepsilon$, so the Hausdorff $\alpha$ measure of $B(0,1)$ is
smaller than
$$\lim_{\varepsilon \rightarrow 0} N_{\varepsilon}
(2\varepsilon)^{\alpha}$$
which is $0$ if $\alpha > Q$. Therefore the Hausdorff dimension is
smaller than $Q$.

Conversely, given any covering of $B(0,1)$ by sets of diameter $\leq
\varepsilon$, there is an associated covering with balls of the same
diameter; the number $M_{\varepsilon}$ of this balls has a lower
bound:
$$M_{\varepsilon} \ \geq \ 1/\varepsilon^{Q}$$
thus there is a lower bound
$$\sum_{cover} \varepsilon^{\alpha} \ \geq  \ \varepsilon^{\alpha} /
\varepsilon^{Q}$$
which shows that if $\alpha < Q$ then $\mathcal{H}^{\alpha}(B(0,1))
= \infty$. Therefore the Hausdorff dimension of the ball is greater
than $Q$.
\quad $\blacksquare$


We collect the important facts discovered until now:
let $N$ be a Carnot group endowed with the left invariant distribution
generated by $V_{1}$ and with an Euclidean norm on $V_{1}$.

\begin{enumerate}
\item[(a)] If $V_{1}$ Lie-generates the whole Lie algebra of $N$
then any two points can be joined by a horizontal path.
\item[(b)] The metric  topology and uniformity of $N$ are the same as Euclidean
topology and uniformity respective.
\item[(c)] The ball $B(0,r)$ looks roughly like the box
$\left\{ x \ = \ \sum_{i=1}^{m} x_{i} \ \mbox{ : }
\| x_{i} \| \leq r^{i} \right\}$.
\item[(d)] the Hausdorff measure $\mathcal{H}^{Q}$ is group
invariant and the Hausdorff dimension of a ball is $Q$.
\item[(e)] there is a one-parameter group of dilatations, where a
dilatation is an isomorphism $\delta_{\varepsilon}$ of $N$ which
transforms the distance $d$ in $\varepsilon d$.
\end{enumerate}

All these facts are true (under slight and almost obvious  modifications) in the
case of a sub-Riemannian manifold.

\paragraph{Pansu derivative. Rademacher theorem}

A Carnot group has it's own concept of differentiability, introduced by Pansu \cite{pansu}.

In Euclidean spaces, given $f: R^{n} \rightarrow R^{m}$ and
a fixed point $x \in R^{n}$, one considers the difference function:
$$X \in B(0,1) \subset R^{n} \  \mapsto \ \frac{f(x+ tX) - f(x)}{t} \in R^{m}$$
The convergence of the difference function as $t \rightarrow 0$ in
the uniform convergence gives rise to the concept of
differentiability in it's classical sense. The same convergence,
but in measure, leads to approximate differentiability. 
This and
another topologies might be considered (see Vodop'yanov
\cite{vodopis}, \cite{vodopis1}).

In the frame of Carnot groups the difference function can be written using only dilatations and the group operation. Indeed, for any function between Carnot groups
$f: G \rightarrow P$,  for  any fixed point $x \in G$ and $\varepsilon >0$  the finite difference function is defined by the formula:
$$X \in B(1) \subset G \  \mapsto \ \delta_{\varepsilon}^{-1} \left(f(x)^{-1}f\left(
x \delta_{\varepsilon}X\right) \right) \in P$$
In the expression of the finite difference function enters $\delta_{\varepsilon}^{-1}$ and $\delta_{\varepsilon}$, which are dilatations in $P$, respectively $G$.

Pansu's differentiability is obtained from uniform convergence of the difference
function when $\varepsilon \rightarrow 0$.

The derivative of a function $f: G \rightarrow P$ is linear in the sense
explained further.  For simplicity we shall consider only the case $G=P$. In this way we don't have to use a heavy notation for the dilatations.

\begin{defi}
Let $N$ be a Carnot group. The function
$F:N \rightarrow N$ is linear if
\begin{enumerate}
\item[(a)] $F$ is a {\it group} morphism,
\item[(b)] for any $\varepsilon > 0$ $F \circ \delta_{\varepsilon} \
= \ \delta_{\varepsilon} \circ F$.
\end{enumerate}
We shall denote by $HL(N)$ the group of invertible linear maps  of
$N$, called the  linear group of $N$.
\label{dlin}
\end{defi}

The condition (b) means that $F$, seen as an algebra morphism,
preserves the grading of $N$.

The definition of Pansu differentiability follows:

\begin{defi}
Let $f: N \rightarrow N$ and $x \in N$. We say that $f$ is
(Pansu) differentiable in the point $x$ if there is a linear
function $Df(x): N \rightarrow N$ such that
$$\sup \left\{ d(F_{\varepsilon}(y), Df(x)y) \ \mbox{ : } \ y \in B(0,1)
\right\}$$
converges to $0$ when $\varepsilon \rightarrow 0$. The functions $F_{\varepsilon}$
are the finite difference functions, defined by
$$F_{t} (y) \ = \ \delta_{t}^{-1} \left( f(x)^{-1} f(x
\delta_{t}y)\right)$$
\end{defi}

The definition says that $f$ is differentiable at $x$ if the
sequence of finite differences $F_{t}$ uniformly converges to a
linear map when $t$ tends to $0$.

We are interested to see how this differential looks like. For
any $f: N \rightarrow N$, $x,y \in N$ $Df(x)y$ means
$$Df(x)y \ =  \ \lim_{t \rightarrow 0} \delta_{t}^{-1}
\left( f(x)^{-1} f(x \delta_{t} y) \right)$$
provided that the limit exists.

\begin{prop}
Let $f: N \rightarrow N$, $y, z \in N$ such that:
\begin{enumerate}
\item[(i)] $Df(x)y$, $Df(x)z$ exist for any $x \in N$.
\item[(ii)] The map $x \mapsto Df(x)z$ is continuous.
\item[(iii)] $x \mapsto
\delta_{t}^{-1} \left( f(x)^{-1} f(x \delta_{t}z)\right)$
converges uniformly to $Df(x)z$.
\end{enumerate}
Then for any $a,b > 0$ and any $w \ = \ \delta_{a} x \delta_{b} y$
the limit $Df(x)w$ exists and we have:
$$Df(x) \delta_{a} x \delta_{b} y \ = \ \delta_{a} Df(x)y \
\delta_{b} Df(x) z$$
\label{ppansuprep}
\end{prop}

If $f$ is Lipschitz then the previous proposition holds almost
everywhere. In  Proposition 3.2, Pansu \cite{pansu} there is a more general
statement:

\begin{prop}
Let $f: N \rightarrow N$ have finite dilatation and suppose that for
almost any $x \in N$ the limits $Df(x)y$ and $Df(x) z$ exist. Then for
almost any $x$ and for any $w \ = \ \delta_{a} y \ \delta_{b} z$
the limit $Df(x) w$ exists and we have:
$$Df(x) \delta_{a} x \delta_{b} y \ = \ \delta_{a} Df(x)y \
\delta_{b} Df(x) z$$
\label{ppansu}
\end{prop}

\paragraph{Proof.}
The idea is to use  Proposition \ref{ppansuprep}. We shall suppose
first that $f$ is Lipschitz. Then any finite difference function
$F_{t}$ is also Lipschitz.

Recall the general Egorov theorem:

\begin{lema} (Egorov Theorem)
Let $\mu$ be a finite measure on $X$ and $(f_{n})_{n}$ a sequence of
measurable functions which converges $\mu$ almost everywhere (a.e.)
to $f$. Then for any $\varepsilon> 0$ there exists a measurable set
$X_{\varepsilon}$ such that $\mu(X \setminus X_{\varepsilon}) <
\varepsilon$ and $f_{n}$ converges uniformly to $f$ on
$X_{\varepsilon}$.
\end{lema}

\paragraph{Proof.}
It is not restrictive to suppose that $f_{n}$ converges pointwise to
$f$. Define, for any $q,p \in \mathbf{N}$, the set:
$$X_{q,p} \ = \ \left\{ x \in X \ \mbox{ : } \  \mid f_{n}(x) - f(x)
\mid \ < \ 1/p \ , \ \forall n \geq q \right\}$$
For fixed $p$ the sequence of sets $X_{q,p}$ is increasing and at
the limit it fills the space:
$$\bigcup_{q \in \mathbf{N}} X_{q,p} \ = \ X$$
Therefore for any $\varepsilon > 0$ and any $p$ there is a $q(p)$
such that $\mu(X \setminus X_{q(p) p}) \ < \ \varepsilon / 2^{p}$.

Define then $$X_{\varepsilon} \ = \ \bigcap_{p \in N} X_{q(p) p}$$
and check that it satisfies the conclusion.
\quad $\blacksquare$

If we could improve the conclusion of Egorov Theorem by claiming that if
$\mu$ is a Borel measure then each  $X_{\varepsilon}$ can be chosen
to be open, then the proof would resume like this. 

We can take $X$ to be a ball in  the group $N$ and
$\mu$ to be the $Q$ Hausdorff measure, which is Borel. Then we are
able to apply Proposition \ref{ppansuprep} on $X_{\varepsilon}$.
When we tend $\varepsilon$ to $0$ we obtain the claim (for $f$
Lipschitz).

Unfortunately, such an improvement in the conclusion of Egorov theorem 
does not hold. In fact, in Pansu paper, during the proof of his Rademacher theorem,  there is a missing step, as Vodop'yanov remarked. This missing step 
describes precisely what happens when $X_{\varepsilon}$ (a Borel measurable set) 
is replaced by an arbitrarily fine open approximation $X'_{\varepsilon}$. 
For this especially subtle step see Vodop'yanov in same proceedings [...].   

\quad $\blacksquare$

A consequence of this proposition is   Corollary 3.3, Pansu \cite{pansu}.

\begin{cor}
If $f: N \rightarrow N$ has finite dilatation and $X_{1}, ... , X_{K}$ is
a basis for the distribution $V_{1}$ such that
$$s \mapsto f(x \delta_{s} X_{i})$$ is differentiable in $s=0$ for
almost any $x \in N$, then $f$ is differentiable almost everywhere
and the differential is linear, i.e.  if
$$y \ = \ \prod_{i = 1}^{M} \delta_{t_{i}} X_{g(i)}$$
then
$$Df(x) y \ = \ \prod_{i = 1}^{M} \delta_{t_{i}} X_{g(i)}(f)(x)$$
\label{corr1}
\end{cor}

A very important result is
theorem 2, Pansu \cite{pansu}, which  contains the
Rademacher theorem for Carnot groups.

\begin{thm}
Let $f: M \rightarrow N$ be a Lipschitz  function between Carnot groups.
Then $f$ is differentiable almost everywhere.
\label{ppansuu}
\end{thm}

This is the kind of result that motivates the study of sub-Riemannian manifolds
as metric spaces endowed with a measure. It is to be noticed that  in the Rademacher theorem enters differentiability, measure and distance. This theorem
is indeed a concentrate of fundamental facts from analysis and measure theory,
hence its importance.

The proof is based on the corollary \ref{corr1} and the technique of development of a curve, Pansu sections 4.3 - 4.6 \cite{pansu}.

Let $c: [0,1] \rightarrow N$ be a Lipschitz curve such that $c(0) = 0$.  To
any division $\Sigma: 0=t_{0} < .... < t_{n} = 1$ of the interval $[0,1]$ is a associated the element $\sigma_{\Sigma} \in N$ (the algebra) given by:
$$\sigma_{\Sigma} \ = \ \sum_{k=0}^{n} c(t_{k})^{-1}c(t_{k+1})$$
Lemma (18) Pansu \cite{pansu1}  implies the existence of a constant $C>0$ such that
\begin{equation}
\| \sigma_{\Sigma} \ - \ c(1) \| \ \leq \ C \left(  \sum_{k=0}^{n}
\|c(t_{k})^{-1}c(t_{k+1})\| \right)^{2}
\label{pu1}
\end{equation}

Take now a finer division $\Sigma'$ and look at the interval
$[t_{k}, t_{k+1}]$ divided further by $\Sigma'$ like this:
$$t_{k} = t'_{l} < ... < t'_{m}=t_{k+1}$$
The estimate \eqref{pu1} applied for each interval $[t_{k}, t_{k+1}]$ lead
us to the inequality:
$$\| \sigma_{\Sigma'} \ - \ \sigma_{\Sigma} \| \ \leq \ \sum_{k=0}^{n} d(c(t_{k}), c(t_{k+1}))^{2}$$
The curve has finite length (being Lipschitz) therefore the right hand side
of the previous inequality tends to 0 with the norm of the division. Set
$$\sigma(s) \ = \ \lim_{\| \Sigma \| \rightarrow 0} \sigma_{\Sigma}(s)$$
where $\sigma_{\Sigma}(s)$ is relative to the curve $c$ restricted to the interval $[0,s]$. The curve such defined is called the development of
the curve $c$. It is easy to see that $\sigma$ has the same length
(measured with the Euclidean norm $\| \cdot \|$) as $c$. If we parametrise $c$
with the length then we have the estimate:
\begin{equation}
\|\sigma(s) \- \ c(s) \| \ \leq \ C s^{2}
\label{impes}
\end{equation}
This shows that $\sigma$ is a Lipschitz curve (with respect to the Euclidean distance). Indeed, prove first that $\sigma$ has finite dilatation in almost
any point, using \eqref{impes} and the fact that $c$ is Lipschitz. Then show that the dilatation is majorised by the Lipschitz constant of $c$.
By the classical Rademacher theorem $\sigma$ is almost everywhere derivable.

\begin{rk}
The development of a curve can be done in an arbitrary Lie connected group,
endowed with a left invariant distribution which generates the algebra.
One should add some logarithms, because in the case of Carnot groups, we have
identified the group with the algebra. The inequality \eqref{impes} still holds.
\label{impr}
\end{rk}

Conversely, given a curve $\sigma$ in the algebra $N$, we can perform the inverse operation to development (called multiplicative integral by Pansu;
we shall call it "lift"). Indeed, to any division $\Sigma$ of the interval
$[0,s]$ we associate
the point
$$c_{\Sigma}(s) \ = \ \prod_{k=0}^{n} \left( \sigma(t_{k+1}) \ - \ \sigma(t_{k}) \right)$$
Define then
$$c(s) \ = \ \lim_{\| \Sigma \| \rightarrow 0} c_{\Sigma}(s)$$
Remark that if $\sigma([0,1]) \subset V_{1}$ and it is almost everywhere differentiable then $c$ is a horizontal curve.

Fix $s \in [0,1)$ and parametrise $c$ by the length. Apply the inequality
\eqref{impes} to the curve $t \mapsto c(s)^{-1}c(s+t)$. We get the fact that
the vertical part of $\sigma(s+t) - \sigma(s)$ is controlled by
$s^{2}$. Therefore, if $c$ is Lipschitz then $\sigma$ is included in $V_{1}$.

Denote the $i$ multiple bracket $[x,[x,...[x,y]...]$ by $[x,y]_{i}$. For any
division $\Sigma$ of the interval $[0,s]$ set:
$$A^{i}_{\Sigma}(s) \ = \ \sum_{k=0}^{n} [\sigma(t_{k}),\sigma(t_{k+1})]_{i}$$
As before, one can show that when the norm of the division tends to 0,
$A^{i}_{\Sigma}$ converges. Denote by
$$A^{i}(s) \ = \ \int_{0}^{s} [\sigma, d\sigma]_{i} \ = \  \lim_{\| \Sigma \| \rightarrow 0} A^{i}_{\Sigma}(s)$$
the i-area function. The estimate corresponding to \eqref{impes}
is
\begin{equation}
\| A^{i}(s) \| \ \leq \ C s^{i+1}
\label{impesai}
\end{equation}

What is the significance of $A^{i}(s)$? The answer is simple, based on the well known formula of derivation of left translations
in a group. Define, in a neighbourhood of $0$ in the Lie algebra
$\mathfrak{g}$ of the Lie group $G$, the operation
$$X \opg Y \ = \ \log_{G}\left( \exp_{G}(X) \exp_{G}(Y) \right)$$
The left translation by $X$ is the function $L_{X}(Y) \ = \ X \opg Y$. It is known that
$$D L_{X}(0) (Z) \ = \ \sum_{i=0} \frac{1}{(i+1)!} [X,Z]_{i}$$
It follows that if the group $G$ is Carnot then $[X,Z]_{i}$ measures the infinitesimal variation of the $V_{i}$ component of $L_{X}(Y)$, for
$Y = 0$. Otherwise said,
$$\frac{d}{dt}\sigma_{i}(t) \ = \  \frac{1}{(i+1)!} \lim_{\varepsilon \rightarrow 0} \varepsilon^{-i} \ \int_{t}^{t + \varepsilon} [\sigma, d\sigma]_{i}$$
Because $\sigma$ is Lipschitz, the left hand side exists for all $i$ for a.e.
$t$. The estimate \eqref{impesai} tells us that the right hand side equals $0$. This implies the following proposition (Pansu, 4.1 \cite{pansu}).

\begin{prop}
If $c$ is Lipschitz  then $c$ is differentiable almost everywhere.
\label{pr41}
\end{prop}

\paragraph{Proof of theorem \ref{ppansuu}.}
The proposition \ref{pr41} implies that we are in the hypothesis of
corollary \ref{corr1}.
\quad $\blacksquare$

In order to prove Stepanov theorem, the 
Pansu-Rademacher theorem has been improved
for Lipschitz functions $f: A \subset M \rightarrow N$, where M, N
are Carnot groups and $A$ is just measurable, in Vodop'yanov,
Ukhlov \cite{voduk}.

 Their technique differs from Pansu. It resembles  with the one used in Margulis, Mostow \cite{marmos1}, where
it is proven that any quasi-conformal map from a Carnot-Carath\'eodory manifold
to another is a.e. differentiable, without using Rademacher theorem. Same result for quasi-conformal maps on Carnot groups has been first proved by Koranyi, Reimann \cite{kore}.
Finally  Magnani \cite{magnani} reproved a.e. differentiability of Lipschitz functions on Carnot groups, defined on measurable sets,  continuing Pansu technique.

For a review of other connected results, such as approximate differentiability,
differentiability in Sobolev topology, see the excellent Vodop'yanov papers
\cite{vodopis}, \cite{vodopis1}.

\subsection{Nilpotentisation of a sub-Riemannian manifold}
\label{nilpo}

In this section we return to the case of a connected real Lie group endowed with a nonintegrable left invariant distribution. We shall associate to such a sub-Riemannian space a Carnot group, using the procedure of nilpotentisation.

The filtration \eqref{filtra} has the straightforward property: if $x \in V^{i}$ and
$y \in V^{j}$ then $[x,y] \in V^{i+j}$ (where $V^{k} = \mathfrak{g}$ for all $k \geq m$ and $V^{0} = \left\{ 0 \right\}$). This allows to construct the Lie algebra
$$\mathfrak{n}(\mathfrak{g}, D) \ = \ \oplus_{i=1}^{m} V_{i} \ \ , \
V_{i} \ = \ V^{i}/V^{i-1}$$
with Lie bracket $ [\hat{x}, \hat{y}]_{n} \ = \ \widehat{[x,y]}$, where
$\hat{x} = x + V^{i-1}$ if $x \in V^{i}\setminus
V^{i-1}$.

\begin{prop}
$\mathfrak{n}(\mathfrak{g},D)$  with distinguished space $V_{1} \ = \ D$
is the Lie algebra of a Carnot group of step $m$. It is called the nilpotentisation of the filtration \eqref{filtra}.
\end{prop}
Remark that $\mathfrak{n}(\mathfrak{g},D)$ and $\mathfrak{g}$ have the same dimension, therefore they are isomorphic as vector spaces.

Let $\left\{ X_{1}, ... , X_{p} \right\}$ be a basis of the vector space $D$.
We shall build a basis of $\mathfrak{g}$ which will give a vector spaces isomorphism with $\mathfrak{n}(\mathfrak{g},D)$.

A word with letters $A \ = \ \left\{X_{1}, ... , X_{p}\right\}$ is a string
$X_{h(1)}...X_{h(s)}$ where $h: \left\{1 , ... , s \right\} \rightarrow \left\{ 1 , ... , p \right\}$. The set of words forms the dictionary $Dict(A)$, ordered lexicographically. We set the function $Bracket: Dict(A) \rightarrow
\mathfrak{g}$ to be
$$Bracket(X_{h(1)}...X_{h(s)}) \ = \ [X_{h(1)},[X_{h(2)},[ ... , X_{h(s)}] ... ] $$
For any $x \in Bracket(Dict(A))$ let
$\hat{x} \in Dict(A)$ be the least word  such that $Bracket(\hat{x}) \ = \ x$.
The collection of all these words is denoted by $\hat{g}$.

The length $l(x) \ = \ length(\hat{x})$ is well defined. The dictionary
$Dict(A)$ admit a filtration made by the length of words function. In the same
way the function $l$ gives the filtration
$$V^{1} \cap Bracket(Dict(A)) \ \subset \ V^{2} \cap Bracket(Dict(A)) \ \subset \ ... \ Bracket(Dict(A))$$
Choose now, in the lexicographic order in $\hat{g}$, a set $\hat{B}$ such that
$B \ = \ Bracket(\hat{B})$ is a  basis for $\mathfrak{g}$. Any  element $X$ in this basis can be written as
$$X \ = \ [X_{h(1)},[X_{h(2)},[ ... , X_{h(s)}] ... ] $$
such that $l(X) \ = \ s$ or equivalently $X \in V^{s} \setminus V^{s-1}$.

It is obvious that the map
$$X \in B \ \mapsto \tilde{X} \ = \ X + V^{l(X)-1} \in \mathfrak{n}(\mathfrak{g},D)$$
is a bijection and that $\tilde{B} \ = \ \left\{ \tilde{X_{j}} \mbox{ : }
j \ = \ 1, ... , dim \ \mathfrak{g} \right\}$ is a basis for
$\mathfrak{n}(\mathfrak{g},D)$. We can identify then $\mathfrak{g}$ with
$\mathfrak{n}(\mathfrak{g},D)$ by the identification $X_{j} \ = \ \tilde{X}_{j}$.

Equivalently we can define the nilpotent Lie bracket $[ \cdot , \cdot ]_{n}$
directly on $\mathfrak{g}$, with the use of the dilatations on $\mathfrak{n}(\mathfrak{g},D)$.

Instead of the filtration \eqref{filtra} let us start with a direct sum decomposition of $\mathfrak{g}$
\begin{equation}
\mathfrak{g} \ = \ \oplus_{i = 1}^{m} W_{i} \ , \ \ V^{i} \ = \ \oplus_{j=1}^{i} W_{j}
\label{deco}
\end{equation}
such that $[V^{i}, V^{j}] \subset V^{i+j}$. The chain $V^{i}$ form a filtration
like \eqref{filtra}.

We set
$$\delta_{\varepsilon} (x) \ = \ \sum_{i=1}^{m} \varepsilon^{i} x_{i}$$
for any $\varepsilon > 0$ and $x \in \mathfrak{g}$, which decomposes according
to \eqref{deco} as $x \ = \ \sum_{i=1}^{m} x_{i}$.

\begin{prop}
The limit
\begin{equation}
[x,y]_{n} \ = \ \lim_{\varepsilon \rightarrow 0} \delta_{\varepsilon}^{-1}
[\delta_{\varepsilon} x , \delta_{\varepsilon} y ]
\label{dnil}
\end{equation}
exists for any $x,y \in \mathfrak{g}$ and $(\mathfrak{g}, [ \cdot , \cdot]_{n})$ is the Lie algebra of a Carnot group with dilatations $\delta_{\varepsilon}$.
\label{limbra}
\end{prop}

\paragraph{Proof.}
Let $x \ = \ \sum_{j=1}^{i} x_{j}$ and $y \ = \ \sum_{k=1}^{l} y_{k}$. Then
$$[\delta_{\varepsilon} x , \delta-{\varepsilon} y ] \ = \
\sum_{s = 2}^{i+l} \varepsilon^{s} \sum_{j+k =s} \sum_{p=1}^{s} [x_{j},y_{k}]_{p}$$
We apply $\delta_{\varepsilon}^{-1}$ to this equality and we obtain:
$$\delta_{\varepsilon}^{-1}
[\delta_{\varepsilon} x , \delta_{\varepsilon} y ] \ = \
\sum_{s = 2}^{i+l}  \sum_{j+k =s} \sum_{p=1}^{s} \varepsilon^{s-p} [x_{j},y_{k}]_{p}$$
When $\varepsilon$ tends to 0 the expression converges to the  limit
\begin{equation}
[x,y]_{n} \ = \ \sum_{s = 2}^{i+l}  \sum_{j+k =s}  [x_{j},y_{k}]_{s}
\label{expre}
\end{equation}
For any $\varepsilon > 0$ the expression
$$[x,y]^{\varepsilon} \ = \ \delta_{\varepsilon}^{-1}
[\delta_{\varepsilon} x , \delta_{\varepsilon} y ]$$
is a Lie bracket (bilinear, antisymmetric and it satisfies the Jacobi identity). Therefore at the limit $\varepsilon \rightarrow 0$ we get a Lie bracket. Moreover, it is straightforward to see from the definition of
$[x,y]_{n}$ that $\delta_{\varepsilon}$ is an algebra isomorphism. We conclude that $(\mathfrak{g}, [ \cdot , \cdot]_{n})$ is the Lie algebra of a Carnot group with dilatations $\delta_{\varepsilon}$.
\quad $\blacksquare$

\begin{prop}
Let $\left\{X_{1}, ... ,X_{dim \ \mathfrak{g}} \right\}$ be a basis
of $\mathfrak{g}$ constructed from a basis $\left\{X_{1}, ... , X_{p}\right\}$
of $D$. Set
$$W_{j} \ = \ span \ \left\{ X_{i} \mbox{ : } l(X_{i}) \ = \ j \right\}$$
Then the spaces $W_{j}$ provides a direct sum decomposition \eqref{deco}. Moreover the identification $\hat{X}_{i} = X_{i}$ gives a Lie algebra  isomorphism between
$(\mathfrak{g}, [ \cdot , \cdot ]_{n})$ and $\mathfrak{n}(\mathfrak{g}, D)$.
\end{prop}

\paragraph{Proof.}
In the basis $\left\{X_{1}, ... ,X_{dim \ \mathfrak{g}} \right\}$ the Lie bracket on $\mathfrak{g}$ looks like this:
$$[X_{i}, X_{j}] \ = \ \sum C_{ijk}X_{k}$$
where $c_{ijk} = 0$ if $l(X_{i}) + l(X_{j}) < l(X_{k})$. From here the first part of the proposition is straightforward. The expression of the Lie bracket generated by the decomposition \eqref{deco} is obtained from  \eqref{expre}. We have
$$[X_{i}, X_{j}]_{n} \ = \ \sum \lambda_{ijk} C_{ijk} X_{k}$$
where $\lambda_{ikj} = 1$ if $l(X_{i}) + l(X_{j}) = l(X_{k})$ and $0$ otherwise.  The Lie algebra isomorphism follows from the expression of the Lie bracket on $\mathfrak{n}(\mathfrak{g},D)$:
$$[\hat{X}_{i}, \hat{X}_{j}] \ = \ [X_{i}, X_{j}] \ + \ V^{l(X_{i}) + l(X_{j})-1}$$
\quad $\blacksquare$

In conclusion  the expression of the nilpotent Lie
bracket depends on the choice of basis $B$ trough the transport of the dilatations group from $\mathfrak{n}(\mathfrak{g},D)$ to $\mathfrak{g}$.

Let $N(G,D)$ be the simply connected Lie group with Lie algebra
$\mathfrak{n}(\mathfrak{g},D)$. As previously, we identify $N(G,D)$ with
$\mathfrak{n}(\mathfrak{g},D)$ by the exponential map.

\subsection{Gromov-Hausdorff convergence and metric tangent cones}

This section contains the general notions and results which allows to say what a
tangent space to a metric space is.

\paragraph{Distances between metric spaces}

The references for this section  are Gromov \cite{gromov}, chapter 3, and Burago \& al. \cite{burago} section 7.4.  There are several definitions of distances between metric spaces. The very fertile idea of introducing such distances belongs to Gromov.

In order to introduce the Hausdorff distance between metric spaces, recall
the Hausdorff distance between subsets of a metric space.

\begin{defi}
For any set $A \subset X$ of a metric space and any $\varepsilon > 0$ set
the $\varepsilon$ neighbourhood of $A$ to be
$$A_{\varepsilon} \ = \ \cup_{x \in A} B(x,\varepsilon)$$
The Hausdorff distance between $A,B \subset X$ is defined as
$$d_{H}^{X}(A,B) \ = \ \inf \left\{ \varepsilon > 0 \mbox{ : } A \subset B_{\varepsilon} \ , \ B \subset A_{\varepsilon} \right\}$$
\end{defi}

By considering all isometric embeddings of two metric spaces $X$, $Y$ into
an arbitrary metric space $Z$ we obtain the Hausdorff distance between $X$, $Y$ (Gromov \cite{gromov} definition 3.4).

\begin{defi}
The Hausdorff distance $d_{H}(X,Y)$ between metric spaces $X$ $Y$ is the infimum of the numbers
$$d_{H}^{Z}(f(X),g(Y))$$
for all isometric embeddings $f: X \rightarrow Z$, $g: Y \rightarrow Z$ in a
metric space $Z$.
\end{defi}

If $X$, $Y$ are compact then $d_{H}(X,Y) < + \infty$. Indeed, let $Z$ be the disjoint union of $X,Y$ and $M \ = \ \max \left\{ diam(X) , diam(Y) \right\}$.
Define the distance on $Z$ to be
$$d^{Z}(x,y) \ = \ \left\{ \begin{array}{ll}
d^{X}(x,y) & x, y \in X \\
d^{Y}(x,y) & x, y \in Y \\
\frac{1}{2} M & \mbox{ otherwise}
\end{array} \right. $$
Then $d_{H}^{Z}(X,Y) < + \infty$.

The Hausdorff distance between isometric spaces equals $0$. The converse is also true (Gromov {\it op. cit.} proposition 3.6) in the class of compact metric spaces.

\begin{thm}
If $X,Y$ are compact metric spaces such that $d_{H}(X,Y) = 0$ then
$X,Y$ are isometric.
\label{tgro}
\end{thm}

For the proof of the theorem we need the Lipschitz distance ({\it op. cit. definition 3.1}) and a criterion for convergence of metric spaces in the Hausdorff distance ( {\it op. cit. proposition 3.5 }). We shall give the definition of Gromov for Lipschitz distance and the first part of the
mentioned proposition.

\begin{defi}
The Lipschitz distance $d_{L}(X,Y)$ between bi-Lipschitz homeomorphic metric spaces
$X,Y$ is the infimum of
$$\mid \log \ dil(f) \mid + \mid \log \ dil(f^{-1}) \mid$$
for all $f: X \rightarrow Y$, bi-Lipschitz homeomorphisms.
\end{defi}

Obviously, if $d_{L}(X,Y) = 0$ then $X,Y$ are isometric. Indeed, by definition we have a sequence $f_{n}: X \rightarrow Y$ such that $dil (f_{n}) \rightarrow 1$ as $n \rightarrow \infty$. Extract an uniformly convergent subsequence; the limit is an isometry.

\begin{defi}
A $\varepsilon$-net in the metric space $X$ is a set $P \subset X$ such that
$P_{\varepsilon} \ = \ X$. The separation of the net $P$ is
$$sep(P) \ = \ \inf \left\{ d(x,y) \mbox{ : } x \not = y \ , \ x,y \in P \right\}$$
A $\varepsilon$-isometry between $X$ and $Y$ is a function $f: X \rightarrow Y$ such that $dis \ f \ \leq \varepsilon$ and $f(X)$ is a $\varepsilon$ net in
$Y$.
\end{defi}

The following proposition gives a connection between  convergence of metric spaces
in the Hausdorff distance and  convergence of $\varepsilon$ nets in the Lipschitz distance.

\begin{defi}
A sequence of metric spaces $(X_{i})$ converges in the sense of Gromov-Hausdorff to the metric space $X$ if $d_{H}(X,X_{i}) \rightarrow 0$ as
$i \rightarrow \infty$. This means that there is a sequence $\eta_{i} \rightarrow 0$ and isometric embeddings $f_{i}: X \rightarrow Z_{i}$,
$g_{i} : X_{i} \rightarrow Z_{i}$ such that
$$d_{H}^{Z_{i}} (f_{i}(X), g_{i}(X_{i})) \ < \ \eta_{i}$$
\label{dhauc}
\end{defi}

The next proposition is corollary 7.3.28 (b), Burago \& al. \cite{burago}. The proof is adapted from Gromov, proof of proposition 3.5 (b).

\begin{prop}
If there exists a $\varepsilon$-isometry between $X,Y$ then $d_{H}(X,Y) < 2 \varepsilon$.
\label{pbur}
\end{prop}

\paragraph{Proof.}
Let $f: X \rightarrow Y$ be the $\varepsilon$-isometry. On the disjoint union
$Z = X \cup Y$ extend the distances $d^{X}$, $d^{Y}$ in the following way. Define the distance between $x \in X$ and $y \in Y$ by
$$d(x,y) \ = \ \inf\left\{ d^{X}(x,u) + d^{Y}(f(u),y) + \varepsilon \right\}$$
This gives a distance $d^{Z}$ on $Z$. Check that
$d^{Z}_{H}(X,Y) < 2 \varepsilon$.
\quad $\blacksquare$

The proof of theorem \ref{tgro} follows as an application of previous propositions.

\paragraph{Metric tangent  cones}

The infinitesimal (not local!) geometry of a path metric space is described with the help of Gromov-Hausdorff convergence of pointed metric spaces.

This is definition 3.14 Gromov \cite{gromov}.
\begin{defi}
The sequence of pointed metric spaces $(X_{n}, x_{n}, d_{n})$ converges in the sense
of Gromov-Hausdorff to the pointed space $(X,x, d_{X})$ if
for any $r>0$, $\varepsilon > 0$ there is $n_{0} \in N$ such that for all
$n \geq n_{0}$ there exists $f_{n} : B_{n}(x_{n},r) \subset X_{n} \rightarrow X$ such that:
\begin{enumerate}
\item[(1)] $f_{n}(x_{n}) \ = \ x$,
\item[(2)] the distorsion of $f_{n}$ is bounded by $\varepsilon$: $dis \ f_{n} \ < \ \varepsilon$,
\item[(3)]
$B_{X}(x,r) \ \subset \ \left(f_{n}(B_{n}(x_{n},r)) \right)_{\varepsilon}$.
\end{enumerate}
\end{defi}

The Gromov-Hausdorff limit is defined up to isometry, in the class of compact metric spaces (Proposition 3.6 Gromov \cite{gromov}) or in the class of locally compact cones. Here it is the definition of a cone.

\begin{defi}
A pointed metric space $(X,x_{0})$ is called a cone if there is a one parameter group of dilatations with center $x_{0}$  $\left\{ \delta_{\lambda}: X \rightarrow X \mbox{ : } \lambda >0 \right\}$
 such that $\delta_{\lambda}(x_{0}) \ =
\ x_{0}$ and for any $x,y \in X$
$$d(\delta_{\lambda}(x), \delta_{\lambda}(y)) \ = \ \lambda \ d(x,y)$$
\end{defi}

The infinitesimal geometry of a metric space $X$ in the neighbourhood of $x_{0} \in X$ is described by the tangent space to $(X,x_{0})$ (if such object exists).

\begin{defi}
The tangent space to $(X,x_{0})$ is the Gromov-Hausdorff limit
$$(T_{x_{0}}, 0, d_{x_{0}})  \ = \ \lim_{\lambda \rightarrow \infty} (X,x_{0},
\lambda d)$$
\end{defi}

Three remarks are in order:
\begin{enumerate}
\item[1.] The tangent space is obviously a cone. We shall see that in a large class of situations is also a group, hence a graded nilpotent group, called
for short Carnot group.
\item[2.] The tangent space comes with a metric inside. This space is
path metric (Gromov \cite{gromov} Proposition 3.8).
\item[3.] The tangent cone is defined up to isometry, therefore there is
no way to use the metric tangent cone definition to construct a tangent bundle.
For a modification of the definitions in this direction  see Margulis, Mostow \cite{marmos2}, or 
Vodop'yanov \& Greshnov \cite{vodopis2}, \cite{vodopis3}.
\end{enumerate}

\subsection{Mitchell's theorem 1. Bella\"{\i}che's construction}

We collect here three key items in the edification of sub-Riemannian geometry.
The first is Mitchell \cite{mit} theorem 1:

\begin{thm}
For a regular sub-Riemannian space $(M,D,g)$, the tangent cone of
$(M,d_{CC})$ at $x \in M$ exists and it is isometric to $(N(D), d_{N})$,
which is a Carnot group with a left invariant distribution and $d_{N}$ is
a induced Carnot-Carath\'eodory distance.
\label{mite1}
\end{thm}
We shall prove this
theorem in the particular case of  a Lie group with left-invariant
distribution. For a detailed  proof of the general case see Vodop'yanov \& Greshnov \cite{vodopis2} \cite{vodopis3}.

The Carnot group $N(D)$ is called the nilpotentisation of the
regular distribution $D$ and it can be constructed from $D$
exclusively. However, the metric on $N(D)$ depends on the choice
of metric $g$ on the sub-Riemannian manifold. 

Recall that the limit in the Gromov-Hausdorff sense is defined up to isometry.
This means it this case that $N(D)$ is a model for the tangent space at $x$ to
$(M,d_{CC})$. In the Riemannian case $D = TM$ and $N(D) = R^{n}$, as a group with
addition. 

This theorem tells us nothing about the tangent bundle. 
There are however other ways to associate a tangent bundle to a metric measure  space
(Cheeger \cite{cheeger}) or to a regular sub-Riemannian manifold (Margulis \& Mostow
\cite{marmos1}, \cite{marmos2}). These bundles differs. As Tyson (interpreting Cheeger) asserts (see Tyson paper in these proceedings (??)), Cheeger tangent bundle can be identified with the distribution $D$ and Margulis-Mostow bundle is  the same as
the usual tangent bundle, but with the fiber isomorphic with $N(D)$, instead of
$R^{n}$.  We shall explain how the Margulis-Mostow tangent bundle is constructed
a bit further (again in the particular case considered here).

Let us not, for the moment, be too ambitious and restrict to the question:
is there a metric derivation of the group operation on $N(D)$? Bella\"{\i}che \cite{bell} writes that 
he asked Gromov this question, who pointed out that the key tool to construct the operation from 
metric is uniformity. Bella\"{\i}che proposed therefore the following construction, which starts from the proof of Mitchell theorem 1, where it can be seen that the Gromov-Hausdorff convergence
to the tangent space {\it is uniform with respect to $x \in M$}. This means that for any $\varepsilon > 0$ there is $R(\varepsilon) > 0$ and map
$\phi_{x,\varepsilon}: B_{CC}(x,R(\varepsilon))  \rightarrow N(D)$
 such that
$$ d_{N}(\phi_{x,\varepsilon}(y), \phi_{x,\varepsilon}(z)) = d_{CC}(y,z) + o(\varepsilon) \ \ \forall y,z \in B_{CC}(x,R(\varepsilon))$$
Let us forget about $\varepsilon$ (Bella\"{\i}che does not mention anything about it further) and take arbitrary  $X,Y \in N(D)$. Pick then   $y \in M$ such that $\phi_{x}(y) = X$.
Denote $\displaystyle \phi_{xy} = \phi_{y} \phi_{x}^{-1}$.
Then the operation in $N(D)$ is defined by:
$$X \opn Y = \lim_{\lambda \rightarrow \infty} \phi_{xy}^{-1} \delta_{\lambda}^{-1}
\phi_{xy} \delta_{\lambda} (Y)$$

It is easier to understand this in the Euclidean case, that is in $R^{n}$. We can take for example
$$\phi_{x}(y) = Q(x)(y-x) = X \ \ , y = x+Q^{-1}X \ \ , x \mapsto Q(x) \in SO(n)
\mbox{ arbitrary}$$
(and we have no dependence on $\varepsilon$)
Let us compute the operation. We get
$$\lim_{\lambda \rightarrow \infty} \phi_{xy}^{-1} \delta_{\lambda}^{-1}
\phi_{xy} \delta_{\lambda} (Y) = \ X+Y$$
as expected. Notice that the arbitrary choice of the rotations $Q(x)$ does not influence the result. The tangent spaces at any point can rotate independently, which is a sign that this construction cannot lead to a tangent bundle.

There are several problems with Bella\"{\i}che's construction:
\begin{enumerate}
\item[a)] when $\lambda \rightarrow \infty$ the expression $\phi_{xy} \delta_{\lambda} (Y)$ might  not make sense,
\item[b)] is not clear how $\varepsilon$ and $\lambda$ interact.
\end{enumerate}

This is the reason why we introduced (first in \cite{buliga1}, then here) the notion of uniform group, which encodes all that one really need to do the construction, again in the case of Lie groups with left invariant distributions.

Another way to transform Bella\"{\i}che's construction into an effective one (and more, to obtain a tangent bundle) is proposed by Margulis and Mostow. We shall explain further their construction. However, there are other problems emerging. as mentioned in the introduction.

\subsection{Margulis \& Mostow tangent bundle}
In this section we shall apply Margulis \& Mostow \cite{marmos2} construction of the tangent bundle to a SR manifold  for the case of a group with left invariant distribution. It will turn that the tangent bundle does not have a group structure, due to the fact that, as previously, the non-smoothness of the right translations is not studied.

The main point in the construction of a tangent bundle is to have a functorial definition of the tangent space. This is achieved by Margulis \& Mostow \cite{marmos2} in a very natural way. One of the geometrical definitions of a tangent
vector $v$ at a point $x$, to a manifold $M$, is the following one: identify $v$ with the class of smooth curves which pass through $x$ and have tangent $v$.
If the manifold $M$ is endowed with a distance then one can  define the equivalence relation based in $x$ by:  $c_{1} \equiv_{x}
c_{2}$ if $c_{1}(0)  = c_{2}(0) = x$ and the distance between $c_{1}(t)$ and
$c_{2}(t)$ is negligible with respect to  $t$ for small $t$.  The set of equivalence classes
is the tangent space at $x$. One has to put then some structure on the tangent
space (as, for example, the nilpotent multiplication).

To put is practice this idea is not so easy though. This is achieved by the following sequence of definitions and theorems. For commodity we shall explain this construction in the case $M=G$ connected Lie group, endowed with a left
invariant distribution $D$. The general case is the one of a regular sub-Riemannian manifold. We shall denote by $d_{G}$ the CC distance on $G$ and
we identify $G$ with $\mathfrak{g}$, as previously. The CC distance induced by
the distribution $D^{N}$, generated by left translations of $G$ using nilpotent
multiplication $\opn$, will be denoted by $d_{N}$.

\begin{defi}
A $C^{\infty}$ curve in $G$ with  $x \ = \ c(0)$ is called rectifiable at $t=0$ if
$d_{G}(x, c(t)) \leq Ct$ as $t \rightarrow 0$.

Two $C^{\infty}$ curves $c', c"$ with $c'(0) = x = c"(0)$ are called equivalent at $x$ if $$t^{-1} d_{G}(c'(t),c"(t)) \rightarrow 0$$ as $t \rightarrow 0$.

The tangent cone to $G$ as $x$, denoted by $C_{x} G$ is the set of equivalence classes of all $C^{\infty}$ paths $c$ with $c(0) = x$, rectifiable at $t=0$.
\end{defi}

Let $c: [-1,1] \rightarrow G$ be a $C^{\infty}$ rectifiable curve, $x = c(0)$ and
\begin{equation}
v \ = \ \lim_{t \rightarrow 0} \delta_{t}^{-1} \left( c(0)^{-1} \opg c(t) \right)
\label{ax1}
\end{equation}
The limit $v$ exists because the curve is rectifiable.

Introduce the curve $c_{0}(t) \ = \ x \exp_{G}(\delta_{t}v)$. Then
$$d(x, c_{0}(t)) \ = \ d(e, x^{-1}c(t)) < \mid v \mid t$$
 as $t \rightarrow 0$  (by the Ball-Box theorem)
The curve  $c$  is equivalent with $c_{0}$. Indeed, we have (for $t>0$):
$$\frac{1}{t} d_{G}(c(t), c_{0}(t)) \ = \ \frac{1}{t} d_{G}(c(t), x \opg \delta_{t} v)  \ = \ \frac{1}{t} d_{G}(\delta_{t}(v^{-1}) \opg x^{-1} \opg c(t) , 0)$$
The latter expression is equivalent (by the Ball-Box Theorem) with
$$\frac{1}{t} d_{N} (\delta_{t}(v^{-1}) \opg x^{-1} \opg c(t) , 0) \ = \
d_{N} \left(\delta_{t}^{-1} \left(  \delta_{t}(v^{-1}) \opg \delta_{t}\left(
\delta_{t}^{-1} \left( x^{-1} \opg c(t)\right) \right) \right) \right) $$
The right hand side (RHS) converges to
$d_{N}(v^{-1} \opn v , 0)$, as $t \rightarrow 0$, as a consequence of
the definition of $v$ and theorem \ref{teore}.

Therefore we can identify $C_{x} G$ with the set of curves
$t \mapsto x \exp_{G}(\delta_{t} v)$, for all $v \in \mathfrak{g}$. Remark that the equivalence relation between curves $c_{1}$, $c_{2}$, such that
$c_{1}(0) \ = \ c_{2}(0) \ = \ x$ can be redefined as:
\begin{equation}
\lim_{t \rightarrow 0} \delta_{t}^{-1} \left( c_{2}(t)^{-1} \opg  c_{1}(t) \right) \ = \ 0
\label{ax2}
\end{equation}

In order to define the multiplication Margulis \& Mostow introduce  the families of segments rectifiable at  $t$.

\begin{defi}
A family of segments rectifiable at $t=0$ is a $C^{\infty}$ map
$$\mathcal{F} : U  \rightarrow G$$
where $U$ is an open neighbourhood of $G \times 0$ in $G \times R$ satisfying
\begin{enumerate}
\item[(a)] $\mathcal{F}(\cdot, 0) \ = \ id $
\item[(b)] the curve $t \mapsto \mathcal{F}(x,t)$ is rectifiable at $t=0$ uniformly for all $x \in G$, that is for every compact  $K$  in $G$ there is a  constant $C_{K}$ and a compact neighbourhood $I$ of $0$ such that $d_{G}(y,\mathcal{F}(y,t)) < C_{K} t $ for all $(y,t) \in K \times I$.
\end{enumerate}

Two families of segments rectifiable at $t=0$ are called equivalent if 
$$t^{-1} d_{G}(\mathcal{F}_{1}(x,t), \mathcal{F}_{2}(x,t)) \rightarrow 0$$ as $t\rightarrow 0$,  uniformly on compact sets in the domain of definition.
\end{defi}

Part (b) from the  definition of a family of segments rectifiable can be restated as:
there exists the limit
\begin{equation}
v(x) \ = \ \lim_{t \rightarrow 0} \delta_{t}^{-1} \left( x^{-1} \opg \mathcal{F}(x,t) \right)
\label{ax3}
\end{equation}
and the limit is uniform with respect to $x \in K$, $K$ arbitrary compact set.

It follows then, as previously, that $\mathcal{F}$ is equivalent to
$\mathcal{F}_{0}$, defined by:
$$\mathcal{F}_{0}(x,t) \ = \ x \opg \delta_{t} v(x)$$
Also, the equivalence between families of segments rectifiable can be redefined
as:
\begin{equation}
\lim_{t \rightarrow 0} \delta_{t}^{-1} \left( \mathcal{F}_{2}(x,t)^{-1} \opg  \mathcal{F}_{1}(x,t) \right) \ = \ 0
\label{ax22}
\end{equation}
uniformly with respect to $x \in K$, $K$ arbitrary compact set.

\begin{defi}
The product of two families $\mathcal{F}_{1}$, $\mathcal{F}_{2}$ of segments rectifiable at $t=0$ is defined by
$$ \left( \mathcal{F}_{1} \circ \mathcal{F}_{2}\right) (x,t) \ = \ \mathcal{F}_{1}(\mathcal{F}_{2}(x,t),t)$$
\end{defi}

The product is well defined by Lemma 1.2 {\it op. cit.}. One of the main results is then the following theorem (5.5).

\begin{thm}
Let $c_{1}$, $c_{2}$ be $C^{\infty}$ paths rectifiable at $t=0$, such that
$c_{1}(0) = x_{0} = c_{2}(0)$. Let
$\mathcal{F}_{1}$, $\mathcal{F}_{2}$ be two families of segments rectifiable
at $t=0$ with:
$$\mathcal{F}_{1}(x_{0}, t) \ = \ c_{1}(t) \ \ , \ \ \mathcal{F}_{2}(x_{0}, t) \ = \ c_{2}(t)$$
Then the equivalence class of
$$t \mapsto \mathcal{F}_{1} \circ \mathcal{F}_{2}(x_{0}, t)$$ depends
only on the equivalence classes of $c_{1}$ and $c_{2}$. This defines the product of the elements of the tangent cone $C_{x_{0}} G$.
\end{thm}

This theorem is the straightforward consequence of the following facts (5.1(5) and 5.2 in Margulis \& Mostow \cite{marmos2}).

 We shall denote by $\mathcal{F} \approx \mathcal{F}'$ the equivalence relation of families of segments rectifiable; the equivalence relation of rectifiable curves based at $x$ will be denoted by $c \stackrel{x}{\approx} c'$.

\begin{lema}
\begin{enumerate}
\item[(a)] Let $\mathcal{F}_{1} \approx \mathcal{F}_{2}$ and
$\mathcal{G}_{1} \approx \mathcal{G}_{2}$. Then $\mathcal{F}_{1}
\circ \mathcal{G}_{1} \approx \mathcal{F}_{2} \circ
\mathcal{G}_{2}$. \item[(b)] The map $\mathcal{F} \mapsto
\mathcal{F}_{0}$ is constant on equivalence classes of families of
segments rectifiable.
\end{enumerate}
\label{lhelps}
\end{lema}

\paragraph{Proof.}
Let $$\mathcal{F}_{0} (x,t) \ = \ x \opg \delta_{t} w_{1}(x) \ \ , \ \
\mathcal{G}_{0} (x,t) \ = \ x \opg \delta_{t} w_{2}(x)$$
For the point (a) it is sufficient to prove that
$$\mathcal{F} \circ \mathcal{G}  \approx  \mathcal{F}_{0} \circ \mathcal{G}_{0}$$
This is true by the following chain of estimates.
$$\frac{1}{t} d_{G}(\mathcal{F} \circ \mathcal{G}(x,t) , \mathcal{F}_{0} \circ \mathcal{G}_{0}(x,t)) \ = $$
$$= \ \frac{1}{t} d_{G}(\delta_{t}w_{1}(\mathcal{G}_{0}(x,t))^{-1} \opg \delta_{t}w_{2}(x)^{-1} \opg x^{-1} \opg \mathcal{F}(\mathcal{G}(x,t),t) ,  0)$$
The RHS of this equality behaves like
$$d_{N}(\delta_{t}^{-1} \left( \delta_{t}w_{1}(\mathcal{G}_{0}(x,t))^{-1}
\opg \delta_{t}w_{2}(x)^{-1} \opg \delta_{t} \left( \delta_{t}^{-1} \left(
x^{-1} \opg \mathcal{G}(x,t)\right) \right) \right. \opg $$
$$\opg \delta_{t} \left( \left.
\delta_{t}^{-1} \left( \mathcal{G}(x,t)^{-1} \opg \mathcal{F}(\mathcal{G}(x,t),t) \right) \right) \right) , 0)$$
This quantity converges (uniformly with respect to $x \in K$, $K$ an arbitrary compact) to
$$d_{N}(w_{1}(x)^{-1} \opn w_{2}(x)^{-1} \opn w_{2}(x) \opn w_{1}(x), 0) \ = \ 0$$

The point (b) is easier: let $\mathcal{F} \approx \mathcal{G}$ and consider
$\mathcal{F}_{0}$, $\mathcal{G}_{0}$, as above. We want to prove that
$\mathcal{F}_{0} \ = \ \mathcal{G}_{0}$, which is equivalent to $w_{1} \ = \
w_{2}$.

Because $\approx $ is an equivalence relation all we have to prove is that
if $\mathcal{F}_{0} \approx \mathcal{G}_{0}$ then $w_{1} \ = \ w_{2}$.  We have:
$$\frac{1}{t} d_{G}(\mathcal{F}_{0}(x,t), \mathcal{G}_{0}(x,t)) \ = \
\frac{1}{t} d_{G}(x \opg \delta_{t} w_{1}(x) , x \opg \delta_{t} w_{2}(x))$$
We use the $\opg$ left invariance of $d_{G}$ and the Ball-Box theorem to
deduce that the RHS behaves like
$$d_{N}(\delta_{t}^{-1} \left( \delta_{t} w_{2}(x)^{-1} \opg \delta_{t} w_{1}(x)^{-1}\right) , 0)$$
which converges to $d_{N}(w_{1}(x), w_{2}(x))$ as $t$ goes to $0$. The two families are equivalent, therefore the limit equals $0$, which implies that
$w_{1}(x) = w_{2}(x)$ for all $x$.
\quad $\blacksquare$

We shall apply this theorem. Let $c_{i}(t) \ = \ x_{0} \exp_{G} \delta_{t} v_{i}$, for $i = 1,2$. It is easy to check that
$\mathcal{F}_{i}(x,t) \ = \ x \exp_{G}(\delta_{t} v_{i})$ are families of segments rectifiable at $t=0$ which satisfy the hypothesis of the theorem.
But then
$$\left( \mathcal{F}_{1} \circ \mathcal{F}_{2}\right) (x,t) \ = \
x_{0} \exp_{G} (\delta_{t} v_{1}) \exp_{G}(\delta_{t} v_{2})$$
which is equivalent with
$$\mathcal{F} \exp_{G}\left( \delta_{t}( v_{1} \opn v_{2} ) \right)$$
Therefore the tangent bundle defined by this procedure is the same as the virtual tangent bundle which we shall define soon, inspired from the construction proposed by Bella\"{\i}che.

 Maybe I misunderstood the notations, but it
seems to me that several times the authors claim that the exponential map which they construct is bi-Lipschitz (as in 5.1(4) and Corollary 4.5). This is false, as explained before. In Bella\"{\i}che \cite{bell}, Theorem 7.32 and
also at the beginning of section 7.6 we find that the exponential map is only
$1/m$ H\"{o}lder continuous (where $m$ is the step of the nilpotentization).
However,  the final results of Margulis \& Mostow hold true, if not entirely proven facts, as statements at least.

\subsection{Vodop'yanov \& Greshnov definition of the derivability. Rademacher theorem}

We choose Vodop'yanov \& Greshnov \cite{vodopis2}, \cite{vodopis3}
way of defining the derivability in order to explain Margulis \& Mostow \cite{marmos1}
Rademacher theorem 10.5 (or Vodop'yanov \& Greshnov theorem 1).
The definition (10.3.1) of derivability in the paper \cite{marmos1} is to be compared with the definition of $\mathcal{P}$ differentiability
from \cite{vodopis2}, first page, which is in my opinion clearer. However, the
reader entering for the first time in this subject might find hard to understand
why such elementary notions as differentiability need so lengthy discussions. It is, I think, another sign of the fact that the foundations of non-Euclidean analysis are still in construction.

We stay, as previously, in the case of Lie groups with left invariant distributions. We put on such groups  Lebesgue measures coming from
arbitrary atlases.

\begin{defi}
A mapping $f: G_{1} \rightarrow G_{2}$ is said to be differentiable at $x \in G_{1}$ if the mapping $ \exp_{G_{2}}^{-1} \circ L_{f(x)} \circ f \circ L_{x} \circ
\exp_{G_{1}}$ is Pansu derivable at $0$, when we identify the algebras
$\mathfrak{g}_{1}$, $\mathfrak{g}_{2}$ with the nilpotentisations of $G_{1}$,
$G_{2}$ respectively.
\label{dpmmv}
\end{defi}
The following theorem then holds. The original (and stronger) versions of this theorem concern 
quasi-conformal mapping and can be found in  Margulis \& Mostow \cite{marmos1},
Vodop'yanov \& Greshnov \cite{vodopis2} and the paper by Vodop'yanov in these proceedings.

\begin{thm} Any Lipschitz map $f: E \subset G_{1} \rightarrow G_{2}$, $E$ measurable, is derivable almost everywhere.
\label{tvodmarmos}
\end{thm}

\section{Uniform groups}

We start with the following setting: $G$ is a topological group endowed with an uniformity such that the operation is uniformly continuous.  More specifically,
we introduce first the double of $G$, as the group $G^{(2)} \ = \ G \times G$ with operation
$$(x,u) (y,v) \ = \ (xy, y^{-1}uyv)$$
The operation on the group $G$, seen as the function
$$op: G^{(2)} \rightarrow G \ , \ \ op(x,y) \ = \ xy$$
is a group morphism. Also the inclusions:
$$i': G \rightarrow G^{(2)} \ , \ \ i'(x) \ = \ (x,e) $$
$$i": G \rightarrow G^{(2)} \ , \ \ i"(x) \ = \ (x,x^{-1}) $$
are group morphisms.

\begin{defi}
\begin{enumerate}
\item[1.]
$G$ is an uniform group if we have two uniformity structures, on $G$ and
$G^{2}$,  such that $op$, $i'$, $i"$ are uniformly continuous.

\item[2.] A local action of a uniform group $G$ on a uniform  pointed space $(X, x_{0})$ is a function
$\phi \in W \in \mathcal{V}(e)  \mapsto \hat{\phi}: U_{\phi} \in \mathcal{V}(x_{0}) \rightarrow
V_{\phi}  \in \mathcal{V}(x_{0})$ such that:
\begin{enumerate}
\item[(a)] the map $(\phi, x) \mapsto \hat{\phi}(x)$ is uniformly continuous from $G \times X$ (with product uniformity)
to  $X$,
\item[(b)] for any $\phi, \psi \in G$ there is $D \in \mathcal{V}(x_{0})$
such that for any $x \in D$ $\hat{\phi \psi^{-1}}(x)$ and $\hat{\phi}(\hat{\psi}^{-1}(x))$ make sense and   $\hat{\phi \psi^{-1}}(x) = \hat{\phi}(\hat{\psi}^{-1}(x))$.
\end{enumerate}

\item[3.] Finally, a local group is an uniform space $G$ with an operation defined in a neighbourhood of $(e,e) \subset G \times G$ which satisfies the uniform group axioms locally.
\end{enumerate}
\label{dunifg}
\end{defi}
Remark that a local group acts locally at left (and also by conjugation) on itself.

This definition deserves an explanation. We shall start by recalling what an uniformity is (consult, for 
example the book James \cite{james}, for an introduction to uniform spaces). We need first some notations. 

A relation R of  a given set $X$ is just a subset $R \subset X \times X$. The inverse of R is 
$$R^{-1} \ = \ \left\{ (y,x) \in X \times X \mbox{ : } (x,y) \in R \right\}$$
The composition of two relations $R,S$ is defined by: 
$$R \circ S \ = \ \left\{ (x,z) \in X \times X \mbox{ : } \exists y \in X \ \ (x,y) \in R \ , \ (y,z) \in S \right\}$$
The composition is an associative operation. 

The diagonal of $X$ is the relation 
$$\Delta(X) \ = \ \left\{ (x,x) \mbox{ : } x \in x \right\} $$

A filter on $X$ is a set $\Omega$ of subsets of $X$ such that: 
\begin{enumerate} 
\item[-] if $A,B \in \Omega$ then $A \cap B \in \Omega$, 
\item[-] if $A  \in \Omega$, $B \subset X$ and $A \subset B$ then $B \in \Omega$.
\end{enumerate}

\begin{defi} An uniformity of a given set $X$ is a filter  $\Omega$ on $X \times X$ such that: 
\begin{enumerate}
\item[i)] for all $D \in \Omega$ $\Delta(X) \subset D$, 
\item[ii)] if $D \in \Omega$ then $D^{-1} \in \Omega$, 
\item[iii)] if $D \in \Omega$ then there is $E \in \Omega$ such that $E \circ E \subset D$.
\end{enumerate}
\end{defi}

For example, if $(X,d)$ is a metric space then the filter on $X \times X$ generated by the sets: 
$$B(r) \ = \ \left\{ (x,y) \in X \times X \mbox{ : } d(x,y)<r \right\}$$
for all $r > 0$ is an uniformity on $X$. Indeed, property i) comes from the fact that $d(x,x) = 0$ for 
any $x \in X$, ii) comes from $d(x,y)=d(y,x)$ for any $x,y \in X$ and iii) comes from the triangle inequality for the distance $d$. 

Let $(X,\Omega_{X})$ and $(Y,\Omega_{Y})$ be two uniform spaces. A function $\phi: X \rightarrow Y$ 
is uniformly continuos if 
$$\left(\phi \times \phi \right)^{-1}(\Omega_{Y}) \ \subset \ \Omega_{X}$$

An uniform group, according to the definition \eqref{dunifg}, is a group $G$ such that left translations are uniformly continuous functions and the left action of $G$ on itself is uniformly continuous too. 
In order two precisely formulate this we need two uniformities: one on $G$ and another on $G \times G$. 

These uniformities should be compatible, which is achieved by saying that $i'$, $i"$ are uniformly continuous. The uniformity of the group operation is achieved by saying that the $op$ morphism is uniformly continuous. 

The particular choice of the operation on $G \times G$ is not essential at this point, but it is 
justified by the cae of a Lie group endowed with the CC distance induced by a left invariant distribution. We shall construct a natural CC distance on $G \times G$, which is left invariant 
with respect to the chosen operation on $G \times G$. These distances induce uniformities which transform $G$ into an uniform group according to definition \eqref{dunifg}. The goal, achieved in 
the next section,  
is proposition \eqref{opsm} which claims that the operation function $op$ is derivable, even if 
right translations are not "smooth", i.e. commutative smooth according to definition 
\eqref{fdcd}.

We prepare now the path to this result. The "infinitesimal version" of an uniform group is a conical 
local uniform group.

\begin{defi}
A conical local uniform group $N$ is a local group with a local action of
$(0,+\infty)$ by morphisms $\delta_{\varepsilon}$ such that
$\displaystyle \lim_{\varepsilon \rightarrow 0} \delta_{\varepsilon} x \ = \ e$ for any
$x$ in a neighbourhood of the neutral element $e$.
\end{defi}

We shall make the following  hypotheses on the local uniform group $G$: there is a local action of $(0, +\infty)$ (denoted by
$\delta$), on $(G, e)$ such that
\begin{enumerate}
\item[H0.] the limit  $\lim_{\varepsilon \rightarrow 0} \delta_{\varepsilon} x \ = \ e$ exists and is uniform with respect to $x$.
\item[H1.] the limit
$$\beta(x,y) \ = \ \lim_{\varepsilon \rightarrow 0} \delta_{\varepsilon}^{-1}
\left((\delta_{\varepsilon}x) (\delta_{\varepsilon}y ) \right)$$
is well defined in a neighbourhood of $e$ and the limit is uniform.
\item[H2.] the following relation holds
$$ \lim_{\varepsilon \rightarrow 0} \delta_{\varepsilon}^{-1}
\left( ( \delta_{\varepsilon}x)^{-1}\right) \ = \ x^{-1}$$
where the limit from the left hand side exists in a neighbourhood of $e$ and is uniform with respect to $x$.
\end{enumerate}

\subsection{Virtual tangent bundle}

\begin{prop}
Under the hypotheses H0, H1, H2 $(G,\beta)$ is a conical local uniform group.
\end{prop}

\paragraph{Proof.}
All the uniformity assumptions permit to change at will the order of taking
limits. We shall not insist on this further and we shall concentrate on the
algebraic aspects.

We have to prove the associativity, existence of neutral element, existence of inverse and the property of being conical. The proof is straightforward.
For the associativity $\beta(x,\beta(y,z)) \ = \ \beta(\beta(x,y),z)$ we compute:
$$\beta(x,\beta(y,z)) \ = \ \lim_{\varepsilon \rightarrow 0 , \eta \rightarrow 0} \delta_{\varepsilon}^{-1} \left\{ (\delta_{\varepsilon}x) \delta_{\varepsilon/\eta}\left( (\delta_{\eta}y) (\delta_{\eta} z) \right) \right\}$$
We take $\varepsilon = \eta$ and we get
$$ \beta(x,\beta(y,z)) \ = \ \lim_{\varepsilon \rightarrow 0}\left\{
(\delta_{\varepsilon}x) (\delta_{\varepsilon} y) (\delta_{\varepsilon} z) \right\}$$
In the same way:
$$\beta(\beta(x,y),z) \ = \ \lim_{\varepsilon \rightarrow 0 , \eta \rightarrow 0} \delta_{\varepsilon}^{-1} \left\{ (\delta_{\varepsilon/\eta}x)\left( (\delta_{\eta}x) (\delta_{\eta} y) \right) (\delta_{\varepsilon} z) \right\}$$
and again taking $\varepsilon = \eta$ we obtain
$$\beta(\beta(x,y),z) \ = \  \lim_{\varepsilon \rightarrow 0}\left\{
(\delta_{\varepsilon}x) (\delta_{\varepsilon} y) (\delta_{\varepsilon} z) \right\}$$
The neutral element is $e$, from H0 (first part): $\beta(x,e) \ = \beta(e,x) \ = \ x$. The inverse of $x$ is $x^{-1}$, by a similar argument:
$$\beta(x, x^{-1})  \ = \ \lim_{\varepsilon \rightarrow 0 , \eta \rightarrow 0} \delta_{\varepsilon}^{-1} \left\{ (\delta_{\varepsilon}x)
\left( \delta_{\varepsilon/\eta}(\delta_{\eta}x)^{-1}\right) \right\}$$
and taking $\varepsilon = \eta$ we obtain
$$\beta(x, x^{-1})  \ = \ \lim_{\varepsilon \rightarrow 0}
\delta_{\varepsilon}^{-1} \left( (\delta_{\varepsilon}x) (\delta_{\varepsilon}x)^{-1}\right) \ = \ \lim_{\varepsilon \rightarrow 0} \delta_{\varepsilon}^{-1}(e) \ = \ e$$
Finally, $\beta$ has the property:
$$\beta(\delta_{\eta} x, \delta_{\eta}y) \ = \ \delta_{\eta} \beta(x,y)$$
which comes from the definition of $\beta$ and commutativity of multiplication
in $(0,+\infty)$. This proves that $(G,\beta)$ is conical.
\quad $\blacksquare$

We arrive at a natural realization of the tangent space to the neutral element.
Let us denote by $[f,g] \ = \ f \circ g \circ f^{-1} \circ g^{-1}$ the commutator of two transformations. For the group we shall denote by
$L_{x}^{G} y \ = \ xy$ the left translation and by $L^{N}_{x}y \ = \ \beta(x,y)$. The preceding proposition tells us that $(G,\beta)$ acts locally by left
translations on $G$. We shall call the left translations with respect to the group operation $\beta$ "infinitesimal". Those infinitesimal translations admit
the very important representation:
\begin{equation}
\lim_{\lambda \rightarrow 0} [L_{(\delta_{\lambda}x)^{-1}}^{G}, \delta_{\lambda}^{-1}] \ = \ L^{N}_{x}
\label{firstdef}
\end{equation}

\begin{defi}
The group $VT_{e}G$ formed by all transformations $L_{x}^{N}$ is called the virtual tangent space at $e$ to $G$.
\end{defi}

The virtual tangent space $VT_{x}G$ at $x \in G$ to $G$ is obtained by translating the group operation and the dilatations from $e$ to $x$. This means: define a new operation on $G$ by
$$y \stackrel{x}{\cdot} z \ = \ y x^{-1}z$$
The group $G$ with this operation is isomorphic to $G$ with old operation and
the left translation $L^{G}_{x}y \ = \ xy$ is the isomorphism. The neutral element is $x$.
Introduce also the dilatations based at $x$ by
$$\delta_{\varepsilon}^{x} y \ = \ x \delta_{\varepsilon}(x^{-1}y)$$
Then $G^{x} \ = \ (G,\stackrel{x}{\cdot})$ with the group of dilatations $\delta_{\varepsilon}^{x}$ satisfy the axioms Ho, H1, H2. Define then the virtual tangent
space $VT_{x}G$ to be: $VT_{x}G \ = \ VT_{x} G^{x}$. A short computation shows that
$$VT_{x} G \ = \ \left\{ L^{N,x}_{y} \ = \ L_{x} L^{N}_{x^{-1}y} L_{x} \mbox{ : } y \in U_{x} \in \mathcal{V}(X) \right\}$$
where
$$L^{N,x}_{y} \ = \ \lim_{ \lambda  \rightarrow 0} \delta_{\lambda}^{-1,x} [\delta_{\lambda}^{x}, L_{(\delta_{\lambda}x)^{x, -1}}^{G}] \delta_{\lambda}^{x}$$

We shall introduce the notion of commutative smoothness, which contains
a derivative resembling with Pansu derivative. This definition is a little bit stronger than the one given by Vodopyanov \& Greshnov \cite{vodopis2}, because  their definition is good for a general CC space, when uniformities are taken according to the distances in CC spaces $G^{(2)}$ and $G$.

\begin{defi}
A function $f: G_{1} \rightarrow G_{2}$ is commutative smooth at $x \in G_{1}$, where
$G_{1}, G_{2}$ are two groups satisfying H0, H1, H2,  if the application
$$u \in G_{1} \ \mapsto \ (f(x), Df(x)u) \in G_{2}^{(2)}$$
exists, where
$$Df(x)u \ = \ \lim_{\varepsilon \rightarrow 0} \delta_{\varepsilon}^{-1}
\left(f(x)^{-1}f(x \delta_{\varepsilon}u)\right)$$
and the convergence is uniform with respect to $u$ in compact sets.
\label{fdcd}
\end{defi}

For example the left translations $L_{x}$ are commutative smooth and
the derivative equals identity. If we want to see how the derivative moves
the virtual tangent spaces we have to give a definition.

Inspired by \eqref{firstdef},  we shall introduce the virtual tangent. We proceed as follows: to $f: G \rightarrow G$ and $x \in G$ let associate the function:
$$\hat{f}^{x}: G \times G \rightarrow G \ , \ \ \hat{f}^{x}(y,z) \ = \ \hat{f}^{x}_{y}(z) \ = \ \left(f(x)\right)^{-1}f(xy)z$$
To this function is associated a flow of left translations
$$ \lambda > 0 \ \mapsto \ \hat{f}^{x}_{\delta_{\lambda} y}: G \rightarrow G$$

\begin{defi}
The function $f: G \rightarrow G$ is virtually  derivable at $x \in G$ if there is
a virtual tangent $VDf(x)$ such that
\begin{equation}
\lim_{\lambda \rightarrow 0} \left[ \left(\hat{f}^{x}_{\delta_{\lambda} y} \right)^{-1} , \delta_{\lambda}^{-1} \right] \ = \ VDf(x)y
\label{2nddef}
\end{equation}
\label{vcd}
and the limit is uniform with respect to $y$ in a compact set.
\end{defi}

\begin{prop}
 Suppose that $f$ is commutative derivable. Then it is also virtually derivable and the virtual tangent to $f$ is given  by:
$$ VD f(x) y \ = \
L_{f(x)} L^{N}_{Df(x)y} L^{-1}_{f(x)}(y)$$
\end{prop}

\begin{rk} All definitions hold for $f: G_{1} \rightarrow G_{2}$ in an obvious way.
\end{rk}

With this definition $L_{x}$ is commutative smooth and it's virtual tangent in any point $y$ is a
group morphism from $VT_{y}G$ to $VT_{xy}G$. More generally, by standard reasoning, we get again the proposition \ref{ppansuprep}, in this more general setting. This time we shall give the proof, for further reference.

\begin{prop}
Let $f: G \rightarrow G$, $y, z \in G$ such that:
\begin{enumerate}
\item[(i)] $Df(x)y$, $Df(x)z$ exist for any $x \in G$.
\item[(ii)] The map $x \mapsto Df(x)z$ is continuous.
\item[(iii)] $x \mapsto
\delta_{t}^{-1} \left( f(x)^{-1} f(x \delta_{t}z)\right)$
converges uniformly to $Df(x)z$.
\end{enumerate}
Then for any (sufficiently small) $a,b > 0$ and any $w \ = \ \beta(\delta_{a}
 x, \delta_{b} y)$
the limit $Df(x)w$ exists and we have:
$$Df(x) \delta_{a} x \delta_{b} y \ = \ \beta( \delta_{a} Df(x)y ,
\delta_{b} Df(x) z)$$
\label{ppansuprepgen}
\end{prop}

\paragraph{Proof.}
The proof is standard. Remark that if $Df(x)y$ exists then for any
$a>0$ the limit $Df(x) \delta_{a} y$ exists and
$$Df(x) \delta_{a} y \ = \ \delta_{a} Df(x) y$$
It is not restrictive therefore to suppose that everything happens in a neighbourhood of the identity. Let $w= \beta(y,z)$.  We also use the notation
$$\beta_{t}(x,y) \ = \ \delta_{t}^{-1} \left( \delta_{t}(x) \delta_{t}(y) \right)$$

We write:
$$\delta_{t}^{-1} \left( f(x)^{-1} f(x \delta_{t}w) \right) \ = \
\beta_{t}\left( (1)_{t} , (2)_{t}\right)$$
with the notation (useful for further reference in proposition \eqref{pupu})
$$(1)_{t} \ = \ \delta_{t}^{-1} \left( f(x)^{-1} f(x \delta_{t}y)
\right)$$
$$(2)_{t} \ = \ \delta_{t}^{-1} \left( f(x \delta_{t}y)^{-1}
f(x \delta_{t}y \delta_{t} z) \right)$$
When $t$ tends to $0$ (i) implies that $(1)_{t}$ tends to $Df(x)y$,
(ii) and (iii) imply that $(2)_{t}$   tends to $Df(x) z$.
\quad $\blacksquare$

Remark that in principle the right translations are not commutative smooth.

Now that we have a model for the tangent space to $e$ at $G$, we can show that
the operation is commutative smooth.

\begin{prop}
Let $G$ satisfy H0, H1, H2 and $\delta_{\varepsilon}^{(2)} : G^{(2)} \rightarrow G^{(2)}$ be defined by
$$\delta_{\varepsilon}^{(2)} (x,u) \ = \ (\delta_{\varepsilon}x,
\delta_{\varepsilon} y)$$
Then $G^{(2)}$ satisfies H0, H1, H2, the operation ($op$ function) is commutative smooth  and we have the relation:
$$D \ op \ (x,u) (y,v) \ = \ \beta(y,v)$$
\label{opsm}
\end{prop}

\paragraph{Proof.}
It is sufficient to use the morphism property of the operation. Indeed, the right hand side of the relation to be proven is
$$RHS \ = \ \lim_{\varepsilon \rightarrow 0}
\delta_{\varepsilon}^{-1} \left( op(x,u)^{-1} op(x,u) op \left(\delta_{\varepsilon}^{(2)}(y,v)\right)\right) \ = $$
$$=  \ \lim_{\varepsilon \rightarrow 0}
\delta_{\varepsilon}^{-1} \left( op(\delta_{\varepsilon}^{(2)}(y,v))\right) \ = \ \beta(y,v)$$
The rest is trivial.
\quad $\blacksquare$

\subsection{The case of sub-Riemannian groups}

The notion of virtual tangent space is not based on the use
of distances, but on the use of dilatations. In fact, any manifold has a tangent space to any of its points, not only the Riemannian manifolds. We shall prove 
further that $VT_{e} G$ is isomorphic to the nilpotentisation
$N(G,D)$.

We start from Euclidean norm on $D$ and we choose an orthonormal basis of $D$. We can then extend the Euclidean norm to $\mathfrak{g}$ by
stating that the basis of $\mathfrak{g}$ constructed, as explained, from the basis on $D$, is orthonormal. By left translating the Euclidean norm on
$\mathfrak{g}$ we endow $G$ with a structure of Riemannian manifold. The
induced Riemannian  distance $d_{R}$ will give an uniform structure on $G$.
This distance is left invariant:
$$d_{R}(xy, xz) \ = \ d_{R}(y,z)$$
for any $x,y,z \in G$.

Any left invariant distance $d$ is uniquely determined if we set $d(x) \ = \ d(e,x)$.

The following lemma is important (compare with lemma \ref{p2.4}).

\begin{lema}
Let $X_{1}, ... , X_{p}$ be a basis of $D$.
Then there are $U \subset G$ and $V \subset N(G,D)$, open neighbourhoods of the neutral elements $e_{G}$, $e_{N}$ respectively,  and a surjective function $g: \left\{1, ... , M \right\} \rightarrow
\left\{ 1, ... ,p\right\}$ such that any $x \in U$, $y \in V$ can be written as
\begin{equation}
x \ = \ \prod_{i = 1}^{M} \exp_{G}(t_{i}X_{g(i)}) \ \ , \ y \ = \ \prod_{i = 1}^{M} \exp_{N}(\tau_{i}X_{g(i)})
\label{fpgen}
\end{equation}
\label{lgen}
\end{lema}

\paragraph{Proof.}
We shall make the proof for $G$; the proof for $N(G,D)$ will follow from the identifications explained before.

Denote by $n$ the dimension of $\mathfrak{g}$. We start the proof with the
remark that the function
\begin{equation}
(t_{1}, ... , t_{n})  \ \mapsto  \ \prod_{i = 1}^{n} \exp_{G}(t_{i}X_{i})
\label{tp1}
\end{equation}
is invertible in a neighbourhood of $0 \in R^{n}$, where the $X_{i}$ are elements of a basis $B$ constructed as before. Remember that each $X_{i} \in B$
is a multi-bracket of elements from the basis of $D$. If we replace a bracket
$\exp_{G} (t [x,y])$ in the expression \eqref{tp1} by $exp_{G} (t_{1}x) \ exp_{G} (t_{2}y) \ exp_{G} (t_{3}x)
 \ exp_{G} (t_{4}y)$ and we replace $t$ by $(t_{1}, ... , t_{4})$ then  the
image of a neighbourhood of $0$ by the obtained function still covers a
neighbourhood of the neutral element. We repeat this procedure a finite number of times and the thesis is proven.
\quad $\blacksquare$

As a corollary we obtain the Chow theorem for our particular example.

\begin{thm}
Any two points $x,y \in G$ can be joined by a horizontal curve.
\end{thm}

Let $d_{G}$ be the Carnot-Carath\'eodory distance induced by the distribution
$D$ and the metric. This distance is also left invariant. We obviously have $d_{R} \ \leq \ d_{G}$. We want to show that $d_{R}$ and $d_{G}$ induce the same uniformity on $G$.

Let us introduce another left invariant distance on $G$
$$d^{1}_{G}(x) \ = \ \inf \left\{ \sum \mid t_{i} \mid \mbox{ : } x \ = \
\prod \exp_{G} (t_{i} Y_{i}) \ , \ Y_{i} \in D \right\}$$
and the auxiliary functions :
$$\Delta^{1}_{G}(x) \ = \ \inf \left\{ \sum_{i=1}^{M} \mid t_{i} \mid \mbox{ : } x \ = \
\prod \exp_{G} (t_{i} X_{g(i)}) \right\}$$
$$\Delta^{\infty}_{G}(x) \ = \ \inf \left\{ \max \mid t_{i} \mid \mbox{ : } x \ = \
\prod \exp_{G} (t_{i} X_{g(i)}) \right\}$$
We can prove that $d^{1}_{G} = d_{G}(e_{G},\cdot)$. (Indeed $d_{G}(e_{G},\cdot) \leq d^{1}_{G}$. On the other part $d^{1}_{G}(x)$ is less equal than the variation of any Lipschitz curve joining $e_{G}$ with $x$. Therefore we have equality.)

The functions $\Delta^{1}_{G}$, $\Delta^{\infty}_{G}$ don't induce left invariant distances. Nevertheless they are useful, because of their equivalence:
\begin{equation}
\Delta^{\infty}_{G}(x) \ \leq \ \Delta^{1}_{G}(x) \ \leq \ M \ \Delta_{G}^{\infty}(x)
\label{dequi}
\end{equation}
for any $x \in G$. This is a consequence of the lemma
\ref{lgen}.

We have therefore the chain of inequalities:
$$d_{R} \ \leq \ d_{G} \ \leq \ \Delta^{1}_{G} \ \leq \ M \Delta^{\infty}_{G}$$
But from the proof of lemma \ref{lgen} we see that $\Delta^{\infty}_{G}$ is uniformly continuous. This proves the equivalence of the uniformities.

Because $\exp_{G}$ does not deform much the Riemannian $d_{R}$ distances near $e$, it can be checked  that the group $G$ with the dilatations
$$\tilde{\delta}_{\varepsilon} (\exp_{G} x) \ = \ \exp_{G}(\delta_{\varepsilon}
x)$$
satisfies H0, H1, H2.

The same conclusion is true for the local uniform group (with the uniformity induced by the Euclidean distance) $\mathfrak{g}$ with the operation:
$$X \opg Y \ = \ \log_{G} \left( \exp_{G} (X) \exp_{G} (Y) \right)$$
for any $X,Y$ in a neighbourhood of $0 \in \mathfrak{g}$. Here the dilatations
are $\delta_{\varepsilon}$. We shall denote this group by $\log G$. These two groups are isomorphic as local uniform groups by the map $\exp_{G}$. Dilatations commute with the isomorphism. They have therefore isomorphic (by $\exp_{G}$) virtual tangent spaces.

\begin{thm}
The virtual tangent space $VT_{e}G$ is isomorphic to $N(G,D)$. More precisely
$N(G,D)$ is equal (as local group) to the virtual tangent space to $\log G$:
$$N(G,D) \ = \ VT_{0} \log G$$
\label{teore}
\end{thm}

\paragraph{Proof.}
The product $XY$ in $\log G$ is given by  Baker-Campbell-Hausdorff formula
$$X \opg Y \ = \ X + Y + \frac{1}{2} [X,Y] + ... $$
Use proposition \ref{limbra} to compute $\beta(X,Y)$ and show that
$\beta(X,Y)$ equals the nilpotent multiplication.
\quad $\blacksquare$

\subsection{Proof of Mitchell's theorem 1}
\label{metriccone}

We shall prove the result of Mitchell \cite{mit} theorem 1, that $G$ admits in any point a metric tangent cone, which is isomorphic with the nilpotentisation
$N(G,D)$. Mitchell theorem is true for regular sub-Riemannian manifolds.
The proof that we give here is based on the lemma \ref{lgen} and Gromov
\cite{gromov} section 1.2.

Because left translations are isometries, it is sufficient to prove that $G$
admits a metric tangent cone in identity and that the tangent cone is isometric
with $N(G,D)$.

For this we transport all in the algebra $\mathfrak{g}$, endowed with two
brackets $[\cdot, \cdot]_{G}$, $[\cdot, \cdot ]_{N}$ and with two operations
$\opg$ and $\opn$. Denote by $d_{G}$, $d_{N}$ the Carnot-Carath\'eodory distances corresponding to the $\opg$, respectively $\opn$ left invariant distributions on (a neighbourhood of $0$ in) $\mathfrak{g}$. $l_{G}$,
$l_{N}$ are the corresponding length functionals. We shall denote by $B_{G}(x,R)$, $B_{N}(x,R)$ the balls centered in $x$ with
radius $R$ with respect to $d_{G}$, $d_{N}$.

We can refine lemma \ref{lgen} in order to obtain the Ball-Box theorem in this
more general situation.

\begin{thm}
(Ball-Box Theorem) Denote by
$$Box^{1}_{G} (\varepsilon) \ = \ \left\{ x \in G \mbox{ : } \Delta^{1}_{G}(x) \ < \ \varepsilon \right\} $$
$$Box^{1}_{N} (\varepsilon) \ = \ \left\{ x \in N \mbox{ : } \Delta^{1}_{N}(x) \ < \ \varepsilon \right\} $$
For small $\varepsilon > 0$ there is a constant $C>1$ such that
$$\exp_{G} \left( Box^{1}_{N} (\varepsilon) \right) \ \subset \
Box^{1}_{G} (C \varepsilon) \ \subset \
\exp_{G} \left( Box^{1}_{N} (C^{2} \varepsilon) \right)$$
\label{bbt}
\end{thm}

\paragraph{Proof.}
Reconsider the proof of lemma \ref{lgen}. This time, instead of the trick of
replacing commutators $\exp_{G}[X,Y](t)$ with four letters words
$$\exp_{G}X(t_{1}) \exp_{G} Y (t_{2}) \exp_{G}X (t_{3}) \exp_{G} Y (t_{4})$$
we shall use a smarter replacement (which, important! , works equally for
the nilpotentisation $N(G,D)$).

We start with a basis $\left\{ Y_{1}, ... , Y_{n} \right\}$ of the algebra
$\mathfrak{g}$, constructed from multibrackets of elements $\left\{ X_{1}, ... , X_{p} \right\}$ which form a basis for the distribution $D$. Introduce an Euclidean distance by declaring the basis of $\mathfrak{g}$ orthonormal. Set
$$[X(t), Y(t)]^{\circ} \ = \ (tX) \opg (tY) \opg (-tX) \opg (-tY)$$
It is known that for any set of vectors $Z_{1}, ... , Z_{q} \in \mathfrak{g}$,
if we denote by $\alpha^{\circ}(Z_{1}(t), ... , Z_{q}(t))$ a $[ \cdot , \cdot]^{\circ}$ multibracket and by $\alpha(Z_{1}, ... , Z_{q})$ the same constructed
$[\cdot , \cdot ]$ multibracket, then we have
$$ \alpha^{\circ}(Z_{1}(t), ... , Z_{q}(t)) \ = \ t^{q} \alpha(Z_{1}, ... , Z_{q}) \ + \ o(t^{q})$$
with respect to the Euclidean norm.

Remark as previously that   the function
\begin{equation}
(t_{1}, ... , t_{n})  \ \mapsto  \ \prod_{i = 1}^{n} (t_{i}Y_{i})
\label{tpu1}
\end{equation}
is invertible in a neighbourhood of $0 \in R^{n}$. Each $X_{i}$ from the
basis of $\mathfrak{g}$ can be written as a multibracket
$$X_{i} \ = \ \alpha_{i}( X_{k_{1}}, ... , X_{k_{j_{i}}})$$
which has the length $l_{i} = j_{i} - 1$. If $l_{i}$ is odd then replace
$(t_{i}X_{i})$ by
$$\alpha^{\circ}( X_{k_{1}}(t^{1/l_{i}}), ... , X_{k_{j_{i}}}(t^{1/l_{i}}))$$
If $l_{i}$ is even then the multibracket $\alpha_{i}$ can be rewritten as
$$\alpha_{i}( X_{k_{1}}, ... , X_{k_{j_{i}}}) \ = \ \beta_{i}(X_{k_{1}} , ... ,
[X_{k_{p}}, X_{k_{p+1}}], ... , X_{k_{j_{i}}})$$
Replace then $(t_{i}X_{i})$ by
$$\left\{ \begin{array}{ll}
\beta_{i}^{\circ}(X_{k_{1}}(\mid t \mid^{1/l_{i}}) , ... ,
[X_{k_{p}}(\mid t \mid^{1/l_{i}}), X_{k_{p+1}}(\mid t \mid^{1/l_{i}})]^{\circ}, ... , X_{k_{j_{i}}}(\mid t \mid^{1/l_{i}})) & \mbox{ if } t \geq 0 \\
\beta_{i}^{\circ}(X_{k_{1}}(\mid t \mid^{1/l_{i}}) , ... ,
[X_{k_{p+1}}(\mid t \mid^{1/l_{i}}), X_{k_{p}}(\mid t \mid^{1/l_{i}})]^{\circ}, ... , X_{k_{j_{i}}}(\mid t \mid^{1/l_{i}})) & \mbox{ if } t \leq 0
\end{array} \right.$$
After this replacements in the  expression \eqref{tpu1} one obtains a function
$E_{G}$ which is still invertible in a neighbourhood of $0$. We obtain a function $E_{N}$ with the same algebraic expression as $E_{G}$, but with
$[\cdot , \cdot ]^{\circ}_{N}$ brackets instead of $[\cdot , \cdot ]^{\circ}_{G}$ ones. Use these functions to (obviously) end the proof of the theorem.
\quad $\blacksquare$

\begin{thm}(Mitchell's theorem 1)
The Gromov-Hausdorff limit of pointed metric spaces $(\mathfrak{g}, 0, \lambda
d_{G})$ as $\lambda \rightarrow \infty$ exists and equals $(\mathfrak{g}, d_{N})$.
\end{thm}

\paragraph{Proof.}
We shall use the proposition \ref{pbur}. For this we shall construct  $\varepsilon$ isometries between $Box^{1}_{N}(1)$ and $Box^{1}_{G}(C \varepsilon)$. These are provided by  the function $E_{G} \circ
E_{N}^{-1} \circ \delta_{\varepsilon}$.

The trick consists in the definition of the nets. We shall exemplify the construction for the case of a 3 dimensional algebra $\mathfrak{g}$. The basis
of $\mathfrak{g}$ is $X_{1}, X_{2}, X_{3} \ = \ [X_{1}, X_{2}]_{G}$. Divide
the interval $[0,X_{1}]$ into $P$ equal parts, same for the interval $[0,X_{2}]$. The interval $[0,X_{3}]$ though will be divided into $P^{2}$ intervals. The net so obtained, seen in $N \ \equiv \ \mathfrak{g}$, turn $E_{G} \circ
E_{N}^{-1} \circ \delta_{\varepsilon}$ into a $\varepsilon$ isometry between $Box^{1}_{N}(1)$ and $Box^{1}_{G}(C \varepsilon)$. To check this is mostly a matter of smart (but still heavy) notations.
\quad $\blacksquare$

The proof of this theorem is basically a refinement of Proposition 3.15, Gromov \cite{gromov}, pages 85--86, mentioned in the introduction of these notes.
Close examination shows that the theorem is a consequence of the following facts:
\begin{enumerate}
\item[(a)] the identity map $id \ : (\mathfrak{g}, d_{G}) \rightarrow
(\mathfrak{g}, d_{N})$ has finite dilatation in $0$ (equivalently
$\exp_{G}: N(D,G) \rightarrow G$ has finite dilatation in $0$),
\item[(b)] the change from the $G$-invariant distribution induced by $D$ to the
$N$ invariant distribution induced by the same $D$ is classically smooth.
\end{enumerate}
It is therefore natural that the result does not feel the non-derivability of
$id$ map in $x \not = 0$. Here one has to think of $id$ map as the exponential map, defined from the group algebra $\mathfrak{g}$, endowed with the nilpotent
distribution, to the group $G$, endowed with initial distribution.

\subsection{Constructions of tangent bundles}

In the previous sections we had fixed one parameter group of dilatations which
satisfy the axioms $H0$, $H1$, $H2$. We shall look now to the class of all
dilatation groups (or dilatation flows, because these groups are one dimensional)
which satisfy those axioms.

$G$ is a real connected Lie group endowed with a left invariant distribution
$D$, as considered previously. We shall denote by $d$ the Carnot-Carath\'eodory distance and by $\delta_{\varepsilon}$ the one parameter group of dilatations
(which are morphisms of the nilpotentisation, but they are transported first
on the Lie algebra $\mathfrak{g}$ and next on a neighbourhood of the
neutral element in $G$). We know now that $G$ endowed with the dilatations
$\delta_{\varepsilon}$ and with the uniformities
\begin{enumerate}
\item[-] the $d$ induced uniformity on $G$,
\item[-] the uniformity induced by the Carnot-Carath\'eodory distance
on $G^{(2)}$. This distance correspond to the left-invariant distribution generated by $D \times D \subset \ Lie (G^{(2)})$,
\end{enumerate}
is an uniform group according to the definition \ref{dunifg}.

The use of the distance $d$ simplifies the writing of the paper. The constructions described further work for uniform groups.

We begin by considering the class of all dilatation  flows $\varepsilon \mapsto \hat{\delta}_{\varepsilon}
: G \rightarrow G$ such that:
\begin{enumerate}
\item[(a)] $\displaystyle\lim_{\varepsilon \rightarrow 0} \hat{\delta}_{\varepsilon}(x) \ = \ e$ uniformly with respect to
$x \in K$, $K$ compact neighbourhood of $e$.
\item[(b)] $\displaystyle\beta(x,y) \ = \ \lim_{\varepsilon \rightarrow 0}
\hat{\delta}_{\varepsilon}^{-1} \left( \hat{\delta}_{\varepsilon}(x) \hat{\delta}_{\varepsilon}(y) \right)$ uniformly with respect to
$x, y \in K$, $K$ compact neighbourhood of $e$.
\item[(c)] $\displaystyle x^{-1} \ = \ \lim_{\varepsilon \rightarrow 0}
\hat{\delta}_{\varepsilon}^{-1} \left( \hat{\delta}_{\varepsilon}(x)^{-1}\right)$
uniformly with respect to
$x \in K$, $K$ compact neighbourhood of $e$.
\end{enumerate}
If we replace $\delta_{\varepsilon}$ by $\hat{\delta}_{\varepsilon}$ then we obtain the same infinitesimal objects (same virtual tangent space, same virtual tangent bundle).

\begin{defi}
\begin{enumerate}
\item[(i)] Two dilatation flows $\delta'$, $\delta"$ are equivalent (we write $\delta' \equiv \delta"$) if
$$ \lim_{\varepsilon \rightarrow 0} \frac{1}{\varepsilon} d(\delta'_{\varepsilon} \circ \delta_{\varepsilon}^{-1}(x), \delta'_{\varepsilon}
\circ \delta_{\varepsilon}^{-1}(x)) \ = \ 0$$
uniformly with respect to $x$.
We denote by $[\delta']$ the equivalence class of $\delta'$.
\item[(ii)] $\mathcal{A}$   is the class
of all dilatation flows which are equivalent to flows of the following
form: $$\hat{\delta}_{\varepsilon}(x) \ = \ y_{\varepsilon}^{-1} \delta_{\varepsilon}(x) \ \ , \ \lim_{\varepsilon \rightarrow 0} \delta_{\varepsilon}^{-1} y_{\varepsilon} \ = \ 0$$
\end{enumerate}
\end{defi}
It is easy to check that any dilatation flow $\delta'_{\varepsilon} \in \mathcal{A}$ satisfies conditions  (a), (b), (c).

Further we shall look to "finite curves" which are  kind of "rectifiable" curves
(in Margulis, Mostow terminology) with respect to a dilatation flow in the class
$\mathcal{A}$.

\begin{defi}
The class of finite curves is the collection of all curves $c: [a,+\infty) \rightarrow G$ such that $\displaystyle\delta_{\varepsilon} c(\frac{1}{\varepsilon})$ converges to a point $x_{c}$, as $\varepsilon$ goes to $0$.

On the class of finite curves we define two equivalence relations (weak and strong):
\begin{enumerate}
\item[-] the curves $c_{1}$, $c_{2}$ are weakly equivalent
$$c_{1} \weak c_{2} \ \Leftrightarrow \ \lim_{\varepsilon \rightarrow 0}
\frac{1}{\varepsilon} d( \delta_{\varepsilon} c_{1}(\frac{1}{\varepsilon}) , \delta_{\varepsilon}c_{2}(\frac{1}{\varepsilon})) \ = \ 0$$
\item[-] the curves $c_{1}$, $c_{2}$ are strongly equivalent
$$c_{1} \strong c_{2} \ \Leftrightarrow \ \lim_{\varepsilon \rightarrow 0}
\frac{1}{\varepsilon} d(\  \delta_{\varepsilon}c_{1}(\frac{1}{\varepsilon}) \  x , \delta_{\varepsilon}c_{2}(\frac{1}{\varepsilon}) \ x) \ = \ 0$$
uniform with respect to $x$ (in a neighbourhood of the neutral element).
\end{enumerate}
\end{defi}
It is an  obvious remark that $c_{1} \strong c_{2}$ implies $c_{1} \weak c_{2}$. If we denote by $[c]_{w}$, $[c]_{s}$, the weak, strong respectively, equivalence class of $c$, then $[c]_{w}$ decomposes in disjoint union of strong classes.

\begin{prop}Let $c_{1}, c_{2}$ be finite curves. Define
$$\delta_{(c_{1},c_{2})_{\varepsilon}} \ = \ \delta_{\varepsilon}c_{1}(\frac{1}{\varepsilon})^{-1} \  \delta_{\varepsilon} c_{2}(\frac{1}{\varepsilon}) \
\delta_{\varepsilon}$$
We have then:
\begin{enumerate}
\item[(a)] $c_{1} \weak c_{2}$ is equivalent to $\displaystyle \delta_{(c_{1},c_{2})} \in \mathcal{A}$.
\item[(b)] $c_{1} \strong c_{2}$ is equivalent to $\displaystyle \delta_{(c_{1},c_{2})} \equiv \delta$.
\item[(c)] Let $\Sigma [c]$ be the spectrum of the finite curve $c$, that is the collection of all $\displaystyle [\delta_{(c,c')}]$, with $c' \weak c$. Then $$[c]_{w} \ = \ \bigcup_{[\delta_{(c,c')}] \in \Sigma [c] } [c']_{s}$$
(disjoint union). In other words the weak class of $c$ decomposes into strong
classes enumerated by the spectrum of $c$.
\end{enumerate}
\label{pro1}
\end{prop}

\paragraph{Proof.}
The points (a), (b) by direct check of definitions. Point (c) follows from
(a), (b). \quad $\blacksquare$

The definition of spectrum in proposition \eqref{pro1} (c) deserves attention. We may choose first a class of curves $\mathcal{C}$. Then weak, strong classes and spectra are defined with respect to the 
class $\mathcal{C}$. For the moment the class $\mathcal{C}$ is the class of finite curves, but in the future we shall restrict this class.

The composition of curves is defined by:
$$c_{1} \circ c_{2}(\frac{1}{\varepsilon}) \ = \ \delta_{\varepsilon}^{-1} \left( \delta_{\varepsilon}c_{1}(\frac{1}{\varepsilon}) \
 \delta_{\varepsilon}c_{2}(\frac{1}{\varepsilon})\right)$$

The following proposition justifies the introduction of strong equivalence classes.

\begin{prop}
 The composition operation induces an operation on strong equivalence
classes.

Moreover, if $c_{1}\strong c_{1}'$ and $c_{2} \weak c_{2}'$ then
$c_{1} \circ c_{2} \weak c_{1}' \circ c_{2}'$.

However, composition does not induce an operation on weak classes, generally.
\label{pro2}
\end{prop}

\paragraph{Proof.}
We shall prove the second assertion: if $c_{1} \strong c_{1}'$ and $c_{2} \weak c_{2}'$ then $c_{1} \circ c_{2} \weak c_{1}' \circ c_{2}'$. The first and third
parts of the proposition have a similar proof.

Let $\varepsilon > 0$. We have then
$$\frac{1}{\varepsilon} d\left( \delta_{\varepsilon}^{-1} (c_{1} \circ c_{2})
\left(\frac{1}{\varepsilon}\right) , \delta_{\varepsilon}^{-1} (c_{1}' \circ c_{2}')
\left(\frac{1}{\varepsilon}\right) \right) =  $$
$$ = \frac{1}{\varepsilon} d\left( \delta_{\varepsilon} c_{1} \left(\frac{1}{\varepsilon}\right) \delta_{\varepsilon} c_{2}
\left(\frac{1}{\varepsilon}\right) , \delta_{\varepsilon} c_{1}' \left(\frac{1}{\varepsilon}\right) \delta_{\varepsilon} c_{2}'
\left(\frac{1}{\varepsilon}\right) \right) \leq $$
$$ \leq \frac{1}{\varepsilon} d\left( \delta_{\varepsilon} c_{2}'
\left(\frac{1}{\varepsilon}\right) , \delta_{\varepsilon} c_{1}' \left(\frac{1}{\varepsilon}\right)\right) + $$
$$ + \leq \frac{1}{\varepsilon} d\left( \delta_{\varepsilon} c_{1}' \left(\frac{1}{\varepsilon}\right) \delta_{\varepsilon} c_{2}
\left(\frac{1}{\varepsilon}\right), \delta_{\varepsilon} c_{1} \left(\frac{1}{\varepsilon}\right) \delta_{\varepsilon} c_{2}
\left(\frac{1}{\varepsilon}\right)\right)$$
But both quantities from the RHS tend to $0$ as $\varepsilon \rightarrow 0$, because of the hypothesis.
\quad $\blacksquare$

On $\mathcal{A}/\equiv$ we have the following operation:
$$[y_{\varepsilon}^{-1} \delta_{\varepsilon}] \circ [u_{\varepsilon}^{-1} \delta_{\varepsilon}]   \ = \ [y_{\varepsilon}^{-1} u_{\varepsilon}^{-1} \delta_{\varepsilon}]$$
(Check that the operation is well defined). This operation makes
$\mathcal{A}/\equiv$ into a group. This group admits an one parameter
family of isomorphisms:
$$\underline{\delta}_{\mu} [[y_{\varepsilon}^{-1} \delta_{\varepsilon}] \ = \
[y_{\varepsilon \mu}^{-1} \delta_{\varepsilon \mu}]$$
Therefore any Lie subgroup of $\mathcal{A}/\equiv$ is Carnot.
A particular subgroup of $\mathcal{A}/\equiv$ is the spectrum of a curve.

Another observation is that a weak class can be seen as a bundle. The bijection between fibers (strong classes) is as follows: take in the strong class
$[c_{1}]_{s}$ a representant $c_{1}\weak c$. Then the function:
$$ c' \in [c]_{s} \ \mapsto c_{1} \circ c^{-1} \circ c' \ \in [c_{1}]_{s}$$
is a bijection.

However,
$\mathcal{A}/\equiv$, or the spectrum of a curve are too big to be finite dimensional groups. That is why we shall restrict to smaller classes.

A particular class of finite curves is given by words.
A word is a string $(x_{1} y_{1} .... x_{p} y_{p})$ with $x_{i}, y_{i} \in G$.
To any word it is associated the curve:
$$c \ = \ \left( \begin{array}{cccccc}
x_{1} & & x_{2} & ... & x_{p} & \\
 & y_{1} &  & y_{2} & ... & y_{p}
\end{array} \right) \ \ , \ c(\frac{1}{\varepsilon}) \ = \ \delta_{\varepsilon}^{-1} \left(  x_{1} \ \delta_{\varepsilon} y_{1} \ ... x_{p} \delta_{\varepsilon} y_{p} \right) $$
Word concatenation correspond to  curves composition.

For any $\mu > 0$ we define $\displaystyle \left(\delta_{\mu} \circ c \right)(\frac{1}{\varepsilon}) \ = \ \delta_{\mu} c(\frac{1}{\varepsilon \mu})$.

\begin{prop}
We have the following relations:
\begin{enumerate}
\item[(a)] $\displaystyle \left( \begin{array}{cccccc}
x_{1} & & x_{2} & ... & x_{p} & \\
 & y_{1} &  & y_{2} & ... & y_{p}
\end{array} \right) \circ \left( \begin{array}{cc}
0 & \\
  & y
\end{array} \right) \ \weak \ \left( \begin{array}{cccccc}
x_{1} & & x_{2} & ... & x_{p} & \\
 & y_{1} &  & y_{2} & ... & \beta (y_{p}, y)
\end{array} \right)$
\item[(b)] $\displaystyle \left( \begin{array}{cc}
x & \\
  & 0
\end{array} \right) \circ \left( \begin{array}{cccccc}
x_{1} & & x_{2} & ... & x_{p} & \\
 & y_{1} &  & y_{2} & ... & y_{p}
\end{array} \right) \ = \ \left( \begin{array}{cccccc}
x \ x_{1} & & x_{2} & ... & x_{p} & \\
 & y_{1} &  & y_{2} & ... & y_{p}
\end{array} \right)$
\item[(c)] $\displaystyle \delta_{\mu} \circ \left( \begin{array}{cccccc}
x_{1} & & x_{2} & ... & x_{p} & \\
 & y_{1} &  & y_{2} & ... & y_{p}
\end{array} \right) \ = \ \left( \begin{array}{cccccc}
x_{1} & & x_{2} & ... & x_{p} & \\
 & \delta_{\mu} y_{1} &  & \delta_{\mu} y_{2} & ... & \delta_{\mu} y_{p}
\end{array} \right)$
\end{enumerate}
\label{pro3}
\end{prop}
This proposition gives an interesting interpretation of the Margulis \& Mostow tangent bundle. 

\begin{cor}
The virtual tangent bundle (or the Margulis \& Mostow
tangent bundle) is  the bundle  of two-letters words, with basis
the group and bundle projection on the first letter.
\end{cor}

The definition of the tangent bundle can be naturally given in at least two ways.
Here is the first.

\begin{defi}
The word tangent bundle $WG$ is the group of strong equivalence classes of words.
The bundle structure of $WG$ is given by the map
$$[c]_{s} \mapsto x_{c} \ = \ \lim_{\varepsilon \rightarrow 0} \delta_{\varepsilon} c\left(\frac{1}{\varepsilon}\right)$$
$W_{0}G$ is the fiber of the neutral element and it is a group.
\label{deftg}
\end{defi}

We shall explain further the structure of $WG$. For this we shall start from  the class $\mathcal{C}^{m}([0,+\infty), \mathfrak{g})$ of all functions in $\mathfrak{g}$ which admit a development of order $m$ around $0$. To any
$f \in \mathcal{C}^{m}([0,+\infty), \mathfrak{g})$ we associate the finite curve
$$c_{f}(\frac{1}{\varepsilon}) \ = \ \delta_{\varepsilon}^{-1}
\left( f(\varepsilon) \right)$$
Conversely, to any finite curve $c$ we associate the function
$$f_{c}(\varepsilon) \ = \ \delta_{\varepsilon} \left( c\left( \frac{1}{\varepsilon} \right) \right)$$
This function is prolongated by continuity to $\varepsilon=0$ (which is possible, because the curve is finite).

For any  curve $c$ associated to a word we have $f_{c} \in \mathcal{C}^{m}([0,+\infty), \mathfrak{g})$. We shall use further the following notations:
$$C(G) \ = \ \left\{ c_{f} \mbox{ : } f \in \mathcal{C}^{m}([0,+\infty), \mathfrak{g}) \right\} / \strong$$
$$WG \ = \ \left\{ c_{f} \mbox{ : } f \mbox{ is a word } \right\} / \strong$$
 The word tangent bundle $WG$ can be seen as a subset of $C(G)$.

The tangent bundle $WG$ is a group with the multiplication operation of strong
classes of curves. We want to show that in fact is a Lie group.

Let $f \in \mathcal{C}^{m}([0,+\infty), \mathfrak{g})$ and
$$f(\varepsilon) \ = \ P(\varepsilon) + S(\varepsilon)$$
with $P \in \mathfrak{g}[\varepsilon]$ is a polynomial of degree at most
$m$ and $S$ is of order $o(\varepsilon^{m})$. We want to check if
$c_{f}$ and $c_{P}$ are strongly equivalent. This is equivalent with $c_{P}^{-1}
\circ c_{f} \strong 0$, which is straightforward.

We conclude therefore that
$$C(G) \ = \ \left( \mathfrak{g}[\varepsilon] / \varepsilon^{m+1} \mathfrak{g}[\varepsilon] \right) / \strong$$
in the sense that the sides of the previous equality can be identified. Let us denote
$$\mathfrak{g}^{(m)}[\varepsilon] \ = \ \mathfrak{g}[\varepsilon] / \varepsilon^{m+1} \mathfrak{g}[\varepsilon]$$
This is a factorisation of a Lie algebra by an ideal, hence it is itself a Lie algebra.

Recall that $\mathfrak{g}$ admits a direct sum decomposition and a filtration:
$$\mathfrak{g} \ = \ V_{1} + ... + V_{m} \ , \ \ V^{i} \ = \ V_{1}+ ... + V_{i}$$
such that $[V^{i}, V^{j}] \subset V^{i+j}$ (here $V^{k} = \mathfrak{g}$ for any
$k \geq m$).

 Any element of $\mathfrak{g}^{(m)}[\varepsilon]$ has a representant
$P \in \mathfrak{g}[\varepsilon]$, $\deg P \leq m$,
$$P(\varepsilon) \ = \ P_{1}(\varepsilon)+...+P_{m}(\varepsilon)$$
such that $P_{i} \in V_{i}[\varepsilon]$  for all $i$. From the degree condition
we see that we can write (for any $i= 1,...,m$):
$$P_{i}(\varepsilon) \ = \ \sum_{k=0}^{m} a_{ik} \varepsilon^{k} \ , \ \ a_{ik} \in V_{i}$$
We can do the following decomposition:
$$P(\varepsilon) \ = \ P^{W}(\varepsilon) + P^{S}(\varepsilon)$$
$$P^{W}(\varepsilon)_{i} \ = \ \sum_{k=0}^{i} a_{ik} \varepsilon^{k}$$
Let $c^{1}$, $c^{2}$ two arbitrary finite curves and $f_{i} = f_{c^{i}}$, $i=1,2$.
Then $c^{1} \weak c^{2}$ if and only if
$$\delta_{\varepsilon}^{-1} \left( f_{1}^{-1}(\varepsilon) f_{2}(\varepsilon) \right) \rightarrow 0$$
as $\varepsilon \rightarrow 0$. From here (and the filtration of $\mathfrak{g}$) we deduce that $c_{P} \weak c_{P^{W}}$, therefore the weak class of $f$ is
$P^{W}$, where $P$ is the development of order $m$ of $f$ around $0$.

We have to identify now the strong classes which compose a weak class. For this
consider the sets:
$$Weak(\mathfrak{g}) \ = \ \left\{ P^{W} \mbox{ : } P \in \mathfrak{g}[\varepsilon] \ , \ \ \deg P \leq m \right\}$$
$$Strong(\mathfrak{g}) \ = \ \left\{ P^{S} \mbox{ : } P \in \mathfrak{g}[\varepsilon] \ , \ \ \deg P \leq m \right\}$$
$$O(\mathfrak{g}) \ = \ \left\{ P \in Strong(\mathfrak{g}) \mbox{ : }
c_{P} \strong 0 \right\}$$
The set $O(\mathfrak{g})$ can be seen as a normal  subgroup of $\mathfrak{g}^{(m)}[\varepsilon]$. Indeed, $P \in O(\mathfrak{g})$ if and only if for
any $x \in G$ we have
$$\delta_{\varepsilon}^{-1} Ad_{x}^{G} P(\varepsilon) \rightarrow 0$$
as $\varepsilon \rightarrow 0$. This implies that $O(\mathfrak{g})$ is
$Ad^{G}$ invariant, hence the claim.

Because of proposition \ref{pro2} it follows that if $P \in
 O(\mathfrak{g})$ and $Q \in Strong(\mathfrak{g})$ then $PQ \in Strong(\mathfrak{g})$. Putting all together, we have obtained the characterisation of $C(G)$.

\begin{thm}
$C(G)$ is a Lie group with Lie algebra isomorphic to
$$Weak(\mathfrak{g}) + \left( Strong(\mathfrak{g})/ O(\mathfrak{g}) \right)$$
For any $f \in \mathcal{C}^{m}([0,+\infty),G)$ with associated $P$, development of
order $m$ around $0$, the weak class $[c_{f}]_{w}$ can be identified with
$P^{W}$. The spectrum of $c_{f}$ (in the same class of curves as $c_{f}$) does not
depend on $f$ and equals $Strong(\mathfrak{g})/ O(\mathfrak{g})$.
\label{tstruc1}
\end{thm}

$WG$ is a closed subgroup of $C(G)$, hence it is a Lie group. $W_{0}G$ is a
subgroup of $WG$.

It will be useful further to work with $C(G)$. We shall denote by $T_{0} G$ the fiber of $0$ in $C(G)$, namely the set:
$$T_{0} G \ = \ \left\{ f \in \mathcal{C}^{m}([0,+\infty),G) \mbox{ : } f(0) = 0 \right\} / \strong$$

The notion which later will lead to derivability is transport.

\begin{defi}
The continuous function $f: N \rightarrow M$ transports finite curves into finite curves according to the definition:
$$f*c (\frac{1}{\varepsilon}) \ = \ \delta_{\varepsilon}^{-1} f( \delta_{\varepsilon} c(\frac{1}{\varepsilon}))$$
\end{defi}

Let us  see how the translations and dilatations act on the tangent bundle.

We start from the observation: translations and dilatations preserve strong classes. Indeed, let $f: G \rightarrow G$ be a left or right translation and $c_{f_{1}} \strong c_{f_{2}}$, with $f_{1}, f_{2} \in \mathcal{C}^{m}([0,+\infty),G)$. We want to prove that
$f * c_{1} \strong f * c_{2}$. Suppose first that $f(x) = y x$. The claim is then true because the distance (or uniformity) in $G$ is left invariant. If $f$ is a right translation ($f(x) = xy$), take $z \in G$ arbitrary. Then
$$\frac{1}{\varepsilon} d(\delta_{\varepsilon} \left( f * c_{1} (\frac{1}{\varepsilon}) \right) z, \delta_{\varepsilon} \left( f * c_{2} (\frac{1}{\varepsilon}) \right)z) \ = \ \frac{1}{\varepsilon} d(f_{1}(\varepsilon)yz, f_{2}(\varepsilon)yz)$$
which converges to $0$ as $\varepsilon \rightarrow 0$ because $c_{1} \strong
c_{2}$.

Recall that we defined a dilatation on the class of curves (and we used the concatenation notation)
$$\delta_{\mu} (c) (\frac{1}{\varepsilon}) = (\delta_{\mu} \circ c) (\frac{1}{\varepsilon}) = \delta_{\mu} c (\frac{1}{\varepsilon \mu})$$
 with $\mu > 0$. This function preserves strong classes. Indeed:
$$\frac{1}{\varepsilon} d(\delta_{\varepsilon} \left( \delta_{\mu} \circ c_{1} (\frac{1}{\varepsilon}) \right) z, \delta_{\varepsilon} \left( \delta_{\mu} \circ c_{2} (\frac{1}{\varepsilon}) \right)z) \ =$$
$$ = \ \mu  \frac{1}{\varepsilon \mu} d(f_{1}(\varepsilon \mu) z , f_{2}(\varepsilon \mu) z)$$
which converges to $0$ as $\varepsilon \rightarrow 0$ because $c_{1} \strong
c_{2}$.

We prove now that classical dilatations transport strong classes into strong classes. For this we shall use the following characterisation of
strong equivalence: $c_{1} \strong c_{2}$ if and only if for any $z$ (in a neighbourhood of the neutral element) we have
$$d_{N}( \delta_{\varepsilon}^{-1} (f_{1}(\varepsilon) \opn z ),
\delta_{\varepsilon}^{-1} (f_{1}(\varepsilon) \opn z )) \rightarrow 0$$
as $\varepsilon \rightarrow 0$, uniformly with respect to $z$. Here
$d_{N}$ is the distance on the tangent cone at the neutral element
(remember that, as in Vodop'yanov \& Greshnov \cite{vodopis2}, we have identified a neighbourhood
of the neutral element of the tangent cone with a neighbourhood of the neutral element of the group $G$; on this neighbourhood we have also the distance $d_{N}$).

 With this characterisation $\delta_{\mu} * c_{f_{1}} \strong \delta_{\mu} * c_{f_{2}}$ if and only if
$$ \frac{\mu}{\varepsilon} d_{N}(f_{1}(\varepsilon) \opn \delta_{\mu}^{-1} z,
f_{1}(\varepsilon) \opn \delta_{\mu}^{-1} z) \rightarrow 0$$
as $\varepsilon \rightarrow 0$, uniformly with respect to $z$, which is true
by hypothesis.

It is easy to see that translation and dilatations (in both senses: $\delta_{\mu} \circ ... $ and $\delta_{\mu} * ...$) preserve the class of words and the class
of curved induced from functions in $\mathcal{C}^{m}([0,+\infty), G)$.

Coming back to transport, we notice that any  Lipschitz maps $f$ transports weak classes into weak classes. We cannot use this to define a tangent map  the tangent bundles in the sense introduced here, because such tangent bundles are constructed with strong classes.

The following lemma (with straightforward proof) shows that Pansu derivative gives only a partial description about how $f$ transports finite curves.   Here by Pansu derivative we mean a derivative in the sense of definition \ref{dpmmv}.

\begin{lema}
If $f: N \rightarrow L$ is Pansu derivable in $x$ then it transports any
two letter word of the form $\displaystyle \left( \begin{array}{cc}
x &  \\
 & y
\end{array} \right)$ to a curve weakly equivalent to the word $\displaystyle f(x) \left( \begin{array}{cc}
0 &  \\
 & Pf(x)y
\end{array} \right)$, where $Pf(x)$ denotes the Pansu derivative of $f$ at $x$.
\end{lema}

From the point of view developed in this section, the meaning of Rademacher theorem is that any Lipschitz function transports  two letter word weak classes into two letter word weak classes, almost everywhere with respect to the first letter.

The problem that we have in order to introduce a derivative notion is that we would like the derivative of $f$ at $x$ to act on strong classes, but  (the Pansu derivative) provide a transformation
of (a particular family of) weak classes.

Before introducing the derivative notion in this setting, we shall look closer
to Carnot groups.

\section{Calculus on Carnot groups}

Let $N$ be a Carnot group. We shall develop the notions introduced in the previous section. But let us first justify the approach by looking at the right translations in the Heisenberg group.

Consider therefore the right translation $f: H(n) \rightarrow H(n)$,
$$f(x,\bar{x}) = (x+y, \bar{x} + \bar{y} + \frac{1}{2}\omega(x,y))$$
We look at $f(\tilde{x}\delta_{\varepsilon}\tilde{z})$, with $\tilde{x}, \tilde{z}
\in H(n)$, $\tilde{x} = (x,\bar{x})$, $\tilde{z} = (z,\bar{z})$. Then
$$f(\tilde{x}\delta_{\varepsilon}\tilde{z}) \ = \ f(x + \varepsilon z,
\bar{x} + \varepsilon\frac{1}{2}\omega(x,z)  +\varepsilon^{2}\bar{z}) \ =$$
$$= \ (x+y + \varepsilon z, \bar{x} + \bar{y} + \frac{1}{2}\omega(x,y) + \varepsilon \left( \frac{1}{2}\omega(z,y) + \frac{1}{2}\omega(x,z)\right) +
\varepsilon^{2} \bar{z})$$
We compute now the expression
$$\left(f(\tilde{x})\right)^{-1}f(\tilde{x}\delta_{\varepsilon}\tilde{z}) \ = \
(\varepsilon z, \varepsilon \omega(z,y) + \varepsilon^{2} \bar{z})$$
Because of the term $\varepsilon \omega(z,y)$, the right translation is not derivable, unless $y = 0$. Nevertheless, what is to be remarked is that the  polynomial expression $(\varepsilon z, \varepsilon^{2}\bar{z})$ has been transported to $(\varepsilon z, \varepsilon \omega(z,y) + \varepsilon^{2} \bar{z})$, which is not so bad.

We go back now to general Carnot groups. We shall see that in this case we can
pass from weak classes to strong classes by  choosing representantives in a
canonical way.

\begin{defi}
Let $C(N)$ the group of polynomial curves
$$\varepsilon > 0 \ \mapsto \ \left( P_{i}(\varepsilon) \right)_{i = 1,...,n}
\ \ , \ \mbox{deg} P_{i}(\varepsilon) \leq i \ , \ \ P_{i}(\varepsilon) \in
V_{i}[\varepsilon]$$
This is a group with respect to the pointwise multiplication in $N$.

To any variational word
$$c \ = \ \left( \begin{array}{cccccc}
x_{1} & & x_{2} & ... & x_{p} & \\
 & y_{1} &  & y_{2} & ... & y_{p}
\end{array} \right) $$
we associate the element $\displaystyle P_{c}(\varepsilon) \ = \ \delta_{\varepsilon} c (\frac{1}{\varepsilon})$. The class of all such polynomials,
denoted by $WN$,  forms a subgroup of $C(N)$.

Conversely, to any $P \in C(N)$ we associate the finite curve
$$c_{P}(\frac{1}{\varepsilon}) \ = \ \delta_{\varepsilon}^{-1} P(\varepsilon)$$
\end{defi}

A curve $c$ is {\it good} if $c \weak c_{P}$ for some $P \in C(N)$. The curve
$c$ is {\it very good} if $c \weak c_{P}$ with $P \in WN$.

The tangent bundle of $N$ is by definition $WN$.  We collect in a
proposition some basic facts about weak and strong classes of good and very good curves.

\begin{prop}
\begin{enumerate}
\item[(a)] Let $P, Q \in C(N)$. $c_{P} \weak c_{Q}$ if and only if $P=Q$. Therefore on the class of good curves the weak classes are equal with strong classes. Otherwise said: the spectrum of a good curve is trivial (in the class of good curves).
\item[(b)] The composition operation between strong classes of good and very good curves descends on the group operation on $C(N)$, respectively $WN$.
\end{enumerate}
\label{pwts}
\end{prop}

To any left or right translation in $N$ there is an associated map on the tangent bundle. Indeed, let $x \in N$; we look at the transport given
by left translations of variational words:
$$L_{x}*\left( \begin{array}{cccccc}
x_{1} & & x_{2} & ... & x_{p} & \\
 & y_{1} &  & y_{2} & ... & y_{p}
\end{array} \right) \ = \ \left( \begin{array}{cccccc}
x x_{1} & & x_{2} & ... & x_{p} & \\
 & y_{1} &  & y_{2} & ... & y_{p}
\end{array} \right)$$
In the same way
$$R_{x} * \left( \begin{array}{cccccc}
x_{1} & & x_{2} & ... & x_{p} & \\
 & y_{1} &  & y_{2} & ... & y_{p}
\end{array} \right) \ = \ \left( \begin{array}{cccccccc}
x_{1} & & x_{2} & ... & x_{p} & & x & \\
 & y_{1} &  & y_{2} & ... & y_{p} & & 0
\end{array} \right)$$
For that dilation $\delta_{\mu}$ we have:
$$\delta_{\mu} * \left( \begin{array}{cccccc}
x_{1} & & x_{2} & ... & x_{p} & \\
 & y_{1} &  & y_{2} & ... & y_{p}
\end{array} \right) \ = \ \left( \begin{array}{cccccc}
\delta_{\mu} x_{1} & & \delta_{\mu}x_{2} & ... & \delta_{\mu}x_{p} & \\
 & \delta_{\mu}y_{1} &  & \delta_{\mu}y_{2} & ... & \delta_{\mu}y_{p}
\end{array} \right)$$
We also have the formula
$$\delta_{\mu} \circ \left( \begin{array}{cccccc}
x_{1} & & x_{2} & ... & x_{p} & \\
 & y_{1} &  & y_{2} & ... & y_{p}
\end{array} \right) \ = \ \left( \begin{array}{cccccc}
 x_{1} & & x_{2} & ... & x_{p} & \\
 & \delta_{\mu}y_{1} &  & \delta_{\mu}y_{2} & ... & \delta_{\mu}y_{p}
\end{array} \right)$$

Remark that the tangent bundle is a Carnot group. We consider the generating
distribution
$$D^{2} \ = \ \left\{ \left( \begin{array}{cc}
x &  \\
 & y
\end{array} \right) \mbox{ : } x,y \in D \right\}$$
and the Euclidean norm
$$\mid \left( \begin{array}{cc}
x &  \\
 & y
\end{array} \right) \mid \ = \ \left( \mid x \mid^{2} + \mid y \mid^{2} \right)^{\frac{1}{2}}$$
(although any $L_{p}$ norm is good as well).  This gives us a distance on $WN$, denoted by $d_{W}$.

\begin{rk}
 The fact that the tangent bundle  is a Carnot group is a key property, because we stay
in the same category after derivation. Hence one can do an indefinite number of
derivatives.
\end{rk}

The tangent space at the neutral element is
$$T_{0}N \ = \ \left\{ P \in C(N) \mbox{ : } P(0) \ = \ 0 \right\}$$
and it is a Carnot group too. The dimension of $T_{0}N$ equals the homogeneous dimension of $N$. The tangent bundle $WN$ is a semidirect
product between $N$ and $W_{0}N$.

The dilatations in $T_{0}N$ are defined by:
$$\delta_{\mu}(P)(\varepsilon) \ = \ P(\mu \varepsilon)$$
for any $\mu >0$.

\begin{defi}
The function $f: N \rightarrow L $ is weakly derivable at $x \in N$ if there is
a function $Tf(x): T_{0}N \rightarrow WL$ with the following property:
\begin{equation}
\frac{1}{\varepsilon} d \left( f(x P(\varepsilon)), Tf(x)(P)(\varepsilon)\right) \rightarrow 0
\label{dwd}
\end{equation}
as $\varepsilon \rightarrow 0$, uniformly with respect to $P \in B(0,1) \subset
T_{0}N$.
\end{defi}

It is straightforward from this definition that there exist $Df(x): T_{0}N \rightarrow T_{0}L$ such that
$$Tf(x)(P) \ = \ f(x) Df(x)(P)$$
and $Df(x)$ commutes with dilatations.

The following is a simple, but important observation.
\begin{prop}
Let $f: N \rightarrow N$ be a function, $f = (f_{1}, .... , f_{m})$, where $m$
is the step of $N$. Suppose that $f_{i}$ admits a development of order $i$ around $x$.  Then $f$ is weakly derivable.
\label{pepe}
\end{prop}

The equivalent of proposition \ref{ppansuprepgen} is not true for weak derivatives, because the quantities from the proof, denoted by
$(1)_{t}$, $(2)_{t}$ might converge to infinity. That is why we state two more
definitions.

\begin{defi}
The function $f: N \rightarrow L $ is  derivable at $x \in N$ if there exists
{\bf a group morphism} $Df(x): T_{0}N \rightarrow T_{0}L$ such that
$$\frac{1}{\varepsilon} d \left( f(x P(\varepsilon)), f(x) Df(x)(P)(\varepsilon)\right) \rightarrow 0 $$
as $\varepsilon \rightarrow 0$, uniformly with respect to $P \in \bar{B(0,1)} \subset T_{0}N$.
\label{mildder}
\end{defi}

\begin{defi}
The function $f: N \rightarrow L $ is  strongly derivable at $x \in N$ if there exists
$Df(x): T_{0}N \rightarrow T_{0}L$ such that
$$\frac{1}{\varepsilon} d \left( f(x P(\varepsilon)) y, f(x) Df(x)(P)(\varepsilon)y \right) \rightarrow 0 $$
as $\varepsilon \rightarrow 0$, for any $y \in N$,   uniformly with respect to $P \in \bar{B(0,1)} \subset T_{0}N$.
\end{defi}

Now comes the equivalent of proposition \ref{ppansuprepgen}.

\begin{prop}
Let $f: N \rightarrow L$ be a function.
\begin{enumerate}
\item[a)] If $f$ is strongly derivable at $x$ then it is derivable.
\item[b)] If $f$ is weakly derivable and
\begin{enumerate}
\item[-] the convergence in the definition of the derivative is uniform with respect to $x$,
\item[-] $f$  satisfies the following condition (A) (from "additive"): for any $P,Q \in B(0,1) \subset
T_{0}N$ and any $x \in N$
\begin{equation}
\mbox{{\bf (A)}} \ \  \ \ \ \frac{1}{\varepsilon} d\left(f(x P(\varepsilon)Q(\varepsilon)), f(xP(\varepsilon))
f(x)^{-1}f(xQ(\varepsilon))\right) \rightarrow 0
\label{A}
\end{equation}
goes to $0$ as $\varepsilon \rightarrow 0$, uniformly with respect to $P,Q$,
\end{enumerate}
then $f$ is derivable.
\end{enumerate}
\label{pupu}
\end{prop}

One can imagine many sufficient conditions for (A). One of them is the following:
there exists $C >0$ such that for any $u,x,z \in N$ we have
\begin{equation}
d(f(u)^{-1}f(uz) , f(x)^{-1} f(xz) ) \ \leq \ C d(u,x) d(0,z)^{m}
\label{Asuf}
\end{equation}

Let us see what \eqref{Asuf} means in the commutative case $N = R^{n}$. Take
$f: R^{n} \rightarrow R^{k}$ and $x, y \in R^{n}$. Construct the function
$\displaystyle g_{xy}(z) = f(x+z) - f(y+z)$. If $f$ has Lipschitz derivatives then
it satisfies condition \eqref{Asuf}, which in this case reads:
$$\| g_{xy}(z) - g_{xy}(0) \| \leq C \|x-y\| \|z\|$$
The constant $C$ equals $Lip(\nabla f)$.

\paragraph{Proof.}(of proposition \ref{pupu})
We shall prove only the point (b), because (a) has an easy proof. We have to prove the morphism property of $Df(x)$.

Let $P,Q \in T_{0}N$. Because $Df(y)$ commutes with dilatations for any $y \in U$, we can suppose that $P,Q, \in \bar{B(0,1)} \subset T_{0}N$.

We start from the inequality:
$$\frac{1}{\varepsilon} d \left( f(x(P)(\varepsilon)Q(\varepsilon)), f(x)Df(x)(P)(\varepsilon) Df(x)(Q)(\varepsilon)\right) \ \leq A(\varepsilon) + B(\varepsilon)$$
with $A(\varepsilon) \rightarrow 0$ because of condition (A). We used the notations:
$$A(\varepsilon) \ = \ \frac{1}{\varepsilon} d\left(f(x P(\varepsilon)Q(\varepsilon)), f(xP(\varepsilon))
f(x)^{-1}f(xQ(\varepsilon))\right)$$
$$B(\varepsilon) \ = \ \frac{1}{\varepsilon} d\left(f(xP(\varepsilon))
f(x)^{-1}f(xQ(\varepsilon)), f(x)Df(x)(P)(\varepsilon) Df(x)(Q)(\varepsilon)\right)$$
We continue:
$$B(\varepsilon) \ \leq \ C(\varepsilon) + D(\varepsilon)$$
with the notations
$$C(\varepsilon) \ = \ \frac{1}{\varepsilon} d\left(f(xP(\varepsilon))
f(x)^{-1}f(xQ(\varepsilon)),f(x P(\varepsilon))Df(x)(Q)(\varepsilon)\right)$$
\begin{equation}
D(\varepsilon) \ = \ \frac{1}{\varepsilon} d\left(f(x P(\varepsilon))Df(x)(Q)(\varepsilon), f(x)Df(x)(P)(\varepsilon) Df(x)(Q)(\varepsilon)\right)
\label{clueref}
\end{equation}
\begin{rk}
The quantity $D(\varepsilon)$ will be important in the section concerning general
groups, in relation to the definition \ref{dmild} of mild equivalence.
\end{rk}
The  quantity $C(\varepsilon)$  converges to $0$ as $\varepsilon \rightarrow 0$ because $f$ is weakly derivable. The second quantity $D(\varepsilon)$ behaves  in the same way again because $f$ is weakly derivable, because  of the uniformity assumption in part (b) of the hypothesis and because of proposition \ref{pwts} (a). We have proved that
$$\frac{1}{\varepsilon} d \left( f(x(P)(\varepsilon)Q(\varepsilon)), f(x)Df(x)(P)(\varepsilon) Df(x)(Q)(\varepsilon)\right)$$
converges to $0$ as $\varepsilon \rightarrow 0$, uniformly with respect to
$P,Q, \in \bar{B(0,1)} \subset T_{0}N$.

The one thing left to do is to use once again the fact that $f$ is weakly derivable, this time along $P(\varepsilon)Q(\varepsilon)$. This fact, together
with the convergence previously established, end the proof of the morphism property.
\quad $\blacksquare$

It seems hard to believe that these three definitions
of derivability are really different. This is a purely non-Euclidean property. We shall see in the section dedicated to the Heisenberg group that the differences are really big and
that the good notion of derivability seems to be \ref{mildder}.

The Rademacher theorem for our notion of derivative is:

\begin{thm}
Let $f: N \rightarrow L$ be a Lipschitz function which satisfies condition (A).
Then $f$ is  derivable a.e. and the derivative of $f$ has the property:
$$Df(x)\left( \begin{array}{cc}
x &  \\
 & y
\end{array} \right) \ = \ \left( \begin{array}{cc}
0 &  \\
 & Pf(x)y
\end{array} \right)$$
\label{tradnew}
\end{thm}

The proof is a direct consequence of  Pansu-Rademacher theorem,  propositions
\ref{pepe},  \ref{pupu} and the observation that the Margulis-Mostow tangent bundle (bundle of two-letter words) generates the word tangent bundle.

We end this section with the remark that
the derivative satisfies the chain rule by the very definition.

\section{Case of the Heisenberg group}

In this case $CH(n) = WH(n)$, which comes from the fact that two letter words generate $CH(n)$. $T_{0}H(n)$ is the group of all polynomials with the form:
$$P(\varepsilon) \ = \ (u,u',u")(\varepsilon) \ = \ (\varepsilon u , \varepsilon u' + \varepsilon^{2} u")$$
with $u \in R^{2n}, u', u" \in R$. This group is isomorphic with
$H(n)\times R$ and the isomorphism is given by:
$$(u,u',u") \in T_{0}H(n) \mapsto ((u,u"), u') \in H(n) \times R$$
The morphisms from $R$ to $T_{0}H(n)$  which commute with dilatations forms a
set called $HL(R,T_{0}H(n))$ ("HL" means "horizontal linear"). Any element $F$ of this set can be seen as an element $(a,a') \in R^{2n} \times R$, meaning that
$F \in HL(R, T_{0}H(n))$ has the form:
$$F(u) \ = \ (u a , u'a', 0)$$

\paragraph{Derivable curves}
Let $\tilde{f}: R \times H(n)$, $f(t) = (f(t), \bar{f}(t))$ derivable. This means that there is $(a(t), a'(t))$ such
that for any $u \in R, \mid u \mid < 1$ we have:
$$\delta_{\varepsilon}^{-1} \left((-f(t+\varepsilon u), -\bar{f}(t+\varepsilon u))
(f(t), \bar{f}(t))(\varepsilon u a(t), \varepsilon u a'{t}) \right) \rightarrow 0$$
as $\varepsilon \rightarrow 0$, uniformly with respect to $u$. We obtain the following characterisation of $\tilde{f}$:
\begin{enumerate}
\item[a)] $f$ is derivable with respect to $t$ and $a(t) = \dot{f}(t)$,
\item[b)] the map $\displaystyle s \mapsto \bar{f}(t+s) - \frac{1}{2} \omega(f(t),
f(t+s))$ has a second order development at $s=0$ of the form:
$$\bar{f}(t+s) - \frac{1}{2} \omega(f(t),
f(t+s)) = \bar{f}(t) + s a'(t) + o(s^{2})$$
\end{enumerate}
Putting a) and b) together  we get: $\tilde{f}$ is classically derivable and
$$\left\{ \begin{array}{l}
a(t) = \dot{f}(t) \\
a'(t) = \dot{\bar{f}} - \frac{1}{2} \omega(f(t),\dot{f}(t)) \\
\frac{d}{dt}a'(t) = 0
\end{array} \right.$$
So the curve $\tilde{f}$ is not horizontal generally, but the difference from horizontality, measured by $a'$,  is constant along the curve. A straightforward
computation says that when $f$ is parametrised by arclength, the $2$ Hausdorff measure of $\tilde{f}$ is $a'L(f)$, where $L(f)$ is the length of $f$.

An analoguous computation shows that $\tilde{f}$ is strongly derivable if and only if it is linear.

\paragraph{Derivative of the right translation} The right translations are derivable. We continue with the motivating example from the beginning of the previous section.

We saw there that the right translation $f: H(n) \rightarrow H(n)$,
$f(\tilde{x}) \ = \ \tilde{x} \tilde{y}$ transports the polynomial
$(z,0,\bar{z})$ to  $(z,\omega(z,y),\bar{z})$. The same computation
shows that
$$Df(\tilde{x}) (u,u',u") \ = \ (u, u'+\omega(u,y), u")$$
We can write this derivative as a matrix:
$$Df(\tilde{x}) \ = \ \left( \begin{array}{ccc}
1 & 0 & 0 \\
-Jy & 1 & 0 \\
0 & 0 & 1
\end{array} \right)$$
and check that $Df(\tilde{x})$ is indeed a morphism of $T_{0}H(n)$ which commute with dilatations. That means $Df(\tilde{x}) \in HL(T_{0}H(n))$, where 
$$HL(T_{0}H(n)) \ = \ \left\{ F \in GL(T_{0}H(n)) \mbox{ : } \forall x,y \in T_{0}H(n) \ [Fx,Fy] = F[x,y] \ , \ \forall \varepsilon > 0 \ \delta_{\varepsilon} F = F \delta_\varepsilon \right\}$$

\begin{rk}
It has been shown in Ambrosio, Rigot \cite{ambrig} that the problem of optimal transportation in the Heisenberg group has solutions. In classical situations the optimal transport map $\psi$ satisfies the Monge-Amp\`ere equation, which contains
$\det \nabla \psi$. In the case of the Heisenberg group one can prove that right translations are optimal transport maps (for a particular choice of the measures
to be transported). But right translations are not Pansu differentiable, hence the problem of understanding the Monge-Amp\`ere equation in this setting. We don't deal in this paper with the calculus of variations associated to the notions
of derivative just introduced. This will be done in a future paper. I think that
it is encouraging though the fact that right translations are derivable, which might lead to an interpretation of the Monge-Amp\`ere equation in the Heisenberg group.
\end{rk}

\paragraph{Derivable homeomorphisms} We look first for the class of isomorphisms
of $T_{0} H(n)$ which commute with dilatations. This class is denoted by $HL(T_{0}H(n))$. We get $\mathcal{A} \in HL(T_{0}H(n))$ if and only if $\mathcal{A}$ is represented by a matrix:
$$\left( \begin{array}{ccc}
A & 0 & 0 \\
c & d & 0 \\
0 & 0 & h
\end{array} \right)$$
such that $A$ is conformally symplectic and $\displaystyle h = (\det A)^{\frac{1}{n}}$ .

We proceed like in the case of curves: let $f: H(n) \rightarrow H(n)$,
$f(x,\bar{x}) \ = \ (\phi(x, \bar{x}), \bar{\phi}(x,\bar{x}))$. If $f$ is derivable at $(x,\bar{x})$ then the first conclusion is that $\phi \ = \ \phi(x)$.
An interesting consequence is obtained if $f$ is also volume preserving. In this case $\phi$ is symplectic and moreover
$$\bar{\phi}(x, \bar{x}) \ = \ \bar{x} + c \cdot x + F(x)$$
where $F$ is the generating function of $\phi$ and $c$ does not depend on $(x, \bar{x})$.

We paste here the section 3.2, 3.3  of \cite{buliga1}, for completeness of the
picture of derivability in the Heisenberg group. Before this we state  Theorem \ref{corect}, which  is
correct version of theorem 3.8. from \cite{buliga1}.

\begin{thm}
Take any   $\tilde{\phi}$ locally bi-Lipschitz volume preserving homeomorphism
 of $H(n)$, {\bf which is derivable}. Then $\tilde{\phi}$ has the form:
$$\tilde{\phi}(x,\bar{x}) \ = \ (\phi(x), \bar{x} \ + \ F(x))$$
Moreover $$\phi \in Sympl(R^{2n}, Lip)$$
  $F: R^{2n} \rightarrow R$ is
Lipschitz and for almost any point  $(x,\bar{x}) \in H(n)$ we have:
$$DF(x) y \ = \ \frac{1}{2} \omega(\phi(x), D\phi(x)y) \ - \ \frac{1}{2}\omega(x,y)$$
\label{corect}
\end{thm}

\paragraph{Proof.}
Let $(u,u',u") \in T_{0}H(n)$ and $(x,\bar{x}) \in H(n)$.
The condition on $f$ to be derivable has the form:
$$E(\varepsilon) = \delta_{\varepsilon}^{-1} \left\{ (-\phi((x,\bar{x})(u,u',u")(\varepsilon)),
-\bar{\phi}((x,\bar{x})(u,u',u")(\varepsilon)))(\phi(x,\bar{x}), \bar{\phi}(x,\bar{x})) \right. $$
$$ \left. (\varepsilon Au, \varepsilon cu + \varepsilon d u' + \varepsilon^{2} h u") \right\} \rightarrow 0$$
as $\varepsilon \rightarrow 0$, uniformly with respect to $(u,u',u")$ such that
$$\mid u\mid + \mid u' \mid + \mid u" \mid^{\frac{1}{2}} \leq 1$$
Here $A$ is conformally symplectic and $\displaystyle h = (\det A)^{\frac{1}{n}}$.
The $R^{2n}$ component of the $E(\varepsilon)$ is
$$\frac{1}{\varepsilon} \left( \phi(x,\bar{x}) - \phi(x+\varepsilon u, \bar{x} +
\varepsilon u' + \varepsilon \frac{1}{2}\omega(x,u) + \varepsilon^{2} u") +
\varepsilon Au \right)$$
and it has to converge to $0$ as $\varepsilon \rightarrow 0$. Take now $u=0$,
$u'=0$ and get that $\phi$ is derivable with respect to $\bar{x}$ and the
derivative equals $0$. Therefore $\phi = \phi(x)$. From here we get that
$\phi$ is classically derivable with respect to $x$ and
$$A \ = \ \frac{\partial \phi}{\partial x}$$
The second component of $E(\varepsilon)$ admits a similar treatment. We leave the details to the reader. The conclusions are:
\begin{enumerate}
\item[1)] $\bar{\phi}$ is derivable with respect to $x$ and $\bar{x}$,
\item[2)] we have $\displaystyle \frac{\partial \bar{\phi}}{\partial \bar{x}} =
d = h$,
\item[3)] we have also the equality
$$c u \ = \ \frac{\partial \bar{\phi}}{\partial x} u - \frac{1}{2} \omega(
\phi, \frac{\partial \phi}{\partial x} u ) + \frac{\partial \bar{\phi}}{\partial \bar{x}} \frac{1}{2} \omega(x,u)$$
for any $u \in R^{2n}$.
\item[4)] the coefficient $c$ has to be constant with respect to $(x,\bar{x})$.
\end{enumerate}
The function $f$ is Lipschitz therefore Pansu derivable almost everywhere. Because
 f is also volume preserving, the area formula (change of variables formula, we
have not touched this subject in this paper, see for example Vodop'yanov, Ukhlov 
\cite{vodopis} \cite{voduk} or  Pauls \cite{pauls1}
and Magnani \cite{magnani}) implies that a.e.  $A$ have to be symplectic, hence $h = d = 1$.

Use 2), 3), 4) to get
$$\frac{\partial \bar{\phi}}{\partial x} \ = \ \frac{1}{2}\omega(\phi, \frac{\partial \phi}{\partial x} u) - \frac{1}{2}\omega(x,u) + c$$
But the generating function $F$ of $\phi$ is defined up to an additive constant by
$$DF(x) u \ = \ \frac{1}{2}\omega(\phi, \frac{\partial \phi}{\partial x} u) - \frac{1}{2}\omega(x,u)$$
(This function exists, as a consequence of the fact that $\phi$ is symplectic a.e.) The proof is finished. \quad $\blacksquare$

\subsection{Hamiltonian dynamics in a nutshell}

The purpose of this section is to use the sub-Riemannian calculus on the Heisenberg group to recover the basic notions in Hamiltonian dynamics. We
show that Hamiltonian dynamics is just the study of flows of volume preserving
derivable homeomorphisms of the Heisenberg group. Notions of generating function,
hamiltonian, Hamilton equation, Hofer distance, appear without premeditation.
One might wonder what is the result of the study of volume preserving derivable
homeomorphisms in a general Carnot group. The case of H-type Carnot groups
may as well be linked to Clifford analysis (because there is a natural correspondence between Clifford algebras and H-type groups,  see Reimann \cite{reim} or Barbano \cite{barbano}).

\paragraph{Classically smooth maps and the derivative: lifts of symplectic diffeomorphisms}

In this section we are interested in the
group of volume preserving classically diffeomorphisms of $H(n)$, which are also derivable in our new sense. Recall a previous proposition which roughly says that a sufficiently regular function is weakly derivable. In order that the function to be derivable some supplementary conditions upon higher order derivatives have to be satisfied. What is interesting
in this subsection is the proof that  classically smooth enough {\bf volume preserving}
diffeomorphisms are also derivable, not only weakly derivable.

\paragraph{Volume preserving diffeomorphisms}

Let $A \subset R^{2n}$, open,  such that $0 \not \in A$. This condition does not restrict the generality; it leads only to a simplified notation.

\begin{defi}
$Diff^{2}(A,vol)$ is  the group of volume preserving
diffeomorphisms $\tilde{\phi}$ of $H(n)$ such that:
\begin{enumerate}
\item[a)] $\tilde{\phi}$ has horizontal compact support in $A \times R \subset H(n)$, in the sense that it differs from a map $(x, \bar{x}) \mapsto (x, \bar{x} + c)$ ($c \in R$, arbitrar), only in a compact set included in $A \times R$,
\item[b)]  $\tilde{\phi}$ and it's inverse have (classical) regularity $C^{2}$.
\end{enumerate}

In the same way we define $Sympl^{2}(A)$ to be the group of
$C^{2}$ symplectomorphisms of $R^{2n}$ with compact support in $A$.
\label{dd}
\end{defi}

For the regularity chosen in the definition \eqref{dd}, the result equivalent with theorem 
\eqref{corect} is given further. We denote further by $\lambda$ a one-form on $R^{2n}$ such that 
$2\ d \lambda \ = \ \omega$. 

\begin{thm}
We have the isomorphism of groups
$$Diff^{2}(A,vol) \ \approx \ Sympl^{2}(A) \times R$$
given by the mapping
$$\tilde{f} = (f,\bar{f}) \  \in \ Diff^{2}(A,vol) \ \mapsto \ \left(
f \in Sympl^{2}(A) , \bar{f}(0,0) \right)$$
The inverse of this isomorphism has the expression
$$\left( f \in Sympl^{2}(A) , a \in R \right) \ \mapsto  \  \tilde{f} = (f,\bar{f}) \  \in \ Diff^{2}(A,vol)$$
$$\tilde{f}(x,\bar{x}) \ = \ (f(x), \ \bar{x} + F(x))$$
where $F(0)= a$ and $dF \ = \ f^{*} \lambda \ - \ \lambda$.
\label{t1}
\end{thm}

\begin{rk}
We need the assumption $0 \not \in A$ only for the morphism property.
\end{rk}

\paragraph{Proof.}
Let $\tilde{f} = (f,\bar{f}) : H(n) \rightarrow H(n)$ be an element of
the group \\
$Diff^{2}(A, vol)$.
We shall compute:
$$D \tilde{f} ((x,\bar{x})) (y,\bar{y}) \ = \
\lim_{\varepsilon \rightarrow 0} \delta_{\varepsilon^{-1}}  \left(
\left(\tilde{f}(x,\bar{x})\right)^{-1}  \tilde{f}\left((x,\bar{x})  \delta_{\varepsilon} (y,\bar{y})\right)\right) $$
We know that $D \tilde{f}(x,\bar{x})$ has to be a linear mapping.

After a short computation we see that we have to pass to the limit
$\varepsilon \rightarrow 0$ in the following expressions (representing the two
components of $D \tilde{f} ((x,\bar{x})) (y,\bar{y})$):
\begin{equation}
\frac{1}{\varepsilon} \left( f\left(x+ \varepsilon y, \bar{x} + \varepsilon^{2}
\bar{y} + \frac{\varepsilon}{2} \omega(x,y)\right) - f(x,\bar{x}) \right)
\label{exp1}
\end{equation}
\begin{equation}
\frac{1}{\varepsilon^{2}} \left( \bar{f}\left(x+ \varepsilon y, \bar{x} + \varepsilon^{2}
\bar{y} + \frac{\varepsilon}{2} \omega(x,y)\right) - \bar{f}(x,\bar{x}) -
\right.
\label{exp2}
\end{equation}
$$\left. -
\frac{1}{2}\omega\left(f(x,\bar{x}), f\left(x+ \varepsilon y, \bar{x} + \varepsilon^{2}
\bar{y} + \frac{\varepsilon}{2} \omega(x,y)\right)\right)\right)
$$

The first component \eqref{exp1} tends to
$$\frac{\partial f}{\partial x} (x,\bar{x}) y + \frac{1}{2} \frac{\partial f}{\partial \bar{x}} (x,\bar{x}) \omega(x,y)$$
The terms of order $\varepsilon$ must cancel in the second component \eqref{exp2}. We obtain the cancellation
condition (we shall omit from now on the argument $(x,\bar{x})$ from all functions):
\begin{equation}
\frac{1}{2} \omega(x,y) \frac{\partial \bar{f}}{\partial \bar{x}} -
\frac{1}{2} \omega(f, \frac{\partial f}{\partial x} y) -
\frac{1}{4} \omega(x,y) \omega(f, \frac{\partial f}{\partial \bar{x}}) +
\frac{\partial \bar{f}}{\partial x} \cdot y \ = \ 0
\label{cancel}
\end{equation}
The second component tends to
$$\frac{\partial \bar{f}}{\partial \bar{x}} \bar{y} - \frac{1}{2} \omega(f,
\frac{\partial f}{\partial \bar{x}}) \bar{y}$$
The group morphism $D \tilde{f}(x,\bar{x})$
is represented by the matrix:
\begin{equation}
d \tilde{f}(x,\bar{x}) \ = \ \left( \begin{array}{cc}
\frac{\partial f}{\partial x} + \frac{1}{2} \frac{\partial f}{\partial \bar{x}}
\otimes Jx & 0 \\
0 & \frac{\partial \bar{f}}{\partial \bar{x}} - \frac{1}{2} \omega(f,
\frac{\partial f}{\partial \bar{x}})
       \end{array} \right)
\label{tang}
\end{equation}
We shall remember now that $\tilde{f}$ is volume preserving. This implies:
\begin{equation}
\frac{\partial f}{\partial x} + \frac{1}{2} \frac{\partial f}{\partial \bar{x}}
\otimes Jx \ \in Sp(n)
\label{c1}
\end{equation}
\begin{equation}
\frac{\partial \bar{f}}{\partial \bar{x}} - \frac{1}{2} \omega(f,
\frac{\partial f}{\partial \bar{x}}) = 1
\label{c2}
\end{equation}
The cancellation condition \eqref{cancel} and relation \eqref{c2} give
\begin{equation}
\frac{\partial \bar{f}}{\partial x} y \ = \  \frac{1}{2} \omega(f,\frac{\partial
f}{\partial x} y ) \ - \ \frac{1}{2} \omega(x,y)
\label{c3}
\end{equation}

These conditions describe completely the class of volume preserving diffeomorphisms
of $H(n)$. Conditions \eqref{c2} and \eqref{c3} are in fact differential equations
for the function $\bar{f}$ when $f$ is given. However, there is a compatibility
condition in terms of $f$ which has to be fulfilled for  \eqref{c3} to have
a solution $\bar{f}$. Let us look closer to \eqref{c3}. We can see the symplectic
form $\omega$ as a closed 2-form. Let $\lambda$ be a 1-form such that
$2 d \lambda = \omega$. If we take the (regular) differential with respect
to $x$ in \eqref{c3} we quickly obtain the compatibility condition
\begin{equation}
\frac{\partial f}{\partial x} \ \in \ Sp(n)
\label{c4}
\end{equation}
and \eqref{c3} takes the form:
\begin{equation}
 \ d \bar{f} \ = \ f^{*} \lambda \ - \ \lambda
\label{c5}
\end{equation}
(all functions seen as functions of $x$ only).

Conditions \eqref{c4} and \eqref{c1} imply: there is a scalar function
$\mu = \mu(x,\bar{x})$ such that
$$\frac{\partial f}{\partial \bar{x}} \ = \ \mu \ Jx $$
Let us see what we have until now:
\begin{equation}
\frac{\partial f}{\partial x} \ \in \ Sp(n)
\label{cc1}
\end{equation}
\begin{equation}
\frac{\partial \bar{f}}{\partial x} \ = \ \frac{1}{2} \left[ \left(
\frac{\partial f}{\partial x}\right)^{T} J f \ - \ J x \right]
\label{cc2}
\end{equation}
\begin{equation}
\frac{\partial \bar{f}}{\partial \bar{x}} \ = \ 1 + \frac{1}{2} \omega(f,
\frac{\partial f}{\partial \bar{x}})
\label{cc3}
\end{equation}
\begin{equation}
\frac{\partial f}{\partial \bar{x}} \ = \ \mu \ Jx
\label{cc4}
\end{equation}
Now, differentiate \eqref{cc2} with respect to $\bar{x}$ and use \eqref{cc4}. In the same time
differentiate \eqref{cc3} with respect to $x$. From the equality
$$\frac{\partial^{2} \bar{f}}{\partial x \partial \bar{x}} \ = \
\frac{\partial^{2} \bar{f}}{\partial \bar{x} \partial x}$$
we shall obtain by straightforward computation $\mu = 0$.

The morphism property is easy to  check.
\quad $\blacksquare$

\paragraph{Hamilton's equations}

For any element $\phi \in Sympl^{2}(A)$ we define the lift $\tilde{\phi}$ to be  the image
of $(\phi, 0)$ by the isomorphism described in the theorem \ref{t1}.

\begin{defi}
For any flow $t \mapsto \phi_{t} \in Sympl^{2}(A)$
denote by $\phi^{h}(\cdot, x))$ the  horizontal flow in
$H(n)$ obtained by the lift of all curves $t \mapsto \phi(t,x)$ and
by $\tilde{\phi}(\cdot, t)$ the flow obtained by the lift of all
$\phi_{t}$. The vertical flow is defined by the expression
\begin{equation}
\phi^{v} \ = \ \tilde{\phi}^{-1} \circ \phi^{h}
\label{hameq}
\end{equation}
\label{dhameq}
\end{defi}

Relation \eqref{hameq} can be seen as Hamilton equation, with known quantity $\phi^{v}$ and 
unknown quantities $\phi^{h}, \tilde{\phi}$. This is explained further.

\begin{prop}
Let $t \in [0,1] \mapsto \phi^{v}_{t}$ be a curve of diffeomorphisms
of $H(n)$ satisfying the equation:
\begin{equation}
\frac{d}{dt} \phi^{v}_{t}(x,\bar{x}) \ = \ (0, H(t,x)) \ \ , \ \
\phi^{v}_{0} \ = \ id_{H(n)}
\label{ham}
\end{equation}
Consider the Hamiltonian flow $t \mapsto \phi_{t}$ generated by $H$ and lift it using theorem \ref{t1} to the flow $t \mapsto \tilde{\phi}_{t}$; lift it also horizontally, as explained in definition  
\ref{dhameq}, to the flow $t \mapsto \phi_{t}^{h}$. Then equation   \eqref{hameq} is satisfied.
 
Conversely, for any Hamiltonian flow $t \mapsto \phi_{t}$, generated
by $H$, the vertical flow $ t \mapsto \phi^{v}_{t}$ (defined this time by the relation \eqref{hameq})  satisfies the equation \eqref{ham}.
\label{pham}
\end{prop}

\paragraph{Proof.}
Write the lifts $\tilde{\phi}_{t}$ and $\phi^{h}_{t}$, compute then
the differential of the quantity
$\dot{\tilde{\phi}}_{t} - \dot{\phi}^{h}_{t}$ and show that it
equals the differential of $H$. Use finally the fact that all quantities in the proof 
have compact support. 
\quad $\blacksquare$

\paragraph{Flows of volume preserving diffeomorphisms}

We want to know if there is any nontrivial smooth (according to Pansu differentiability)
flow of volume preserving diffeomorphisms.

\begin{prop}
Suppose that $t \mapsto \tilde{\phi}_{t} \in Diff^{2}(A,vol)$ is a flow such that
\begin{enumerate}
\item[-] is $C^{2}$ in the classical sense with respect to $(x,t)$,
\item[-] is horizontal, that is $t \mapsto \tilde{\phi}_{t}(x)$ is a horizontal curve
for any $x$.
\end{enumerate}
Then the flow is constant.
\label{pho}
\end{prop}

\paragraph{Proof.}
By direct computation, involving second order derivatives. Indeed, let
$\tilde{\phi}_{t}(x,\bar{x}) \ = \ (\phi_{t}(x), \bar{x} + F_{t}(x))$.
From the condition $\tilde{\phi}_{t} \in Diff^{2}(A,vol)$ we obtain
$$\frac{\partial F_{t}}{\partial x} y \ = \  \frac{1}{2} \omega(\phi_{t}(x),\frac{\partial
\phi_{t}}{\partial x}(x) y ) \ - \ \frac{1}{2} \omega(x,y)$$
and from the hypothesis that $t \mapsto \tilde{\phi}_{t}(x)$ is a horizontal curve
for any $x$ we get
$$\frac{d F_{t}}{dt}(x) \ = \ \frac{1}{2} \omega(\phi_{t}(x), \dot{\phi}_{t}(x))$$
Equal now the derivative of the RHS of the first relation with respect to $t$ with the derivative of the RHS of the second relation with respect to $x$. We get the equality, for any $y \in R^{2n}$:
$$ 0  \ = \ \frac{1}{2} \omega(\frac{ \partial \phi_{t}}{\partial x}(x) y, \dot{\phi}_{t}(x))$$
therefore $\dot{\phi}_{t}(x) \ = \ 0$.
\quad $\blacksquare$

One should expect such a result to be true, based on two remarks. The first, general remark: take a
flow of left translations in a Carnot group, that is a flow $t \mapsto \phi_{t}(x) \ = \ x_{t} x$.
We can see directly that each $\phi_{t}$ is smooth, because the distribution is
left invariant. But the flow is not horizontal, because the distribution is not
right invariant.  The second, particular remark: any  flow which satisfies the hypothesis
of proposition \ref{pho} corresponds to a Hamiltonian flow with null Hamiltonian function, hence the flow is constant.

\subsection{Symplectomorphisms, capacities and Hofer distance}

Symplectic capacities are invariants under the action of the symplectomorphisms group.
Hofer geometry is the geometry of the group of Hamiltonian diffeomorphisms, with respect
to the Hofer distance. For an introduction into the subject see Hofer, Zehnder
\cite{hozen} chapters 2,3 and 5, and Polterovich \cite{polte}, chapters 1,2.

 A symplectic capacity is a map which associates to any symplectic
manifold $(M,\omega)$ a number $c(M,\omega) \in [0,+\infty]$. Symplectic capacities
are special cases of conformal symplectic invariants, described by:
\begin{enumerate}
\item[A1.] Monotonicity: $c(M,\omega) \leq c(N,\tau)$ if there is a symplectic
embedding  from $M$ to $N$,
\item[A2.] Conformality: $c(M,\varepsilon \omega) = \mid \varepsilon \mid c(M,\omega)$
for any $\alpha \in R$, $\alpha \not = 0$.
\end{enumerate}

We can see a conformal symplectic invariant from another point of view. Take
a symplectic manifold $(M,\omega)$ and consider the invariant defined over the class
of Borel sets $B(M)$,  (seen as embedded submanifolds). In the particular case of
$R^{2n}$ with the standard symplectic form, an invariant is a function
$c: B(R^{2n}) \rightarrow [0,+\infty]$ such that:
\begin{enumerate}
\item[B1.] Monotonicity: $c(M) \leq c(N)$ if there is a symplectomorphism $\phi$ such
that $\phi(M) \subset N$,
\item[B2.] Conformality: $c(\varepsilon M ) = \varepsilon^{2} c(M)$ for any
$\varepsilon \in R$.
\end{enumerate}

An invariant is nontrivial if it takes finite values on sets with infinite volume,
like cylinders:
$$Z(R) = \left\{ x \in R^{2n} \mbox{ : } x^{2}_{1} + x_{2}^{2} < R \right\}$$

There exist highly nontrivial invariants, as the following theorem shows:

\begin{thm} (Gromov's squeezing theorem) The ball $B(r)$ can be symplectically embedded in the cylinder $Z(R)$   if and only if $r \leq R$.
\end{thm}

This theorem permits to define the invariant:
$$c(A) \ = \ \sup \left\{ R^{2} \mbox{ : } \exists \phi(B(R)) \subset A \right\}$$
called Gromov's capacity.

Another important invariant is Hofer-Zehnder capacity. In order to introduce this
we need the notion of a Hamiltonian flow.

A flow of symplectomorphisms $t \mapsto \phi_{t}$ is Hamiltonian if there is
a function $H: M \times R \rightarrow R$ such that for any time $t$ and place $x$ we have
$$\omega( \dot{\phi}_{t}(x), v) \ = \ dH(\phi_{t}(x),t) v$$
for any $v \in T_{\phi_{t}(x)}M$.

Let $H(R^{2n})$ be the set of compactly supported Hamiltonians. Given
a set $A \subset R^{2n}$, the class of admissible Hamiltonians is $H(A)$, made by all  compactly supported maps in $A$ such that the generated  Hamiltonian flow does not have closed orbits of periods smaller than 1. Then the Hofer-Zehnder capacity
is defined by:
$$hz(A) \ = \ \sup \left\{ \| H \|_{\infty} \mbox{ : } H \in H(A) \right\}$$

Let us denote by $Ham(A)$ the class of Hamiltonian diffeomorphisms compactly supported in $A$. A Hamiltonian diffeomorphism is the time one value of a Hamiltonian flow. In the case which interest us, that is $R^{2n}$, $Ham(A)$ is the connected component of the identity in the group of compactly supported symplectomorphisms.

A curve of Hamiltonian diffeomorphisms (with compact support) is a Hamiltonian flow.
For any such curve $t \mapsto c(t)$ we shall denote by $t\mapsto H_{c}(t, \cdot)$ the associated Hamiltonian  function (with compact support).

On the group of Hamiltonian diffeomorphisms there is a bi-invariant distance introduced by Hofer. This is given by the expression:
\begin{equation}
d_{H}(\phi,\psi) \ = \ \inf \left\{ \int_{0}^{1} \| H_{c}(t) \|_{\infty, R^{2n}} \mbox{ d}t \mbox{ : } c: [0,1] \rightarrow Ham(R^{2n}) \right\}
\label{hodis}
\end{equation}

It is easy to check that $d$ is indeed bi-invariant and it satisfies the triangle property. It is a deep result that $d$ is non-degenerate, that is
$d(id, \phi) \ = \ 0 $ implies $\phi = \ id$.

With the help of the Hofer distance one can define another symplectic invariant, called displacement energy. For a set $A \subset R^{2n}$ the displacement energy is:
$$de(A) \ = \ \inf \left\{ d_{H}(id, \phi) \ \mbox{ : } \phi \in Ham(R^{2n}) \ \ , \ \phi(A) \cap A = \emptyset \right\}$$

A displacement energy can be associated to any group of transformations endowed with a bi-invariant distance (see Eliashberg \& Polterovich \cite{elia}, section 1.3). The fact that the displacement energy is a nontrivial invariant is equivalent with the
non-degeneracy of the Hofer distance. Section 2.  Hofer \& Zehnder \cite{hozen}
is dedicated to the implications of the existence of a non-trivial capacity.    All in all the non-degeneracy of the Hofer distance (proved for example in the Section 5. Hofer \& Zehnder \cite{hozen}) is the cornerstone
of symplectic rigidity theory.

\subsection{Hausdorff measure $2$ and Hofer distance}

We shall give in this section a metric proof of the non-degeneracy of the
Hofer distance.

A flow of symplectomorphisms $t \mapsto \phi_{t}$  with compact support in $A$ is Hamiltonian if it is the projection onto $R^{2n}$ of a flow in
$Hom(H(n),vol, Lip)(A)$.

Consider such a flow which joins identity $id$ with $\phi$. Take the lift
of the flow $t \in [0,1] \mapsto \tilde{\phi}_{t} \in Hom(H(n),vol, Lip)(A)$.

\begin{prop}
The curve $t \mapsto \tilde{\phi}_{t}(x,0)$ has Hausdorff dimension $2$ and measure (up to a multiplicative constant)
\begin{equation}
\mathcal{H}^{2} \left( t \mapsto \tilde{\phi}_{t}(x,0)\right) \ = \ \int_{0}^{1} \mid H_{t}(\phi_{t}(x)) \mid \mbox{ d}t
\label{aham}
\end{equation}
\label{ppham}
\end{prop}

\paragraph{Proof.}
The curve is not horizontal. The tangent has a vertical part (equal to the Hamiltonian). Therefore 
the curve has Hausdorff measure greater than one. We shall compute its Hausdorff measure 2. 

By definition this is: 
$$\mathcal{H}^{2}\left( t \mapsto \tilde{\phi}_{t}(x,0)\right) \ = \ \sup \left\{ Q_{\delta} 
\mbox{ : } \delta > 0 \right\}$$
$$Q_{\delta} \ = \ \inf \left\{ \sum_{i=1}^{n} d^{2}(\tilde{\phi}_{t_{i}}(x,0), \tilde{\phi}_{t_{i+1}}(x,0)) \mbox{ : } 0=t_{1}< ... <t_{N}=1 \ , \ \forall i \ d(\tilde{\phi}_{t_{i}}(x,0), \tilde{\phi}_{t_{i+1}}(x,0)) < \delta \right\}$$
The Ball-Box theorem \ref{bbprop} and the (classical) derivability of the curve ensures us that approximatively 
$$d^{2}(\tilde{\phi}_{t_{i}}(x,0), \tilde{\phi}_{t_{i+1}}(x,0)) \approx \mid H(\phi_{t_{i}}(x))\mid (t_{i+1} - t_{i}) \ + \ \mid \dot{\phi}_{t_{i}}(x) \mid^{2}(t_{i+1}-t_{i})^{2}$$
This implies (almost) the desider inequality \eqref{aham}. If we would use instead of the CC distance 
another distance which is more computable, for example the distance induced by the homogeneous norm 
$$\|(x, \bar{x}) \|^{2} \ = \ \mid x\mid^{2} + \frac{1}{4}\mid \bar{x} \mid$$
then we would have exactly the relation \eqref{aham}, modulus a multiplicative constant $C= \frac{1}{4}$, namely: 
\begin{equation}
\mathcal{H}^{2} \left( t \mapsto \tilde{\phi}_{t}(x,0)\right) \ = \ C \ \int_{0}^{1} \mid H_{t}(\phi_{t}(x)) \mid \mbox{ d}t
\label{ahamc}
\end{equation}
The Hausdorff measure correspond to the chosen homogeneous norm. 

In our case, that is for the Hausdorff measure which comes from the CC distance, the constant $C$ is simply 
$$C \ = \ \lim_{t \rightarrow 0} \frac{d^{2}((0,0), (0,t))}{t}$$
\quad $\blacksquare$

We shall prove now that the Hofer distance \eqref{hodis} is non-degenerate.
For given $\phi$ with compact support in $A$ and generating function
$F$, look at the one parameter family:
$$\tilde{\phi}^{a}(x,\bar{x})  \ = \ (\phi(x), \bar{x} + F(x) + a)$$
The volume of the set 
$$\left\{ (x,z) \mbox{ : } x \in A , \ z \mbox{ between } 0 \mbox{ and }
F(x)+a \right\}$$
attains the minimum   $V(\phi,A)$ for an $a_{0} \in R$.
Let us take  an arbitrary flow $t \in [0,1] \mapsto \tilde{\phi}_{t} \in Hom(H(n),vol, Lip)(A)$ such that $\tilde{\phi}_{0} \ = \ id$ and
$\tilde{\phi}_{1} \ = \ \tilde{\phi}^{a}$.
The family of curves $t \mapsto \tilde{\phi}_{t}(x,0)$ provides a foliation of the set
$$B \ = \ \left\{ \tilde{\phi}_{t}(x,0) \mbox{ : } x \in A \right\}$$

We want to prove the following proposition: 

\begin{prop}
There is a constant $C > 0$ such that 
\begin{equation}
\mathcal{H}^{Q}(B) \ \leq \ C \ vol(A)  \int_{0}^{1} \| H_{t} \|_{\infty,A} \mbox{ d}t
\label{dudu1}
\end{equation}
\label{pdudu1}
\end{prop}

Suppose that \eqref{dudu1} is true.  From the  the definition of the Hofer distance \eqref{hodis} and from the obvious inequality $\mathcal{H}^{Q}(B) \geq V(\phi,A)$ we get 
\begin{equation}
 V(\phi,A) \ \leq \ vol(A) \ C \ d_{H}(id, \phi)
\label{inef1}
\end{equation}
This proves the non-degeneracy of the Hofer distance, because if the
RHS of \eqref{inef1} equals 0 then $V(\phi,A)$ is 0, which means that the generating function of $\phi$ is almost everywhere constant, therefore $\phi$ is the identity everywhere in $R^{2n}$.

\paragraph{Proof of proposition \ref{pdudu1}}

We shall use proposition \ref{ppham} and the definition of Hausdorff measures. For $\varepsilon>0$ pick a $\varepsilon$-dense (in Euclidean distance) set 
$$\left\{ x_{i} \mbox{ : } i \in I_{\varepsilon} \right\}$$
in $A$ and a $\varepsilon$  division $0= t_{0}< ... < t_{N(\varepsilon)}= 1$. It is clear that one can 
(almost)  cover the set $B$ with sets 
$$B_{ij}(\varepsilon) \ = \ L_{\tilde{\phi}_{t_{j}}(x_{i},0)} \left( D\phi_{t_{j}}(x_{i})B^{eucl}_{R^{2n}}(0,\frac{3}{4}\varepsilon) \ \times \ [-\varepsilon,\varepsilon] \right)$$
Denote by $D(\varepsilon)$ this cover and define a $D(\varepsilon)$ nested cover to be  an arbitrary (almost) cover 
$$\mathcal{D}(\varepsilon) \left\{ B^{CC}(\tilde{x}_{\alpha},r_{\alpha}) \mbox{ : } \alpha \in K_{\varepsilon} \right\}$$
 of the set $B$ with CC balls such that any ball $B^{CC}(\tilde{x}_{\alpha},r_{\alpha})$ from the nested cover $\mathcal{D}(\varepsilon)$ is included in exactly one set $B_{ij}$ and moreover 
the center 
$\tilde{x}_{\alpha}$ of the ball  has the form 
$$L_{\tilde{\phi}_{t_{j}}(x_{i},0)}(0,z_{\alpha})$$
We shall denote the fact that $\mathcal{D}(\varepsilon)$ is $D(\varepsilon)$ nested by 
$\mathcal{D}(\varepsilon) < D(\varepsilon)$. 
Let $P(\varepsilon)$ be the number: 
$$P(\varepsilon) \ = \ \inf \left\{ \sum r_{\alpha}^{Q} \mbox{ : } \mathcal{D}(\varepsilon) < D(\varepsilon) \right\} $$
and let 
$$S \ = \ \limsup_{\varepsilon \rightarrow 0} P(\varepsilon)$$
Use proposition \ref{ppham} to estimate $S$, by first "counting" the balls in $\mathcal{D}(\varepsilon) < D(\varepsilon)$ which cover the fiber 
$$\left\{ \tilde{\phi}_{t}(x_{i},0) \mbox{ : } t \in [0,1] \right\}$$
and then "counting" by the elements $x_{i}$ of the $\varepsilon$ net of $A$. We get the upper bound 
$$ C \ vol(A) \ \int_{0}^{1} \| H_{t} \|_{\infty,A} \mbox{ d}t$$
On the other part, not any cover of $B$ with CC balls is nested, therefore from the definition of the Hausdorff measure $\mathcal{H}^{Q}$ we get 
$$S \geq \mathcal{H}^{Q}(B)$$
These two inequalities provide the desired result. \quad $\blacksquare$

We close with the translation of the inequality \eqref{inef1} in symplectic terms.

\begin{prop}
Let $\phi$ be a Hamiltonian diffeomorphism with compact support in
$A$ and $F$ its generating function, that is $dF \ = \ \phi^{*} \lambda -
\lambda$, where $d\lambda \ = \omega$. Consider a Hamiltonian flow
$t \mapsto H_{t}$, with compact support in $A$, such that the time one
map equals $\phi$. Then the following inequality holds:
$$ C \ \inf \left\{ \int_{A} \mid F(x) - c \mid \mbox{ d}x \mbox{ : } c \in
R \right\} \ \leq \ vol(A) \ \int_{0}^{1} \| H_{t}\|_{\infty,A} \mbox{ d}t$$
with $C >0$ an universal constant.

\end{prop}

\subsection{Submanifolds in the Heisenberg group}

We can imagine that a submanifold in the Heisenberg group should be a set
which is locally connected by derivable curves. More than this, suppose that we have a coordinate system made of derivable curves. This means that a submanifold is locally
 a set:
$$\tilde{x} \ = \ \tilde{\phi}(t_{1}, ... , t_{k})$$
such that all curves $$\tau \mapsto \tilde{f}_{j}(\tau) = \tilde{\phi}(... t_{j}+\tau...)$$ (and all other
components fixed) are derivable. Recall that a  curve $\tilde{f}$ is derivable
if it is classically derivable and moreover the quantity
$$\dot{\bar{f}}(t) - \frac{1}{2} \omega(f(t), \dot{f}(t))$$
is constant along the curve. This quantity has an interpretation: it is (up
to a constant) the constant density of the Hausdorff $2$ measure of the curve.

We can imagine therefore (at least) three types of submanifolds. The type I and II  submanifolds are (regular versions of) rectifiable sets:
\begin{enumerate}
\item[I] in the sense of Franchi, Serapioni, Serra Cassano \cite{fraseca}, as
sets which locally are level sets of Lipschitz maps from $H(n)$ to $R$,
\item[II] in the sense of Pauls \cite{pauls1}, as sets which locally are images
of Lipschitz maps from a nilpotent group to $H(n)$.
\end{enumerate}
For the notion of rectifiable sets see Garofalo, Nhieu \cite{ganh}.

We have to comment on type II submanifolds. There is, up to now, no satisfactory
theory of currents in Carnot groups. The Ambrosio, Kirchheim \cite{ambkir} theory of currents in metric spaces does not apply for Carnot groups. It is the belief of the author that a meaningful theory of currents in Carnot groups can be obtained as a theory of representations of the group of "smooth" homeomorphisms of the Carnot group. (This proposition will be the subject of a forthcoming paper; the axioms of this theory, in the case of the Heisenberg group,  have been presented in a talk at the Mathematical Institute,
University of Bern, in 2002. I use the opportunity to thank Martin Reimann for
the kind invitation). It seems reasonable to think that:
\begin{enumerate}
\item[-] in such a theory we should  have, instead of a gradation in the space of currents, a more 
complex relation connected to the ways in in which one can split the Carnot group in smaller Carnot
groups, and
\item[-]  that submanifolds are (associated to) a kind of currents, even in the setting of Carnot groups.
\end{enumerate}
The Ambrosio-Kirchheim theory can be described as a theory associated to $R^{n}$,
which is a Carnot group which split successively in one dimensional subgroups. Hence the usual gradation of currents which lead to the usual classification of manifolds according to their dimension.

In general Carnot group one should consider also the splitting given by the well-
known Kirillov lemma (see, for example \cite{corgre} lemma 1.1.12):

\begin{lema}
Let $\mathfrak{g}$ be a noncommutative nilpotent Lie algebra with one dimensional
center $Z(\mathfrak{g}) = R Z$. Then $\mathfrak{g}$ can written as
$$\mathfrak{g} \ = \ R Z \oplus RY \oplus RX \oplus \mathfrak{w} \ = \ RX \oplus \mathfrak{g}_{0}$$
where $[X,Y] = Z$ and $\mathfrak{g}_{0}$ is the centralizer of $Y$ and an ideal.
\end{lema}
One can split a Carnot group $N$ in two ways: if $Z(N)$ has dimension greater than $1$ then it might be possible to split $N$ using an ideal in $Z(N)$, or, if
$Z(N)$ is one dimensional, we can use Kirillov lemma.

Now I am getting to the point: type II submanifolds in the group $N$ should then be associated to factorisations by  ideals of $N$ which appear in a splitting of
$N$. For example, in the Heisenberg group one can have as type II submanifolds Heisenberg submanifolds or $M$ submanifolds, where $M \subset R^{2n}$ is a Lagrangian subspace (the symplectic form is identically $0$ on $M$).

We shall introduce here yet another type of submanifolds.
\begin{defi}
A type III $k$-submanifold in the Heisenberg group is a set $M \subset H(n)$ such that for
any $\tilde{x} \in M$ there is an open neighbourhood $U_{\tilde{x}} \cap M$ which has the form
$$U_{\tilde{x}} \cap M \ = \ \left\{ \tilde{\phi}_{t}(\tilde{x}) \mbox{ : }
t \in T^{k} \right\} \cap M$$
where $T^{k}$ is the $k$ real torus and $t \in T^{k} \mapsto \tilde{\phi}_{t}$ is the lift (up to multiplication by a constant vertical element of $H(n)$) of a
Hamiltonian action of $T^{k}$ on $R^{2n}$.
\end{defi}

Let $\tilde{x}= (x,\bar{x}) \in M$ and $U_{x}$ be as in the previous definition. Let also $V_{x} \subset R^{2n}$ be the projection of $U_{x}$ on $R^{2n}$.
The moment map of the torus action which defines the manifold around $\tilde{x}$
is denoted by $J: V_{x} \rightarrow R^{k}$. The coordinate curve
$$t_{j} \mapsto \tilde{f}_{j}(t_{j}) = \tilde{\phi}_{(...t_{j} ...)}(x)$$ (all other components are $0$) is then derivable and the derivative of this curve in $t_{j} =
0$ is $$ (\frac{\partial \phi}{\partial t_{j}}(0), J_{j}(x), 0)$$
Therefore a submanifold in the Heisenberg group can be seen as a lift of a
symplectic sub-manifold in $R^{2n}$ endowed with a Hamiltonian torus action. The moment map of the action is visible in the tangent space of the lift.

\begin{defi}
The tangent space at $M$ in  $\tilde{x} \in M$ is the class of all $v \in T_{0}H(n)$ which are tangent to derivable curves in $M$, at $\tilde{x}$. This  space has the description:
$$T_{\tilde{x}}M \ = \ span \left\{(\frac{\partial \phi}{\partial t_{j}}(0), J_{j}(x), 0) \mbox{ : } j = 1,...,k \right\}$$
\end{defi}

A general submanifold should be a combination of these types. The difficult task is to classify these combinations, maybe in the same way irreducible representations of nilpotent Lie groups are classified.

\section{Calculus on general groups}

In this section we come back to a general group $G$ and we try to find the good notion of derivative. The problem that we have is that the Pansu derivative
of a function transforms weak classes of two letter words into similar weak classes. We defined the word tangent bundle to be a set of strong classes. In the case of Carnot groups we had the possibility to choose polynomial representantives
of strong classes, one in each weak class. We cannot do the same in the general case.

Until now we have a Pansu type derivative, which is good for Lipschitz functions.
We also remarked that left, right translations preserve the word tangent bundle, but right translations are not generally Pansu derivable. We need to find a notion of derivative in between these situations.

\subsection{Mild equivalence}

Our starting point is the proposition \ref{pro2}, more specifically the following
useful fact: strong classes act at left on weak classes. This means
the function
$$([c]_{s} , c']_{w} \ \mapsto \ [c]_{s} * [c']_{w} = [c \circ c']_{w}$$
is well defined and it is an action:
$$[c^{1}]_{s}* \left( [c^{2}]_{s}* [c]_{w} \right) \ = \ \left( [c^{1}]_{s} \circ
[c^{2}]_{s} \right) * [c]_{w}$$
Instead of working with strong classes (elements of the word tangent bundle), we shall work with the action of these on weak classes. This is the same idea as in the case of the virtual derivative. If we were in an Euclidean situation, these
choices and differences are not at all significant. Here they are.

We introduced weak and strong equivalence relations. Now we want to identify strong classes which acts the same on weak classes of finite curves $c_{q}$,
$q \in C^{m}_{0}([0,+\infty),G)$, where
$$ C^{m}_{0}([0,+\infty),G) \ = \ \left\{ s \in C^{m}([0,+\infty),G) \mbox{ : } s(0) = 0 \right\} $$
This is the same as introducing yet another equivalence relation on finite curves, which is between the weak and strong relations. We shall call it "mild".

\begin{defi}
Two finite curves $c_{f}, c_{p}$, $f, p \in C^{m}([0,+\infty),G)$,  are mildly
equivalent $c_{f} \mild c_{p}$ if for any $q \in C^{m}_{0}([0,+\infty),G)$,   we have $[c_{f}]_{s} * [c_{q}]_{w} = [c_{p}]_{s} * [c_{q}]_{w}$.

This equivalence relation can be written as $c_{f} \circ c_{q} \weak c_{p} \circ
c_{q}$, for any $q \in C^{m}_{0}([0,+\infty),G)$.
\label{dmild}
\end{defi}

Otherwise said, two elements $P,R \in \mathfrak{g}^{(m)}[\varepsilon]$,   are mildly equivalent if and only if for any $Q \in Weak(\mathfrak{g})$, $Q(0) = 0$ we have
$$\frac{1}{\varepsilon} d(P(\varepsilon) Q(\varepsilon), R(\varepsilon) Q(\varepsilon)) \rightarrow 0$$
as $\varepsilon \rightarrow 0$.

Here comes the second way to define a nice tangent bundle.

\begin{defi}
The mild tangent bundle $MG$ is the group of mild equivalence classes of words, with the bundle structure given by the projection on the first letter.
\label{dmildtan}
\end{defi}

To understand the structure of this bundle we introduce
$$MC(G) \ = \ \left\{ c_{f} \mbox{ : } f \in C^{m}([0,+\infty), G) \right\} / \mild$$
As previously, $MC(G)$ is a Lie group, with the algebra identified with
$$\left(Weak(\mathfrak{g}) + Strong(\mathfrak{g})\right)/M(\mathfrak{g})$$
where $M(\mathfrak{g})$ is the subalgebra of $Strong(\mathfrak{g})$ of elements mildly equivalent with $0$.

\begin{rk}
This algebra is nilpotent (see remarks about the spectrum), like $O(\mathfrak{g})$, but it is not an ideal. Therefore the operation on $MC(G)$ cannot be read directly from the factorisation. The reason for this is that {\bf mild equivalence
is not invariant under conjugation}.
\end{rk}

 The fiber of the neutral element in $MC(G)$ will be denoted by $MC_{0}(G)$. Similarly, $M_{0}G$ denotes the fiber of the neutral element in $MG$.

We can take only the factorisation
$$MC^{0}(G) \ = \ \left\{ c_{f} \mbox{ : } f \in C^{m}_{0}([0,+\infty), G) \right\} / \mild$$
We shall call this the true tangent space of $G$ at the neutral element. This is  also a group and the Lie algebra structure of it is obtained by a factorisation. Indeed, consider the class
$$\mathfrak{g}^{(m)}_{0}[\varepsilon] \ = \ \left\{ P \in \mathfrak{g}^{(m)} \mbox{ : } P(0) =
0 \right\}$$
Then $MC(G) \cap \mathfrak{g}^{(m)}_{0}[\varepsilon]$ is an ideal in
$\mathfrak{g}^{(m)}_{0}[\varepsilon]$ and $MC^{0}(G)$ is a Lie group with Lie algebra obtained by factorisation.

We arrive now to the definition of the  semi-derivative.

\begin{defi}
The function $f$ is semi-derivable at $x$ if there is a function
$$Df(x): \mathfrak{g}^{(m)}_{0}[\varepsilon] \rightarrow
Weak(\mathfrak{g})$$
such that
$$ \frac{1}{\varepsilon} d(f(x P(\varepsilon)), f(x) Df(x)(P)(\varepsilon)) \rightarrow 0$$
as $\varepsilon \rightarrow 0$, uniformly with respect to $P$.
\end{defi}

The  notion of weak  derivative follows.

\begin{defi}
The function $f: G \rightarrow G'$
is  weakly derivable if there is a function  $Df(x): \mathfrak{g}^{(m)}_{0}[\varepsilon] \rightarrow M^{0}G'$ such that for any $Q \in Weak(\mathfrak{g'})$, $Q(0) = 0$, we have
$$\frac{1}{\varepsilon} d(f(x P(\varepsilon)) Q(\varepsilon), Df(x)(P)(\varepsilon) Q(\varepsilon)) \rightarrow 0$$
uniformly with respect to $P$.
\end{defi}

All these algebraic construction seem complicated. There are a lot of differences
(like the one between $MC_{0}(G)$ and $MC^{0}(G)$) which are unseen in the
Euclidean case, and even in the nilpotent case. That is why in the next section we shall take an example and we shall compute everything.

We close this section with the Rademacher theorem. This will be the identical with
the theorem \ref{tradnew}, but from slightly different reasons. Namely, instead  of proposition \ref{pupu} we have:

\begin{prop}
If $f$ satisfies the condition (A) and it is weakly derivable then $Df(x)$ is
a group morphism. Moreover, for any $P,R \in \mathfrak{g}^{(m)}_{0}[\varepsilon]$,
if $P \mild R$ then $Df(x)(P) = Df(x)(Q)$.
\label{popo}
\end{prop}

\paragraph{Proof.} The proof is the same as the one of proposition
\ref{pupu}, with a single difference. The reason for which the quantity $D(\varepsilon)$,  defined by the relation (\ref{clueref}),  converges to $0$ as $\varepsilon \rightarrow 0$ is the derivability of $f$.\quad $\blacksquare$

This suggest the good definition of derivative.

\begin{defi}
The function $f: G \rightarrow G'$
is   derivable if there is a group morphism  $Df(x): MC^{0}(G) \rightarrow MC^{0}(G')$ such that for any $P \in MC^{0}(G)$  we have
$$\frac{1}{\varepsilon} d(f(x P(\varepsilon)) Q(\varepsilon), Df(x)(P)(\varepsilon) Q(\varepsilon)) \rightarrow 0$$
uniformly with respect to $Q \in Weak(\mathfrak{g'})$, $Q(0) = 0$.
\label{dlast}
\end{defi}

The Rademacher theorem in this context is simply this.

\begin{thm}
Let $f: G \rightarrow H$ be a Lipschitz function which satisfies condition (A).
Then $f$ is  a.e. derivable in the sense of definition \ref{dlast}.
\label{trnew}
\end{thm}

If $f: G \rightarrow G'$ is derivable then $Df$ is not a function from
$MC(G)$ to $MC(G')$. This poses a problem concerning upper order derivatives.

The chain rule for this type of derivatives holds. Moreover, if $f$ preserves
the word tangent bundles and $g$ is just derivable then $f\circ g$ is derivable.
We leave to the reader to check this fact.

All in all, we have obtained reasonable notions for calculus in general groups.
But the real hard work of understanding this calculus is still ahead. Several phenomena which are non-Euclidean have been discovered. We have to learn how to control better these new behaviours and this can be done only after deep study
of many particular cases, in order to make the difference between important and
marginal features of non-Euclidean analysis. Only after this we are free to forget
this intricate mechanism and use it with same ease as we do with the classical
calculus.

\subsection{An example: $SO(3)$}

In this section we shall look to the group $G = SO(3)$ of orthogonal matrices in $R^{3}$. The Lie algebra $\mathfrak{so}(3)$ is spanned by the matrices
$$X_{1} \ = \ e_{1} \otimes e_{2} - e_{2} \otimes e_{1} \ \ \ \
 X_{1} \ = \ e_{2} \otimes e_{3} - e_{3} \otimes e_{2} \ \ \ \
X_{1} \ = \ e_{1} \otimes e_{3} - e_{3} \otimes e_{1}$$
with the bracket relations
$$[X_{1}, X_{2}] \ = X_{3} \ \ \ \ [X_{2}, X_{3}] \ = X_{1} \ \ \ \ [X_{3}, X_{1}] \ = X_{2}$$
The space spanned by $X_{1}, X_{2}$ generate the algebra. We shall consider therefore the dilatation flow $\delta_{\varepsilon}: \mathfrak{so}(3) \rightarrow
\mathfrak{so}(3)$
$$\delta_{\varepsilon} X_{1} \ = \ \varepsilon X_{1} \ \ \ \
\delta_{\varepsilon} X_{2} \ = \ \varepsilon X_{2} \ \ \ \
\delta_{\varepsilon} X_{3} \ = \ \varepsilon^{2} X_{3}$$
Then $\mathfrak{so}(3)$ decomposes as
$$\mathfrak{so}(3) \ = \ V_{1} + V_{2} \ \ \ \ V_{1} \ = \ span \left\{ X_{1}, X_{2} \right\} \ \ \ V_{2} \ = \ R X_{3}$$
The nilpotentisation $N$ of $(SO(3), V_{1})$ is the Heisenberg group. The step $m$ of the group is equal to $2$. The algebra $\mathfrak{so}(3)^{(2)}$ is 12 dimensional and it is spanned by $X_{i}$, $\varepsilon X_{i}$, $\varepsilon^{2} X_{i}$, $i= 1,2,3$,  with bracket relations modulo $o(\varepsilon^{2})$. This is the algebra of
the word tangent bundle (check that $WSO(3) = CSO(3)$). It contains as subalgebras
$$\mathfrak{so}(3) \ = \ span \left\{ X_{1}, X_{2}, X_{3} \right\}$$
$$H(1) \  = \ N(SO(3)) \ = \ span\left\{ \varepsilon X_{1}, \varepsilon X_{2},
\varepsilon^{2} X_{3} \right\}$$
(Check that the space $O(\mathfrak{so}(3))$ is trivial.)

The tangent at $0$, called $T_{0}SO(3)$, is 6 dimensional, spanned by
$\varepsilon X_{i}$, $\varepsilon^{2} X_{i}$, $i= 1,2,3$.

If we want to look to $MC^{0}SO(3)$, we find that it is 4 dimensional, isomorphic to the algebra $T_{0}H(1)$. This can be done by knowing that for any $q \in SO(3)$ $Ad_{q}$ is an orthogonal matrix, obtained as a permutation of $q$. This allows us to compute $M(\mathfrak{so}(3))$ to be the group with Lie algebra spanned by $\varepsilon^{2}
X_{1}$, $\varepsilon^{2} X_{2}$, which is abelian (hence Carnot, as forecasted).

Finally, the group $MC(SO(3))$, the mild tangent bundle, can be recovered by looking to the actions of $Ad^{G}$, $Ad^{N}$, on $Weak_{0}(SO(3))$,
 spanned by $\varepsilon X_{i}$, $\varepsilon^{2} X_{3}$, $i= 1,2,3$. Attention, this action is done by applying $Ad^{G}$ or $Ad^{N}$ and then taking projection on
$Weak_{0}(SO(3))$. We notice that $\varepsilon^{2}X_{3}$ is a fixed point of these actions and that the group generated  by $Ad^{G}$, $Ad^{N}$ is simply $SL(3,R)$, which is 8 dimensional. We infer that the mild tangent bundle is $SL(3,R)$.

\end{document}